\documentclass[11pt]{article}
\usepackage{amsmath}
\usepackage{amssymb}
\usepackage{hyperref}

\usepackage{theorem}
\numberwithin{equation}{section}

\newtheorem{theorem}{Theorem}[section]
\newtheorem{proposition}[theorem]{Proposition}
\newtheorem{corollary}[theorem]{Corollary}
\newtheorem{conjecture}[theorem]{Conjecture}
\newtheorem{lemma}[theorem]{Lemma}
\newtheorem{definition}{Definition}[section]

\begin{document}
\title{Global well-posedness and scattering for the defocusing, $L^{2}$-critical, nonlinear Schr{\"o}dinger equation when $d = 2$}
\date{\today}
\author{Benjamin Dodson}
\maketitle

\noindent \textbf{Abstract:} In this paper we prove that the defocusing, cubic nonlinear Schr{\"o}dinger initial value problem is globally well-posed and scattering for $u_{0} \in L^{2}(\mathbf{R}^{2})$. The proof uses the bilinear estimates of \cite{PV} and a frequency localized interaction Morawetz estimate similar to the high frequency estimate of \cite{CKSTT4} and especially the low frequency estimate of \cite{D2}.

\section{Introduction}
In this paper we study the cubic, defocusing, two dimensional initial value problem,

\begin{equation}\label{0.1}
\aligned
i u_{t} + \Delta u &= F(u) = |u|^{2} u, \\
u(0,x) &= u_{0}.
\endaligned
\end{equation}

\noindent Observe that a solution to $(\ref{0.1})$ actually gives a family of solutions. First, if $u(t,x)$ solves $(\ref{0.1})$ with initial data $u_{0}$, then $e^{i \theta} u(t,x)$ solves $(\ref{0.1})$ with initial data $e^{i \theta} u_{0}(x)$. Next, solutions to $(\ref{0.1})$ can be translated in space. If $u(t,x)$ solves $(\ref{0.1})$ with initial data $u_{0}$, then $u(t,x - x_{0})$ solves $(\ref{0.1})$ with initial data $u_{0}(x - x_{0})$.\vspace{5mm}

\noindent It is also possible to translate solutions to $(\ref{0.1})$ in frequency.

\begin{theorem}[Galilean transformation]\label{t1.6}
Suppose $u(t,x)$ solves

\begin{equation}\label{1.21}
\aligned
i u_{t} + \Delta u &= F(u), \\
u(0,x) &= u_{0}.
\endaligned
\end{equation}

\noindent Then $v(t,x) = e^{-it |\xi_{0}|^{2}} e^{ix \cdot \xi_{0}} u(t, x - 2 \xi_{0} t)$ solves the initial value problem

\begin{equation}\label{1.22}
\aligned
i v_{t} + \Delta v &= F(v), \\
v(0,x) &= e^{ix \cdot \xi_{0}} u(0,x).
\endaligned
\end{equation}
\end{theorem}

\noindent \emph{Proof:} This follows by direct calculation. $\Box$\vspace{5mm}

\noindent Finally, for any $\lambda > 0$, if $u(t,x)$ is a solution to $(\ref{0.1})$ on some interval $I$, then

\begin{equation}\label{0.2}
\lambda u(\lambda^{2} t, \lambda x)
\end{equation}

\noindent is a solution to $(\ref{0.1})$ on $\frac{I}{\lambda^{2}} = \{ \frac{t}{\lambda^{2}} : t \in I \}$ with initial data $\lambda u(0, \lambda x)$. $(\ref{0.1})$ is called $L^{2}$ - critical since

\begin{equation}\label{0.3}
 \| \lambda u(0, \lambda x) \|_{L^{2}(\mathbf{R}^{2})} = \| u(0, x) \|_{L^{2}(\mathbf{R}^{2})}.
\end{equation}

\noindent More generally,

\begin{definition}[$\dot{H}^{s_{c}}$ - critical]\label{d0.1}
The defocusing, initial value problem

\begin{equation}\label{0.4}
\aligned
i u_{t} + \Delta u &= F(u) = |u|^{p} u, \\
u(0,x) &= u_{0},
\endaligned
\end{equation}

\noindent is called $\dot{H}^{s_{c}}$ - critical if $p = \frac{4}{d - 2 s_{c}}$, where $u : I \times \mathbf{R}^{d} \rightarrow \mathbf{C}$, $0 \in I \subset \mathbf{R}$.
\end{definition}

\noindent \textbf{Remark:} Notice that the Galilean transformation only preserves the $L^{2}$ norm of $u$, and not the $\dot{H}^{s_{c}}$ norm when $s_{c} > 0$.\vspace{5mm}

\noindent This introduces some additional complications to the $L^{2}$ - critical problem. Criticality plays an important role in the question, ``for which Sobolev spaces is an initial value problem well - posed, either locally or globally?" $L^{2}$ - critical problems with initial data in $L^{2}$ has been an area of intense interest, and $(\ref{0.1})$ in particular.

\begin{conjecture}\label{c0.2}
The initial value problem $(\ref{0.4})$ with $p = \frac{4}{d}$ is globally well - posed and scattering for any $u_{0} \in L^{2}(\mathbf{R}^{d})$.
\end{conjecture}

\noindent \textbf{Remark:} \cite{ChrColTao1} and \cite{ChrColTao2} proved that this conjecture is sharp.\vspace{5mm}

\noindent \cite{TVZ2} proved conjecture $\ref{c0.2}$ is true for dimensions $d \geq 3$ with radially symmetric initial data and \cite{KTV} proved conjecture $\ref{c0.2}$ is true for $u_{0}$ radial, $d = 2$. \cite{D2} then removed the radial symmetry restriction for $d \geq 3$. In this paper we remove the radial symmetry requirement when $d = 2$ and prove

\begin{theorem}\label{t0.3}
The initial value problem $(\ref{0.1})$ is globally well - posed and scattering for any $u_{0} \in L^{2}(\mathbf{R}^{2})$.
\end{theorem}

\noindent This paper is the second in a series of four papers. In an upcoming paper (see preprint \cite{D3}) we prove that conjecture $\ref{c0.2}$ is true when $d = 1$.\vspace{5mm}

\noindent There is also the focusing initial value problem (in this case $F(u) = -|u|^{\frac{4}{d}} u$). There are known counterexamples to global well - posedness and scattering for this problem (see \cite{Glassey}, \cite{Merle}, \cite{TerryTao2}, \cite{We2}, \cite{ZK}). However, all such counterexamples have $L^{2}$ norm above the $L^{2}$ norm of the ground state. \cite{KTV} and \cite{KVZ} showed that conjecture $\ref{c0.2}$ is true when the $L^{2}$ norm of the initial data is below the $L^{2}$ norm of the ground state for $d \geq 2$ and $u_{0}$ radial. \cite{D5} removes the symmetry condition and show that conjecture $\ref{c0.2}$ holds for any data whose $L^{2}$ norm is below the $L^{2}$ norm of the ground state, for any $d \geq 1$.\vspace{5mm}

\noindent Throughout this paper the term solution to $(\ref{0.4})$ on $I \subset \mathbf{R}$, $0 \in I$, is a function $u \in L_{t}^{\infty} L_{x}^{2}(I \times \mathbf{R}^{2})$ that satisfies Duhamel's principle for all $t \in I$.

\begin{definition}[Duhamel's principle]\label{d0.4} $u$ satisfies Duhamel's formula on $I \subset \mathbf{R}$ if $u \in L_{t, loc}^{4}(I; L_{x}^{4}(\mathbf{R}^{2}))$ and for all $t \in I$,



\begin{equation}\label{0.7}
 u(t,x) = e^{it \Delta} u_{0} - i \int_{0}^{t} e^{i(t - \tau) \Delta} F(u(\tau)) d\tau.
\end{equation}



\end{definition}


\begin{definition}[Scattering]\label{d0.5}
A solution to $(\ref{0.4})$ is said to scatter forward in time if there exists $u_{+} \in \dot{H}^{s_{c}}(\mathbf{R}^{d})$ such that

\begin{equation}\label{0.5}
e^{-it \Delta} u(t,x) \rightarrow u_{+},
\end{equation}

\noindent strongly in $\dot{H}^{s_{c}}(\mathbf{R}^{d})$ as $t \rightarrow +\infty$. A solution to $(\ref{0.4})$ is said to scatter backward in time if there exists $u_{-} \in \dot{H}^{s_{c}}(\mathbf{R}^{d})$ such that

\begin{equation}\label{0.6}
e^{-it \Delta} u(t,x) \rightarrow u_{-},
\end{equation}

\noindent strongly in $\dot{H}^{s_{c}}(\mathbf{R}^{d})$ as $t \rightarrow -\infty$.
\end{definition}

\noindent A solution to $(\ref{0.4})$ conserves the quantities mass,

\begin{equation}\label{0.5}
M(u(t)) = \int |u(t,x)|^{2} dx = M(u(0)),
\end{equation}

\noindent and energy,

\begin{equation}\label{0.6}
E(u(t)) = \int (\frac{1}{2} |\nabla u(t,x)|^{2} + \frac{1}{p + 2} |u(t,x)|^{p + 2}) dx = E(u(0)).
\end{equation}

\noindent It is for this reason that research on $\dot{H}^{s_{c}}$ - critical problems has generally concentrated on the $L^{2}$ or mass - critical problems ($p = \frac{4}{d}$) and the $\dot{H}^{1}$ or energy - critical problems ($p = \frac{4}{d - 2}$). However see \cite{KM2}, \cite{Murphy1}, and \cite{Murphy2} for global well - posedness and scattering when a solution to $(\ref{0.4})$ has an assumed bound on the $\dot{H}^{s_{c}}$ norm, $0 < s_{c} < 1$.\vspace{5mm}

\noindent The local theory for $(\ref{0.4})$ has been worked out in \cite{Caz}, \cite{Caz1}, \cite{CaWe}, \cite{CaWe2}, \cite{CaWe1}, and \cite{TV}. For simplicity, only the results that pertain to $(\ref{0.1})$ will be presented here.

\begin{definition}[Blowup criterion]\label{d0.6}
A solution $u : I \times \mathbf{R}^{2} \rightarrow \mathbf{C}$ is said to blow up forward in time if for some $t_{0} \in I$,

\begin{equation}\label{0.7}
\int_{t_{0}}^{\sup(I)} \int |u(t,x)|^{4} dx dt = \infty.
\end{equation}

\noindent Blowing up backward in time is similarly defined.
\end{definition}

\begin{theorem}[Local well - posedness]\label{t0.7}

\begin{enumerate}
\item Given $u_{0} \in L^{2}(\mathbf{R}^{2})$ there exists a unique solution $u$ to $(\ref{0.1})$ with initial data $u_{0}$ on some open interval. Then if $I$ is the maximal interval of existence, $I$ is open. Moreover, for any compact $J \subset I$ the solution map $u(0) \mapsto u(t,x) \in L_{t,x}^{4}(J \times \mathbf{R}^{2})$ is also continuous in an $L^{2}$ - neighborhood of $u_{0}$.

\item If $\sup(I)$ is finite then the solution $u$ blows up forward in time. If $\inf(I)$ is finite then the solution $u$ blows up backward in time.

\item If $u$ does not blow up forward in time then $\sup(I) = \infty$ and $u$ scatters forward in time. If $u$ does not blow up backward in time then $\inf(I) = -\infty$ and $u$ scatters backward in time.

\item If $M(u_{0})$ is sufficiently small then the solution $u$ is global and does not blow up either forward or backward in time.
\end{enumerate}
\end{theorem}

\noindent \emph{Proof:} See \cite{CaWe}, \cite{CaWe1}. $\Box$\vspace{5mm}

\noindent \textbf{Remark:} Observe that combining the pseudoconformal conservation law (see \cite{Caz1} and the references therein),

\begin{equation}\label{0.6.1}
\| (x - 2it \nabla) u(t) \|_{L^{2}(\mathbf{R}^{d})}^{2} + \frac{8 t^{2}}{d + 2} \| u(t) \|_{L_{x}^{\frac{2(d + 2)}{d}}(\mathbf{R}^{d})}^{\frac{2(d + 2)}{d}} = \| x u_{0} \|_{L^{2}(\mathbf{R}^{d})}^{2},
\end{equation}

\noindent with conservation of mass $(\ref{0.6})$ and energy $(\ref{0.6.1})$, and the Sobolev embedding theorem,

\begin{equation}\label{0.6.2}
\| u \|_{L_{t,x}^{\frac{2(d + 2)}{d}}(\mathbf{R} \times \mathbf{R}^{d})}^{\frac{2(d + 2)}{d}} \lesssim \| xu_{0} \|_{L^{2}(\mathbf{R}^{d})} \| u_{0} \|_{H^{1}(\mathbf{R}^{d})} \| u_{0} \|_{L^{2}(\mathbf{R}^{d})}^{4/d}.
\end{equation}

\noindent The left and right hand sides of $(\ref{0.6.2})$ are invariant under the scaling $(\ref{0.2})$.\vspace{5mm}

\noindent The proof of theorem $\ref{t0.3}$ uses the concentration compactness method. This method has been utilized since at least the 1980's (\cite{BL}, \cite{BC}) in a wide variety of partial differential equations, including elliptic, hyperbolic, parabolic, and geometric partial differential equations. The concentration compactness method, and other related methods such as induction on energy and concentration in a backward light cone of a wave equation blowup solution, has also proved to be extremely fruitful to the study of dispersive partial differential equations.\vspace{5mm}


\noindent Historically, such progress has gone energy - critical wave, energy - critical Schr{\"o}dinger, and then mass - critical Schr{\"o}dinger, although there are exceptions to this. For the defocusing, semilinear, energy - critical wave equation, \cite{Struwe} proved global well - posedness and scattering in the radial case when $d = 3$, \cite{Grillakis} in the general case when $d = 3$. \cite{GVS}, \cite{Grillakis1}, \cite{ShatahStruwe}, \cite{ShatahStruwe2}, \cite{ShatahStruwe3}, \cite{Kapitanski} generalized this result to higher dimensions. See also \cite{BahSha}, \cite{TerryTao3}. See \cite{KM3} for a treatment of the focusing wave equation problem.\vspace{5mm}


\noindent The defocusing, energy - critical Schr{\"o}dinger problem is also complete. \cite{B4} proved global well - posedness and scattering for the energy - critical problem in three and four dimensions with radial data. See also \cite{B3}. Independently, \cite{Gril} proved global well - posedness of the radial, energy - critical problem in three dimensions. \cite{TerryTao} extended this result to higher dimensions. \cite{CKSTT4}, \cite{RhV}, \cite{Visan}, \cite{V2} then proved global well - posedness and scattering for the defocusing energy - critical initial value problem with nonradial data. See \cite{KV2}, \cite{V2} for alternate proofs using the methods of this paper. \cite{KM1}, \cite{KV1}, and \cite{TaoT1} treated the focusing energy - critical problem.\vspace{5mm}

\noindent \cite{KM2} treated the $\dot{H}^{1/2}$ - critical problem. See also \cite{Murphy1}, \cite{Murphy2}, \cite{TerryTao1} for additional work on inter - critical ($0 < s_{c} < 1$) problems.\vspace{5mm}

\noindent Theorem $\ref{t0.3}$ follows directly from two results:

\begin{theorem}\label{t0.8}
If theorem $\ref{t0.3}$ fails to hold then there exists a nonzero solution to $(\ref{0.1})$ that lies in a compact subset of $L^{2}(\mathbf{R}^{2})$ modulo scaling, translation, and Galilean symmetries for the entire time of its existence.
\end{theorem}

\noindent \textbf{Remark:} The group $S^{1} = \mathbf{R} / 2 \pi \mathbf{Z}$, which maps $u \mapsto e^{i \theta} u$, is a compact group.
\begin{definition}[Almost periodicity]\label{d0.10}
A solution that lies in a compact subset of $L^{2}(\mathbf{R}^{2})$ modulo scaling, translation, and Galilean symmetries for the entire time of its existence is called an almost periodic solution.
\end{definition}

\begin{theorem}\label{t0.9}
The only solution to $(\ref{0.1})$ that is almost periodic for the entire time of its existence is $u \equiv 0$.
\end{theorem}

\noindent Theorem $\ref{t0.8}$ was finally proved by \cite{TVZ1} in all dimensions $d \geq 1$, building on the seminal work of \cite{MerleVega} (which used a Strichartz result of \cite{MVV}), as well as the work of \cite{Keraani}, \cite{TaoT}, and \cite{BegV}. For the energy - critical problem, \cite{Keraani1} proved a profile decomposition for the energy - critical Schr{\"o}dinger problem in $\mathbf{R}^{3}$ in the same vein as the profile decomposition \cite{Gerard} proved for the Sobolev embedding and \cite{BahGer} proved for the wave equation.\vspace{5mm}

\noindent Now let $G$ be the group generated by the Galilean, translation, and scaling symmetries, and for $g_{\xi_{0}, x_{0}, N_{0}} \in G$, let

\begin{equation}
g_{\xi_{0}, x_{0}, N_{0}} u(x) = \frac{1}{N_{0}} e^{ix \cdot \xi_{0}} u(\frac{x - x_{0}}{N_{0}}).
\end{equation}

\noindent By theorem $\ref{t0.8}$, there exists a compact $K \subset L^{2}(\mathbf{R}^{2})$ such that for each $t \in I$, there exists $g(t) \in G$ such that $g(t) u(t) \in K$. Equivalently, by the Arzela - Ascoli theorem there exist $\xi(t) : I \rightarrow \mathbf{R}^{2}$, $x(t) : I \rightarrow \mathbf{R}^{2}$, $N(t) : I \rightarrow (0, \infty)$, $C : (0, \infty) \rightarrow (0, \infty)$ such that for all $\eta > 0$,

\begin{equation}\label{1.19}
\int_{|x - x(t)| \geq \frac{C(\eta)}{N(t)}} |u(t,x)|^{2} dx < \eta,
\end{equation}

\begin{equation}\label{1.20}
\int_{|\xi - \xi(t)| \geq C(\eta) N(t)} |\hat{u}(t,\xi)|^{2} d\xi < \eta.
\end{equation} 

\noindent By the uncertainty principle and the symmetries discussed in theorem $\ref{t1.6}$ and $(\ref{0.2})$, this is as much as one could possibly hope for.\vspace{5mm}

\noindent \textbf{Remark:} Notice that by $(\ref{1.19})$ and $(\ref{1.20})$, we have the freedom to multiply $N(t)$ by a constant, translate $\xi(t)$ by a distance $\lesssim N(t)$, or translate $x(t)$ by a distance $\lesssim \frac{1}{N(t)}$.\vspace{5mm}

\noindent \emph{Sketch of the proof of $\ref{t0.9}$:} The proof of theorem $\ref{t0.9}$ will occupy the rest of the paper. Since $u(t) \in C(I ; L^{2}(\mathbf{R}^{2}))$ lies in a compact subset of $L^{2}(\mathbf{R}^{2}) / G$, where $G$ is the group generated by translation, Galilean, and scaling symmetries,

\begin{lemma}\label{l1.8}
For any nonzero almost periodic solution $u$ to $(\ref{0.1})$ there exists $\delta(u) > 0$ such that for any $t_{0} \in I$,

\begin{equation}\label{1.26}
\| u \|_{L_{t,x}^{4}([t_{0}, t_{0} + \frac{\delta}{N(t_{0})^{2}}] \times \mathbf{R}^{2})} \sim \| u \|_{L_{t,x}^{4}([t_{0} - \frac{\delta}{N(t_{0})^{2}}, t_{0}] \times \mathbf{R}^{2})} \sim 1.
\end{equation}
\end{lemma}

\noindent \emph{Proof:} See lemma $5.18$ of \cite{KilVis}. $\Box$\vspace{5mm}

\noindent Lemma $5.18$ of \cite{KilVis} also proved

\begin{lemma}\label{l1.7}
If $J$ is an interval with

\begin{equation}\label{1.23}
\| u \|_{L_{t,x}^{4}(J \times \mathbf{R}^{2})} = 1,
\end{equation}

\noindent then for $t_{1}, t_{2} \in J$,

\begin{equation}\label{1.24}
N(t_{1}) \sim N(t_{2}),
\end{equation}

\noindent and

\begin{equation}\label{1.25}
|\xi(t_{1}) - \xi(t_{2})| \lesssim \sup_{t \in J} N(t).
\end{equation}
\end{lemma}

\noindent Combining lemmas $\ref{l1.8}$ and $\ref{l1.7}$, we can choose $\xi(t)$, $N(t)$ such that

\begin{equation}\label{1.29}
|\xi'(t)| + |N'(t)| \lesssim N(t)^{3}.
\end{equation}

\noindent It is convenient to use the notation

\begin{definition}
If $J$ is an interval then let

\begin{equation}\label{1.27}
 N(J) = \sup_{t \in J} N(t).
\end{equation}
\end{definition}

\noindent By lemmas $\ref{l1.7}$ and $\ref{l1.8}$, if $\| u \|_{L_{t,x}^{4}(J \times \mathbf{R}^{2})} = 1$,

\begin{equation}\label{1.28}
N(J) \sim \int_{J} N(t)^{3} dt \sim \inf_{t \in J} N(t).
\end{equation}

\noindent Next, since $u(t)$ lies in a compact subset of $L^{2} / G$ for all $t \in I$, where $I$ is the maximal interval of existence of $(\ref{0.1})$, we can utilize the perturbation theory of \cite{TVZ} and take limits of subsequences of $u(t_{n})$, $t_{n} \in I$.

\begin{theorem}\label{t1.9}
If there exists an almost periodic solution to $(\ref{0.1})$ with $0 < \epsilon \leq \| u(0) \|_{L^{2}} < \infty$, then there exists an almost periodic solution to $(\ref{0.1})$ satisfying $(\ref{1.19})$ and $(\ref{1.20})$, $0 < \epsilon \leq \| u(0) \|_{L^{2}(\mathbf{R}^{2})} < \infty$, 

\begin{equation}\label{1.30}
 \int_{0}^{\infty} \int |u(t,x)|^{4} dx dt = \infty,
\end{equation}

\noindent $N(0) = 1$, $\xi(0) = x(0) = 0$, and $N(t) \leq 1$ on $[0, \infty)$, $|N'(t)| + |\xi'(t)| \lesssim_{u} N(t)^{3}$.
\end{theorem}

\noindent \emph{Proof:} This was proved in \cite{TVZ2} and \cite{TVZ1}. Since $u(t)$ lies in a compact set in $L^{2}(\mathbf{R}^{2})$ modulo symmetries, combining the fact that $N(t)$ is continuous with time reversal symmetry, we can take a limit of $u(t_{n})$, $t_{n} \in I$ under the various symmetries, and obtain $N(t) \leq 1$ for $t \geq 0$, $N(0) = 1$. By the perturbation theory of \cite{TVZ} the limit of $u(t_{n})$ modulo symmetries will be the initial data of an almost periodic solution $u$ with $x(0) = \xi(0) = 0$, $N(0) = 1$. Then by lemma $\ref{l1.8}$, $u$ is defined on $[0, \infty)$. Finally by $(\ref{1.29})$ the theorem holds. $\Box$\vspace{5mm}

\noindent \textbf{Remark:} Lemmas $\ref{l1.8}$ and $\ref{l1.7}$ also imply that if $u$ is a nonzero, almost periodic solution to $(\ref{0.1})$ on a maximal interval $I \subset \mathbf{R}$,

\begin{equation}\label{1.26.1}
\int_{0}^{\sup(I)} \int |u(t,x)|^{4} dx dt = \int_{\inf(I)}^{0} \int |u(t,x)|^{4} dx dt = \infty.
\end{equation}

\noindent Therefore, an almost periodic solution blows up forward and backward in time, which implies that a nonzero almost periodic solution must have mass greater than the small data threshold of theorem $\ref{t0.7}$. \vspace{5mm}

\noindent At this point, we will now take an almost periodic solution satisfying theorem $\ref{t1.9}$, and prove that such a solution must satisfy $u \equiv 0$. For the rest of the paper any expression of the form $A \lesssim_{u} B$ will be abbreviated $A \lesssim B$.\vspace{5mm}

\noindent Proving $u \equiv 0$ follows an argument similar to the argument in \cite{D2}, namely it suffices to exclude $N(t) \leq 1$ on $[0, \infty)$, where $N(t)$ is the scale parameter. There are two possible scenarios, the rapid frequency cascade scenario,

\begin{equation}\label{0.8}
 \int_{0}^{\infty} N(t)^{3} dt < \infty,
\end{equation}

\noindent and the quasisoliton,

\begin{equation}\label{0.9}
\int_{0}^{\infty} N(t)^{3} dt = \infty.
\end{equation}

\noindent The main new ingredient of \cite{D2} was a long time Strichartz estimate, which was used to exclude both $(\ref{0.8})$ and $(\ref{0.9})$. This long time Strichartz estimate relied heavily on the endpoint $L_{t}^{2} L_{x}^{\frac{2d}{d - 2}}$ Strichartz estimate of \cite{KT}, which holds for $d \geq 3$. However, when $d = 2$ the $L_{t}^{2} L_{x}^{\infty}$ estimate does not hold (see \cite{MS}). This makes the long time Strichartz estimates considerably more difficult to both define and to prove. In section three we construct our function spaces $\tilde{X}_{k_{0}}$ out of the $U_{\Delta}^{2}$ spaces of \cite{KoTa}.\vspace{5mm}

\noindent Then, in section four, we will prove theorem $\ref{t3.1}$, which bounds the long time Strichartz norms on our almost periodic solution. As in \cite{D2}, the idea is that $(\ref{1.19})$ and $(\ref{1.20})$ imply that an almost periodic solution to $(\ref{0.1})$ must be mostly concentrated in balls of radius $\sim N(t)$ in frequency and $\sim \frac{1}{N(t)}$ in space. Outside this frequency ball, the solution to $(\ref{0.1})$ will be dominated by solutions to $(i \partial_{t} + \Delta) u = 0$ for long periods of time.\vspace{5mm}

\noindent Although the intuition is also true here, the analysis is much more technically complicated due to a lack of endpoint Strichartz estimates. Indeed, the proof of theorem $\ref{t3.1}$ occupies all of section four, by far the longest section of the paper, as well as the appendix in section six.\vspace{5mm}

\noindent The proof utilizes bilinear Strichartz estimates, both the bilinear Strichartz estimates of \cite{B2} (proposition $\ref{p1.2}$) and \cite{TaoT} (proposition $\ref{p1.4}$), as well as theorems $\ref{t3.2.0}$, $\ref{t3.2.1}$, and $\ref{t3.2.2}$. Theorems $\ref{t3.2.0}$, $\ref{t3.2.1}$, and $\ref{t3.2.2}$ make use of the interaction Morawetz estimates of \cite{PV}.\vspace{5mm}

\noindent Finally, in section five we complete the proof of theorem $\ref{t0.9}$, proving a rigidity result that if $u$ is an almost periodic solution in the form of theorem $\ref{t1.9}$, then $u \equiv 0$. First, we make use of the long time Strichartz estimates and show that if $u$ is an almost periodic solution to $(\ref{0.1})$ that is a rapid frequency cascade, then $u \in L_{t}^{\infty} \dot{H}^{3}$. Then by conservation of energy, this implies that $u \equiv 0$.\vspace{5mm}

\noindent Next, we show that if $u$ is a quasi - soliton, then $u \equiv 0$. We show this with a frequency localized interaction Morawetz estimates. Morawetz estimates have long been useful to proving scattering results for dispersive equations (see \cite{B6}, \cite{GV1}, \cite{LinStrauss}, \cite{Morawetz}, \cite{Nakanishi}), particularly for radial data. For nonradial data the interaction Morawetz estimate has also proved to be quite useful.

\begin{theorem}[Interaction Morawetz estimate]\label{t0.11}
A solution to $(\ref{0.4})$ has the bounds

\begin{equation}\label{0.10}
\| |\nabla|^{\frac{3 - d}{2}} |u(t,x)|^{2} \|_{L_{t,x}^{2}(I \times \mathbf{R}^{d})}^{2} \lesssim \| u \|_{L_{t}^{\infty} L_{x}^{2}(I \times \mathbf{R}^{d})}^{2} \| u \|_{L_{t}^{\infty} \dot{H}^{1/2}(I \times \mathbf{R}^{d})}^{2}
\end{equation}

\noindent for all $d \geq 1$.
\end{theorem}

\noindent \emph{Proof:} See \cite{CKSTT2} when $d = 3$, \cite{TVZ2} for $d \geq 4$, and \cite{CGT1}, \cite{PV} for $d = 1, 2$. \cite{PV} improves $(\ref{0.10})$ to a Galilean invariant version. Galilean invariance of $(\ref{0.10})$ will be utilized heavily in section five. $\Box$\vspace{5mm}

\noindent Since there is no a priori bound on $\| u(t) \|_{\dot{H}^{1/2}}$, we truncate to low frequency. In the energy - critical problem, see \cite{CKSTT4} for example, a solution is truncated to high frequencies and then the interaction Morawetz estimate is computed. Here, as in \cite{D2}, the errors produced by truncating to low frequenies are successfully estimated using the long time Strichartz estimates. Morawetz estimates truncated to low frequencies are very closely related to the almost Morawetz estimates used in the I - method. (See \cite{CGT}, \cite{CR}, \cite{D7} for the almost Morawetz estimates in conjunction with the I - method, \cite{B2} for the Fourier truncation method, and \cite{CKSTT1} for the I - method).\vspace{5mm}

\noindent Finally, in section six we prove lemma $\ref{appendix}$, an unpaid debt from section four.

\section{Linear Estimates}
This section serves to present some linear and bilinear Strichartz estimates that will be used in subsequent sections.

\begin{definition}[Admissible pair]\label{d1.1}
A pair $(p,q)$ will be called an admissible pair for $d = 2$ if $\frac{2}{p} = 2(\frac{1}{2} - \frac{1}{q})$ and $p > 2$.
\end{definition}

\begin{theorem}\label{t1.1}
If $u(t,x)$ solves the initial value problem

\begin{equation}\label{1.1}
\aligned
i u_{t} + \Delta u &= F(t), \\
u(0,x) &= u_{0},
\endaligned
\end{equation}

\noindent on an interval $I$, then

\begin{equation}\label{1.2}
\| u \|_{L_{t}^{p} L_{x}^{q}(I \times \mathbf{R}^{2})} \lesssim_{p,q,\tilde{p},\tilde{q}} \| u_{0} \|_{L^{2}(\mathbf{R}^{2})} + \| F \|_{L_{t}^{\tilde{p}'} L_{x}^{\tilde{q}'}(I \times \mathbf{R}^{2})},
\end{equation}

\noindent for all admissible pairs $(p,q)$, $(\tilde{p}, \tilde{q})$. $\tilde{p}'$ denotes the Lebesgue dual of $\tilde{p}$.
\end{theorem}

\noindent \emph{Proof:} See \cite{Stri} for $p = q = 4$, \cite{ChrKis}, \cite{GV}, \cite{Tao}, \cite{Yaj} for the general result. $\Box$\vspace{5mm}

\noindent \textbf{Remark:} Some endpoint results are available for radial data, (see \cite{Stef} and \cite{TaoT2}), however, these results will not be discussed here since this paper is concerned with nonradial results.

\begin{proposition}\label{p1.2}
If $\hat{u}_{0}$ is supported on $|\xi| \sim N$, $\hat{v}_{0}$ is supported on $|\xi| \sim M$, $M << N$,

\begin{equation}\label{1.3}
\| (e^{it \Delta} u_{0})(e^{it \Delta} v_{0}) \|_{L_{t,x}^{2}(\mathbf{R} \times \mathbf{R}^{2})} \lesssim (\frac{M}{N})^{1/2} \| u_{0} \|_{L^{2}(\mathbf{R}^{2})} \| v_{0} \|_{L^{2}(\mathbf{R}^{2})}.
\end{equation}


\end{proposition}

\noindent \emph{Proof:} See \cite{B2}. $\Box$

\begin{corollary}\label{c1.3}
\begin{equation}\label{1.5}
\| (e^{it \Delta} u_{0})(e^{it \Delta} v_{0}) \|_{L_{t}^{p} L_{x}^{q}(\mathbf{R} \times \mathbf{R}^{2})} \lesssim (\frac{M}{N})^{1/p} \| u_{0} \|_{L^{2}(\mathbf{R}^{2})} \| v_{0} \|_{L^{2}(\mathbf{R}^{2})},
\end{equation}

\noindent for $\frac{1}{p} + \frac{1}{q} = 1$, $2 \leq p \leq \infty$. 


\end{corollary}

\noindent \emph{Proof:} Interpolate the elementary $L_{t}^{\infty} L_{x}^{1}$ bilinear estimate from theorem $\ref{t1.1}$ with proposition $\ref{p1.2}$. $\Box$\vspace{5mm}

\noindent Proposition $\ref{p1.2}$ is $L^{2}$ - critical since it is invariant under scaling symmetry $(\ref{0.2})$. The subcritical bilinear Strichartz estimate of \cite{TaoT} is also an extremely useful result.

\begin{proposition}\label{p1.4}
Suppose $q > \frac{d + 3}{d + 1}$. If the Fourier supports of $\hat{u}_{0}$ and $\hat{v}_{0}$ are separated by distance $\geq cN$ and $\hat{u}_{0}(\xi)$ and $\hat{v}_{0}(\xi)$ are supported in $|\xi| \leq N$, then

\begin{equation}\label{1.8}
\| (e^{it \Delta} u_{0})(e^{it \Delta} v_{0}) \|_{L_{t,x}^{q}(\mathbf{R} \times \mathbf{R}^{d})} \leq C(c) N^{d - \frac{d + 2}{q}} \| u_{0} \|_{L^{2}(\mathbf{R}^{d})} \| v_{0} \|_{L^{2}(\mathbf{R}^{d})}.
\end{equation}
\end{proposition}

\noindent \emph{Proof:} See \cite{TaoT}. $\Box$\vspace{5mm}

\noindent \textbf{Remark:} Propositions $\ref{p1.2}$ and $\ref{p1.4}$ also hold for $(e^{it \Delta} u_{0})(\overline{e^{it \Delta} v_{0}})$. Indeed, since $|e^{it \Delta} u_{0}|^{2} = (e^{it \Delta} u_{0})(\overline{e^{it \Delta} u_{0}})$,

\begin{equation}\label{1.9}
\aligned
 \| (e^{it \Delta} u_{0}) (\overline{e^{it \Delta} v_{0}}) \|_{L_{t,x}^{q}(\mathbf{R} \times \mathbf{R}^{d})}^{q} = \int_{\mathbf{R}} \int_{\mathbf{R}^{d}} |e^{it \Delta} u_{0}|^{q} |e^{it \Delta} v_{0}|^{q} dx dt = \| (e^{it \Delta} u_{0})(e^{it \Delta} v_{0}) \|_{L_{t,x}^{q}(\mathbf{R} \times \mathbf{R}^{d})}^{q}.
\endaligned
\end{equation}

\noindent Interpolating $(\ref{1.8})$ with the elementary $L_{t}^{\infty} L_{x}^{1}$ bilinear estimate gives

\begin{equation}\label{1.10}
\| (e^{it \Delta} u_{0})(e^{it \Delta} v_{0}) \|_{L_{t}^{2} L_{x}^{3/2+}(\mathbf{R} \times \mathbf{R}^{2})} \lesssim N^{-1/3+} \| u_{0} \|_{L^{2}(\mathbf{R}^{2})} \| v_{0} \|_{L^{2}(\mathbf{R}^{2})}.
\end{equation}

\noindent When $d = 1$ the endpoint of $(\ref{1.8})$ does hold,

\begin{equation}\label{1.11}
 \| (e^{it \Delta} u_{0})(e^{it \Delta} v_{0}) \|_{L_{t,x}^{2}(\mathbf{R} \times \mathbf{R})} \lesssim N^{-1/2} \| u_{0} \|_{L^{2}(\mathbf{R})} \| v_{0} \|_{L^{2}(\mathbf{R})}.
\end{equation}

\noindent This gives an improvement over proposition $\ref{p1.2}$ in certain situations.

\begin{proposition}\label{p1.5}
Let $\xi \in \mathbf{R}^{2}$ be the Fourier variable, $\xi = (\xi_{1}, \xi_{2})$. Also suppose that the $\xi_{1}$ support of $\hat{u}_{0}(\xi)$ and $\hat{v}_{0}(\xi)$ are separated by distance $\sim N$. Then by H{\"o}lder's inequality, if $P$ is a projection onto a strip $|\xi_{2}| \leq M$,

\begin{equation}\label{1.12}
\| P((e^{it \Delta} u_{0})(e^{it \Delta} v_{0})) \|_{L_{t,x}^{2}(\mathbf{R} \times \mathbf{R}^{2})} \lesssim \frac{M^{1/2}}{N^{1/2}} \| u_{0} \|_{L^{2}(\mathbf{R}^{d})} \| v_{0} \|_{L^{2}(\mathbf{R}^{2})}.
\end{equation}
\end{proposition}

\noindent A similar result is proved in \cite{CKSTT3}. $\Box$\vspace{5mm} 

\noindent The linear and bilinear Strichartz estimates also hold under convolutions with $L^{1}$ kernels, an important fact since a nonradial solution will be localized in frequency around some $\xi(t) \in \mathbf{R}^{2}$, where $\xi(t)$ is free to move around.\vspace{5mm}

\noindent Suppose $g(t,x - y)$ and $h(t, x - z)$ are convolution kernels with the bounds 

\begin{equation}\label{1.12.1}
\| \sup_{t \in \mathbf{R}} |g(t, x)| \|_{L^{1}(\mathbf{R}^{2})}, \hspace{5mm} \text{and} \hspace{5mm} \| \sup_{t \in \mathbf{R}} |h(t,x)| \|_{L^{1}(\mathbf{R}^{2})} \lesssim 1.
\end{equation}

\noindent Then

\begin{equation}\label{1.13}
 \| (g(t, \cdot) \ast e^{it \Delta} u_{0})(h(t, \cdot) \ast e^{it \Delta} v_{0}) \|_{L_{t,x}^{q}(\mathbf{R} \times \mathbf{R}^{d})}
\end{equation}

\begin{equation}\label{1.14}
\aligned
= \| \int \int g(t, x - y) h(t, x - z) (e^{it \Delta} u_{0})(y) (e^{it \Delta} v_{0})(z) dy dz \|_{L_{t,x}^{q}(\mathbf{R} \times \mathbf{R}^{d})} \\
= \| \int \int g(t, y) h(t, z) (e^{it \Delta} u_{0})(-y + x) (e^{it \Delta} v_{0})(-z + x) dy dz \|_{L_{t,x}^{q}(\mathbf{R} \times \mathbf{R}^{d})} \\
 \lesssim \sup_{y, z} \| (e^{it \Delta} u_{0})(-y + x) (e^{it \Delta} v_{0})(-z + x) \|_{L_{t,x}^{q}(\mathbf{R} \times \mathbf{R}^{d})}.
\endaligned
\end{equation}

\noindent Since theorem $\ref{t1.1}$ and propositions $\ref{p1.2}$ - $\ref{p1.4}$ hold under translations of the initial data, if $\hat{u}_{0}$ is supported on $|\xi| \sim M$, $\hat{v}_{0}$ is supported on $|\xi| \sim N$,

\begin{equation}\label{1.14.1}
\| (g \ast e^{it \Delta} u_{0})(h \ast e^{it \Delta} v_{0}) \|_{L_{t}^{p} L_{x}^{q}} \lesssim c(M, N) \| u_{0} \|_{L^{2}} \| v_{0} \|_{L^{2}},
\end{equation}

\noindent where $c(M, N)$ is the bilinear constant found in theorem $\ref{t1.1}$, propositions $\ref{p1.2}$ - $\ref{p1.4}$. The kernels of $P_{\xi(t), j}$, $P_{\xi(t), \leq j}$, and $P_{\xi(t), \geq j}$ (see $(\ref{2.6})$) all satisfy $(\ref{1.12.1})$.

\section{A function space adapted to the long time Strichartz estimates}
\noindent The reader should compare the function spaces introduced in this section to the long time Strichartz estimates of \cite{D2}. In both cases, the estimates are based on the idea that when $|\xi - \xi(t)| >> N(t)$, an almost periodic solution to $(\ref{0.1})$ is dominated by a solution to the linear problem $i v_{t} + \Delta v = 0$ for long periods of time. The crucial difference is that the endpoint norms $L_{t}^{2} L_{x}^{\frac{2d}{d - 2}}$ are replaced by the $U_{\Delta}^{2}$ spaces.

\begin{definition}[Littlewood - Paley decomposition]\label{d2.1}
Let $\phi \in C_{0}^{\infty}(\mathbf{R}^{2})$ be a radial, decreasing function,

\begin{equation}\label{2.1}
\phi(x) = \left\{
            \begin{array}{ll}
              1, & \hbox{$|x| \leq 1$;} \\
              0, & \hbox{$|x| > 2$.}
            \end{array}
          \right.
\end{equation}

\noindent Define the partition of unity

\begin{equation}\label{2.2}
1 = \phi(x) + \sum_{j = 1}^{\infty} [\phi(2^{-j} x) - \phi(2^{-j + 1} x)] = \psi_{0}(x) + \sum_{j = 1}^{\infty} \psi_{j}(x).
\end{equation}

\noindent For any integer $j \geq 0$, let

\begin{equation}\label{2.3}
P_{j} f = \mathcal F^{-1}(\psi_{j}(\xi) \hat{f}(\xi)) = \int K_{j}(x - y) f(y) dy.
\end{equation}

\noindent $K_{j}$ is an $L^{1}$ kernel. When $j$ is an integer less than zero let $P_{j} f = 0$. Finally let

\begin{equation}\label{2.4}
P_{j_{1} \leq \cdot \leq j_{2}} f = \sum_{j_{1} \leq j \leq j_{2}} P_{j} f.
\end{equation}

\noindent We also define the frequency truncation

\begin{equation}\label{2.5}
 P_{\leq N} f = \mathcal F(\phi(\frac{\xi}{N}) \hat{f}(\xi)).
\end{equation}
\end{definition}

\noindent The Littlewood - Paley decomposition respects $L^{p}$ norms for $1 < p < \infty$.

\begin{lemma}[Littlewood - Paley theorem]\label{l2.1.1}
For $1 < p < \infty$,

\begin{equation}\label{2.5.1}
\| f \|_{L^{p}(\mathbf{R}^{2})} \sim_{p} \| (\sum_{j = 0}^{\infty} |P_{j} f|^{2})^{1/2} \|_{L^{p}(\mathbf{R}^{2})}.
\end{equation}
\end{lemma}

\noindent \emph{Proof:} See \cite{St1}, \cite{St}, \cite{T3} or \cite{T1}. $\Box$\vspace{5mm}

\noindent Because $\xi(t)$ is free to move around in the nonradial case, we will also define a Littlewood - Paley projection centered around $\xi_{0} \in \mathbf{R}^{2}$, $\xi_{0} \neq 0$.

\begin{definition}\label{d2.2}
Let $\xi_{0} \in \mathbf{R}^{2}$. Then define

\begin{equation}\label{2.6}
P_{\xi_{0}, j} u = e^{ix \cdot \xi_{0}} P_{j} (e^{-ix \cdot \xi_{0}} u).
\end{equation}



\noindent Also for $1 \leq p \leq \infty$ define the norm

\begin{equation}\label{2.7.1}
\| P_{\xi(t), j} f \|_{L_{t}^{p} L_{x}^{q}(I \times \mathbf{R}^{2})} = \| \| P_{\xi(t), j} f(t) \|_{L_{x}^{q}(\mathbf{R}^{2})} \|_{L_{t}^{p}(I)}.
\end{equation}
\end{definition}

\noindent \textbf{Remark:} Notice that

\begin{equation}\label{2.7.2}
P_{\xi_{0}, j} f = \int K_{j}(x - y) e^{i(x - y) \cdot \xi_{0}} f(y) dy,
\end{equation}

\noindent which also has an $L^{1}$ kernel, and $P_{\xi(t), j}$ satisfies $(\ref{1.12.1})$. Thus we can use $(\ref{1.13})$ and $(\ref{1.14})$. The same also holds for $P_{\xi(t), \leq j}$ and $P_{\xi(t), \geq j}$.\vspace{5mm}

\noindent To define our long time Strichartz spaces we utilize a class of function spaces first introduced in \cite{KoTa}. \cite{KoTa1}, \cite{KoTa2} applied these spaces to nonlinear Schr{\"o}dinger problems. See \cite{HHK} for a general description of these spaces. For critical problems, $U_{\Delta}^{p}$ spaces are more useful than the $X^{s, b}$ spaces, since the $X^{s,b}$ spaces of \cite{B}, \cite{B1} (see also \cite{Gin}) are not scale invariant except at $b = \frac{1}{2}$, which has the same difficulty as the failure of the embedding $\dot{H}^{1/2}(\mathbf{R}) \subset L^{\infty}(\mathbf{R})$. See \cite{Tat} for a more detailed discussion of this fact.

\begin{definition}[$U_{\Delta}^{p}$ spaces]\label{d2.3}
Let $1 \leq p < \infty$. Let $U_{\Delta}^{p}$ be an atomic space whose atoms are piecewise solutions to the linear equation. We call $u_{\lambda}$ a $U_{\Delta}^{p}$ atom if there exists an increasing sequence $\{ t_{k} \}_{k = 1}^{N}$, $N$ may be finite or infinite, $t_{0} = -\infty$, $N$ is finite, $t_{N + 1} = +\infty$, and

\begin{equation}\label{2.8}
u_{\lambda} = \sum_{k = 0}^{N} 1_{[t_{k}, t_{k + 1})} e^{it \Delta} u_{k}, \hspace{5mm} \text{and} \hspace{5mm} \sum_{k} \| u_{k} \|_{L^{2}}^{p} = 1.
\end{equation}

\noindent If $J \subset \mathbf{R}$ is an interval we say that $u_{\lambda}$ is a $U_{\Delta}^{p}(J)$ atom if $t_{k} \in J$ for all $1 \leq k \leq N$. Then for any $1 \leq p < \infty$, let

\begin{equation}\label{2.9}
\| u \|_{U_{\Delta}^{p}(J \times \mathbf{R}^{2})} = \inf \{ \sum_{\lambda} |c_{\lambda}| : u = \sum_{\lambda} c_{\lambda} u_{\lambda}, \hspace{5mm} \text{$u_{\lambda}$ are $U_{\Delta}^{p}(J)$ atoms} \}.
\end{equation}

\end{definition}

\noindent For any $1 \leq p < \infty$, $U_{\Delta}^{p}(J \times \mathbf{R}^{2}) \subset L^{\infty} L^{2}(J \times \mathbf{R}^{2})$. Additionally, $U_{\Delta}^{p}$ functions are continuous except at countably many points and right continuous everywhere.

\begin{definition}[$V_{\Delta}^{p}$ spaces]\label{d2.4}
Let $1 \leq p < \infty$. Then $V_{\Delta}^{p}$ is the space of functions $u \in L^{\infty}(L^{2})$ such that

\begin{equation}\label{2.10}
\| v \|_{V_{\Delta}^{p}}^{p} = \| v \|_{L_{t}^{\infty} L_{x}^{2}}^{p} + \sup_{\{ t_{k} \} \nearrow} \sum_{k} \| e^{-it_{k} \Delta} v(t_{k}) - e^{-it_{k + 1} \Delta} v(t_{k + 1}) \|_{L_{x}^{2}}^{p}.
\end{equation}

\noindent The supremum is taken over increasing sequences $t_{k}$. If $J \subset \mathbf{R}$ then

\begin{equation}\label{2.10}
\| v \|_{V_{\Delta}^{p}(J \times \mathbf{R}^{2})}^{p} = \| v \|_{L_{t}^{\infty} L_{x}^{2}(J \times \mathbf{R}^{2})}^{p} + \sup_{\{ t_{k} \} \nearrow} \sum_{k} \| e^{-it_{k} \Delta} v(t_{k}) - e^{-it_{k + 1} \Delta} v(t_{k + 1}) \|_{L_{x}^{2}}^{p},
\end{equation}
where each $t_{k}$ lies in $J$. $\{ t_{k} \}$ may be a finite or infinite sequence.
\end{definition}

\noindent \textbf{Remark:} \cite{KoTa} required that the $V_{\Delta}^{2}$ norms be taken over right continuous functions, so as to eliminate functions that were zero almost everywhere in time. This distinction is unnecessary here.

\begin{theorem}\label{t2.5}
For functions that are right continuous in time, the function spaces $U_{\Delta}^{p}$ and $V_{\Delta}^{q}$ obey the embeddings

\begin{equation}\label{2.11}
U_{\Delta}^{p} \subset V_{\Delta}^{p} \subset U_{\Delta}^{q} \subset L^{\infty} (L^{2}), \hspace{5mm} p < q.
\end{equation}

\noindent These spaces are also closed under truncation in time. If $I = [a, b)$,

\begin{equation}\label{2.15}
\aligned
\chi_{I} : U_{\Delta}^{p} \rightarrow U_{\Delta}^{p}, \\
\chi_{I} : V_{\Delta}^{p} \rightarrow V_{\Delta}^{p}.
\endaligned
\end{equation}

\noindent Formally, let $DU_{\Delta}^{p}$ be the space of functions

\begin{equation}\label{2.12}
DU_{\Delta}^{p} = \{ (i \partial_{t} + \Delta)u ; u \in U_{\Delta}^{p} \},
\end{equation}

\noindent and then $DU_{\Delta}^{p} = (V_{\Delta}^{p'})^{\ast}$, with $\frac{1}{p} + \frac{1}{p'} = 1$. Then


\begin{equation}\label{2.14}
\| \int_{0}^{t} e^{i(t - \tau) \Delta} F(\tau) d\tau \|_{U_{\Delta}^{p}(J \times \mathbf{R}^{2})} \lesssim \sup \{ \int_{J} \langle v, F \rangle dt : \| v \|_{V_{\Delta}^{p'}(J \times \mathbf{R}^{2})} = 1 \}.
\end{equation}
\end{theorem}

\noindent \emph{Proof:} See \cite{KochTataruVisan} for a more detailed description of these spaces as well as proofs. $\Box$\vspace{5mm}

\noindent \textbf{Remark:} Since a solution to $(\ref{0.1})$ is continuous in time, we will often utilize the embeddings $(\ref{2.11})$.

\begin{lemma}\label{l2.6}
Suppose $J = I_{1} \cup I_{2}$, $I_{1} = [a, b]$, $I_{2} = [b, c]$, $a \leq b \leq c$. Then,

\begin{equation}\label{2.16}
\| u \|_{U_{\Delta}^{p}(J \times \mathbf{R}^{d})}^{p} \leq \| u \|_{U_{\Delta}^{p}(I_{1} \times \mathbf{R}^{d})}^{p} + \| u \|_{U_{\Delta}^{p}(I_{2} \times \mathbf{R}^{d})}^{p}.
\end{equation}


\end{lemma}

\noindent \emph{Proof:} (See $(29)$ of \cite{KoTa2}) Let $u$ be a $U_{\Delta}^{p}$ atom supported on $I_{1}$ and let $v$ be a $U_{\Delta}^{p}$ atom supported on $I_{2}$. If $$u = \sum_{k = 1}^{K} 1_{[t_{k}, t_{k + 1})} e^{it \Delta} u_{k},$$ and $$v = \sum_{l = 1}^{L} 1_{[t_{l}, t_{l + 1})} e^{it \Delta} v_{l},$$ then $$w = \frac{c_{1}}{(|c_{1}|^{p} + |c_{2}|^{p})^{1/p}} \sum_{k = 1}^{K} 1_{[t_{k}, t_{k + 1})} e^{it \Delta} u_{k} + \frac{c_{2}}{(|c_{1}|^{p} + |c_{2}|^{p})^{1/p}} \sum_{l = 1}^{L} 1_{[t_{l}, t_{l + 1})} e^{it \Delta} v_{l}$$

\noindent is also a $U_{\Delta}^{p}$ atom. Moreover, $$c_{1} u + c_{2} v = (|c_{1}|^{p} + |c_{2}|^{p})^{1/p} w.$$

\noindent Now take $$u = \sum_{\lambda_{1}} c_{\lambda_{1}} u_{\lambda_{1}}, \hspace{5mm} v = \sum_{\lambda_{2}} c_{\lambda_{2}} v_{\lambda_{2}}.$$ Suppose $\sum_{\lambda_{1}} |c_{\lambda_{1}}| = C_{1}$ and $\sum_{\lambda_{2}} |c_{\lambda_{2}}| = C_{2}$.

$$\sum_{\lambda_{1}} c_{\lambda_{1}} u_{\lambda_{1}} + \sum_{\lambda_{2}} c_{\lambda_{2}} v_{\lambda_{2}} = \sum_{\lambda_{1}, \lambda_{2}} \frac{c_{\lambda_{1}} |c_{\lambda_{2}}|}{C_{2}} u_{\lambda_{1}} + \sum_{\lambda_{1}, \lambda_{2}} \frac{|c_{\lambda_{1}}| c_{\lambda_{2}}}{C_{1}} v_{\lambda_{2}}.$$

\noindent Therefore,

$$\| u + v \|_{U_{\Delta}^{p}} \leq \sum_{\lambda_{1}, \lambda_{2}} |c_{\lambda_{1}}| |c_{\lambda_{2}}| \| \frac{u_{\lambda_{1}}}{C_{2}} + \frac{v_{\lambda_{2}}}{C_{1}} \|_{U_{\Delta}^{p}} \leq C_{1} C_{2} (\frac{1}{C_{2}^{p}} + \frac{1}{C_{1}^{p}})^{1/p} = (C_{1}^{p} + C_{2}^{p})^{1/p}.$$

\noindent This proves the lemma. $\Box$\vspace{5mm}

\noindent We next prove a lemma that will be used extensively in section four.

\begin{lemma}\label{l2.6.1}
Suppose $J = \cup_{m = 1}^{k} J^{m}$, where $J^{m}$ are consecutive intervals, $J^{m} = [a_{m}, b_{m}]$, $a_{m + 1} = b_{m}$. Also suppose that $F \in L_{t}^{1} L_{x}^{2}(J \times \mathbf{R}^{2})$ (however our bound will not depend on $\| F \|_{L_{t}^{1} L_{x}^{2}}$.) Then for any $t_{0} \in J$,

\begin{equation}\label{2.17.1}
\aligned
\| \int_{t_{0}}^{t} e^{i(t - \tau) \Delta} F(\tau) d\tau \|_{U_{\Delta}^{2}(J \times \mathbf{R}^{d})} \lesssim \sum_{m = 1}^{k} \| \int_{J^{m}} e^{-i \tau \Delta} F(\tau) d\tau \|_{L^{2}(\mathbf{R}^{d})} \\ + (\sum_{m = 1}^{k} (\sup_{\| v_{m} \|_{V_{\Delta}^{2}(J^{m} \times \mathbf{R}^{2})} = 1} \int_{J^{m}} \langle F(\tau), v_{m}(\tau) \rangle d\tau)^{2})^{1/2}.
\endaligned
\end{equation}
\end{lemma}

\noindent \emph{Proof:} The proof is the same in any dimension. Suppose $t_{0} \in J^{m^{\ast}}$, $1 \leq m^{\ast} \leq k$. Then for $t > t_{0}$, $t \in J$,

\begin{equation}\label{2.17.2}
\aligned
\int_{t_{0}}^{t} e^{i(t - \tau) \Delta} F(\tau) d\tau = 1_{[b_{m^{\ast}}, \infty)}(t) e^{it \Delta} \int_{t_{0}}^{b_{m^{\ast}}} e^{-i \tau \Delta} F(\tau) d\tau + 1_{[a_{m^{\ast}}, b_{m^{\ast}}]}(t) \int_{t_{0}}^{t} e^{i(t - \tau) \Delta} F(\tau) d\tau \\
+ \sum_{m^{\ast} < n < k} 1_{[b_{n}, \infty)}(t) e^{it \Delta} \int_{J^{n}} e^{-i \tau \Delta} F(\tau) d\tau + \sum_{m^{\ast} < n \leq k} 1_{[a_{n}, b_{n}]}(t) \int_{a_{n}}^{t} e^{i(t - \tau) \Delta} F(\tau) d\tau.
\endaligned
\end{equation}

\noindent By lemma $\ref{l2.6}$, and $(\ref{2.14})$,

\begin{equation}\label{2.17.3}
\aligned
\| 1_{[a_{m^{\ast}}, b_{m^{\ast}}]}(t) \int_{t_{0}}^{t} e^{i(t - \tau) \Delta} F(\tau) d\tau + \sum_{m^{\ast} < n \leq k} 1_{[a_{n}, b_{n}]}(t) \int_{a_{n}}^{t} e^{i(t - \tau) \Delta} F(\tau) d\tau \|_{U_{\Delta}^{2}(J \times \mathbf{R}^{d})} \\ \lesssim (\sum_{m^{\ast} \leq n \leq k} ( \sup_{\| v_{n} \|_{V_{\Delta}^{2}(J^{n} \times \mathbf{R}^{2})} = 1} \int_{J^{n}} \langle v_{n}(\tau), F(\tau) \rangle d\tau )^{2})^{1/2}.
\endaligned
\end{equation}

\noindent Next,

\begin{equation}\label{2.17.4}
\| \sum_{m^{\ast} < n < k} 1_{[b_{n}, \infty)}(t) e^{it \Delta} \int_{J^{n}} e^{-i \tau \Delta} F(\tau) d\tau \|_{U_{\Delta}^{2}(J \times \mathbf{R}^{d})} \leq \sum_{m^{\ast} < n < k} \| \int_{J^{n}} e^{-i \tau \Delta} F(\tau) d\tau \|_{L^{2}(\mathbf{R}^{d})}.
\end{equation}

\noindent Finally using $U_{\Delta}^{2} \subset L_{t}^{\infty} L_{x}^{2}$ in the last inequality,

\begin{equation}\label{2.17.5}
\aligned
\| 1_{[b_{m^{\ast}}, \infty)}(t) e^{it \Delta} \int_{t_{0}}^{b_{m^{\ast}}} e^{-i \tau \Delta} F(\tau) d\tau \|_{U_{\Delta}^{2}(J \times \mathbf{R}^{d})} \lesssim \| \int_{t_{0}}^{b_{m^{\ast}}} e^{-i \tau \Delta} F(\tau) d\tau \|_{L^{2}(\mathbf{R}^{d})} \\ 
\lesssim \sup_{\| v \|_{V_{\Delta}^{2}(J^{m^{\ast}} \times \mathbf{R}^{2})} = 1} \int_{J^{m^{\ast}}} \langle F(\tau), v(\tau) \rangle d\tau.
\endaligned
\end{equation}

\noindent A similar computation can be made for $t < t_{0}$, $t \in J$. $\Box$\vspace{5mm}

\noindent The $U_{\Delta}^{2}$ spaces respect linear and bilinear Strichartz estimates. Indeed, checking individual atoms shows that if $p, q$ is an admissible pair then

\begin{equation}\label{2.17.6}
\| u \|_{L_{t}^{p} L_{x}^{q}(I \times \mathbf{R}^{2})} \lesssim_{p,q} \| u \|_{U_{\Delta}^{p}(I \times \mathbf{R}^{2})}.
\end{equation}

\begin{lemma}\label{l2.6.2}
Suppose that under some condition on the supports of $\hat{u}_{0}(\xi)$, $\hat{v}_{0}(\xi)$, $p < \infty$,

\begin{equation}\label{2.17.7}
\| (e^{it \Delta} u_{0})(e^{it \Delta} v_{0}) \|_{L_{t}^{p} L_{x}^{q}(I \times \mathbf{R}^{2})} \lesssim N^{-\alpha} \| u_{0} \|_{L^{2}(\mathbf{R}^{2})} \| v_{0} \|_{L^{2}(\mathbf{R}^{2})}.
\end{equation}

\noindent Then if $\hat{u}(t, \xi)$ and $\hat{v}(t, \xi)$ are under the same conditions,

\begin{equation}\label{2.17.8}
\| uv \|_{L_{t}^{p} L_{x}^{q}(I \times \mathbf{R}^{2})} \lesssim N^{-\alpha} \| u \|_{U_{\Delta}^{p}(I \times \mathbf{R}^{2})} \| v \|_{U_{\Delta}^{p}(I \times \mathbf{R}^{2})}.
\end{equation}
\end{lemma}

\noindent \emph{Proof:} It suffices to check individual atoms. Suppose $u$ is an atom, $t_{k} \nearrow$, and $u = \sum_{k} 1_{[t_{k}, t_{k + 1})}(t) e^{it \Delta} u_{k}$. Then

\begin{equation}\label{2.17.9}
 \| uv \|_{L_{t}^{p} L_{x}^{q}(I \times \mathbf{R}^{2})}^{p} = \sum_{k} \| (e^{it \Delta} u_{k}) v \|_{L_{t}^{p} L_{x}^{q}([t_{k}, t_{k + 1}] \times \mathbf{R}^{2})}^{p}.
\end{equation}

\noindent Then if $v$ is also an atom, there exists $t_{l}' \nearrow$ with $v = \sum_{l} 1_{[t_{l}', t_{l + 1}')}(t) e^{it \Delta} v_{l}$,

\begin{equation}\label{2.17.10}
\| (e^{it \Delta} u_{k}) v \|_{L_{t}^{p} L_{x}^{q}(I \times \mathbf{R}^{2})}^{p} = \sum_{l} \| (e^{it \Delta} u_{k})(e^{it \Delta} v_{l}) \|_{L_{t}^{p} L_{x}^{q}([t_{l}', t_{l + 1}'] \times \mathbf{R}^{2})}^{p}
\end{equation}

\begin{equation}\label{2.17.11}
 \lesssim N^{-\alpha} \sum_{l} \| u_{k} \|_{L^{2}}^{p} \| v_{l} \|_{L^{2}}^{p} \lesssim N^{-\alpha p} \| u_{k} \|_{L^{2}}^{p}.
\end{equation}

\noindent This implies

\begin{equation}\label{2.17.12}
 \sum_{k} \| (e^{it \Delta} u_{k}) v \|_{L_{t}^{p} L_{x}^{q}([t_{k}, t_{k + 1}] \times \mathbf{R}^{2})}^{p} \lesssim N^{-\alpha p}.
\end{equation}

\noindent $\Box$\vspace{5mm}

\noindent Now we are ready to define the long time Strichartz norm $\tilde{X}_{k_{0}}([0, T] \times \mathbf{R}^{2})$. Fix three constants

\begin{equation}\label{2.17.13}
0 < \epsilon_{3} << \epsilon_{2} << \epsilon_{1} < 1.
\end{equation}

\noindent Fix an integer $k_{0} \in \mathbf{Z}_{\geq 0}$ and suppose that $M = 2^{k_{0}}$. Let $[a, b]$ be an interval such that

\begin{equation}\label{2.18}
\int_{a}^{b} \int |u(t,x)|^{4} dx dt = M,
\end{equation}

\noindent and

\begin{equation}\label{2.19}
\int_{a}^{b} N(t)^{3} dt = \epsilon_{3} M.
\end{equation}

\noindent Notice that $(\ref{2.18})$ is invariant under the scaling $(\ref{0.2})$, while $(\ref{2.19})$ is not. Therefore, given an interval that satisfies $(\ref{2.18})$ it is always possible to rescale so that $(\ref{2.19})$ is also satisfied.\vspace{5mm}

\noindent Now by $(\ref{1.19})$, $(\ref{1.20})$, and $(\ref{1.29})$ it is possible to choose $\epsilon_{1}$, $\epsilon_{2}$, and $\epsilon_{3}$ which satisfy $(\ref{2.17.13})$, and also

\begin{equation}\label{2.20}
 |\xi'(t)| + |N'(t)| \leq 2^{-20} \frac{N(t)^{3}}{\epsilon_{1}^{1/2}},
\end{equation}

\begin{equation}\label{2.21}
\int_{|x - x(t)| \geq \frac{2^{-20} \epsilon_{3}^{-1/4}}{N(t)}} |u(t,x)|^{2} dx + \int_{|\xi - \xi(t)| \geq 2^{-20} \epsilon_{3}^{-1/4} N(t)} |\hat{u}(t,\xi)|^{2} d\xi \leq \epsilon_{2}^{2},
\end{equation}

\noindent and

\begin{equation}\label{2.22}
 \epsilon_{3} < \epsilon_{2}^{10}.
\end{equation}

\noindent Next, partition $[a, b]$ in two different ways:

\begin{definition}[Small intervals]\label{d2.7}
Let $[a, b] = \cup_{l = 0}^{M - 1} J_{l}$, with $\| u \|_{L_{t,x}^{4}(J_{l} \times \mathbf{R}^{2})} = 1$. We will call the intervals $J_{l}$ the small intervals. Recall from $(\ref{1.27})$ the notation $N(J_{l}) = \sup_{t \in J_{l}} N(t)$.
\end{definition}

\noindent \textbf{Remark:} Theorem $\ref{t1.1}$ combined with conservation of mass implies that for any admissible pair $(p, q)$,

\begin{equation}\label{2.22.1}
\| u \|_{L_{t}^{p} L_{x}^{q}(J_{l} \times \mathbf{R}^{2})} \lesssim_{p,q} 1.
\end{equation}

\begin{definition}[$J^{\alpha}$ intervals]\label{d2.8}
Let $[a, b] = \cup_{\alpha = 0}^{M - 1} J^{\alpha}$, $\alpha = 0, ..., M - 1$, such that

\begin{equation}\label{2.24}
\int_{J^{\alpha}} (N(t)^{3} + \epsilon_{3} \| u(t) \|_{L_{x}^{4}(\mathbf{R}^{2})}^{4}) dt = 2 \epsilon_{3}.
\end{equation}
\end{definition}

\begin{definition}\label{d2.8}
For an integer $0 \leq j < k_{0}$, $0 \leq k < 2^{k_{0} - j}$, let

\begin{equation}\label{2.25}
 G_{k}^{j} = \cup_{\alpha = k 2^{j}}^{(k + 1) 2^{j} - 1} J^{\alpha}.
\end{equation}

\noindent For $j \geq k_{0}$ let $G_{k}^{j} = [a, b]$.\vspace{5mm}

\noindent Now suppose that $[t_{0}, t_{1}] = G_{k}^{j}$. Let $\xi(G_{k}^{j}) = \xi(t_{0})$, and define $\xi(J^{\alpha})$, and $\xi(J_{l})$ in a similar manner.
\end{definition}

\noindent \textbf{Remark:} By $(\ref{2.20})$, $(\ref{2.24})$, and $(\ref{2.25})$, for all $t \in G_{k}^{j}$,

\begin{equation}\label{2.26}
 |\xi(t) - \xi(G_{k}^{j})| \leq \int_{G_{k}^{j}} 2^{-20} \epsilon_{1}^{-1/2} N(t)^{3} dt \leq 2^{j - 19} \epsilon_{3} \epsilon_{1}^{-1/2}.
\end{equation}

\noindent Therefore, for all $t \in G_{k}^{j}$,

\begin{equation}\label{2.27}
 \{ \xi : 2^{j - 1} \leq |\xi - \xi(t)| \leq 2^{j + 1} \} \subset \{ \xi : 2^{j - 2} \leq |\xi - \xi(G_{k}^{j})| \leq 2^{j + 2} \} \subset \{ \xi : 2^{j - 3} \leq |\xi - \xi(t)| \leq 2^{j + 3} \},
\end{equation}

\noindent and

\begin{equation}\label{2.28}
 \{ \xi : |\xi - \xi(t)| \leq 2^{j + 1} \} \subset \{ \xi : |\xi - \xi(G_{k}^{j})| \leq 2^{j + 2} \} \subset \{ \xi : |\xi - \xi(t)| \leq 2^{j + 3} \}.
\end{equation}

\noindent Now we define our spaces on which we compute the long time Strichartz estimates. The $U_{\Delta}^{2}$ norm replaces the $L_{t}^{2} L_{x}^{\frac{2d}{d - 2}}$ norm in \cite{D2} for dimensions $d > 2$. On a first reading it may be convenient to skip ahead to the sketch of the proof at the beginning section four, and consult definition $\ref{d2.9}$ while reading the sketch.

\begin{definition}[$\tilde{X}_{k_{0}}$ spaces]\label{d2.9}





\noindent For any $G_{k}^{j} \subset [a, b]$ let

\begin{equation}\label{2.32}
\| u \|_{X(G_{k}^{j} \times \mathbf{R}^{2})}^{2} \equiv \sum_{0 \leq i < j} 2^{i - j} \sum_{G_{\alpha}^{i} \subset G_{k}^{j}} \| P_{\xi(G_{\alpha}^{i}), i - 2 \leq \cdot \leq i + 2} u \|_{U_{\Delta}^{2}(G_{\alpha}^{i} \times \mathbf{R}^{2})}^{2} + \sum_{i \geq j} \| P_{\xi(G_{k}^{j}), i - 2 \leq \cdot \leq i + 2} u \|_{U_{\Delta}^{2}(G_{k}^{j} \times \mathbf{R}^{2})}^{2}.
\end{equation}

\noindent Then define $\tilde{X}_{k_{0}}$ to be the supremum of $(\ref{2.32})$ over all intervals $G_{k}^{j} \subset [a, b]$ with $k \leq k_{0}$.

\begin{equation}\label{2.33}
 \| u \|_{\tilde{X}_{k_{0}}([a, b] \times \mathbf{R}^{2})}^{2} \equiv \sup_{0 \leq j \leq k_{0}} \sup_{G_{k}^{j} \subset [a, b]} \| u \|_{X(G_{k}^{j} \times \mathbf{R}^{2})}^{2}.
\end{equation}

\noindent Also for $0 \leq k_{\ast} \leq k_{0}$, let

\begin{equation}\label{2.34}
\| u \|_{\tilde{X}_{k_{\ast}}([a, b] \times \mathbf{R}^{2})}^{2} \equiv \sup_{0 \leq j \leq k_{\ast}} \sup_{G_{k}^{j} \subset [a, b]} \| u \|_{X(G_{k}^{j} \times \mathbf{R}^{2})}^{2}.
\end{equation}

\noindent $\| u \|_{\tilde{X}_{k_{\ast}}(G_{k}^{j} \times \mathbf{R}^{2})}$, $k_{\ast} \leq j$ is defined in a similar manner.\vspace{5mm}
\end{definition}

\noindent \textbf{Remark:} By $(\ref{2.26})$, for any $t \in G_{\alpha}^{i}$,

\begin{equation}\label{2.32.1}
P_{\xi(t), i} P_{\xi(G_{\alpha}^{i}), i - 2 \leq \cdot \leq i + 2} = P_{\xi(t), i}.
\end{equation}

\begin{definition}[$\tilde{Y}_{k_{0}}$ spaces]\label{d2.10}
The $\tilde{Y}_{k_{0}}$ norm measures the $\tilde{X}_{k_{0}}$ norm of $u$ at scales much higher than $N(t)$. This norm provides some crucial ``smallness", closing a bootstrap argument in the next section. Let

\begin{equation}\label{2.35}
\aligned
 \| u \|_{Y(G_{k}^{j} \times \mathbf{R}^{2})}^{2} = \sum_{0 < i < j} 2^{i - j} \sum_{G_{\alpha}^{i} \subset G_{k}^{j} : N(G_{\alpha}^{i}) \leq 2^{i - 5} \epsilon_{3}^{1/2}} \| P_{\xi(G_{\alpha}^{i}), i - 2 \leq \cdot \leq i + 2} u \|_{U_{\Delta}^{2}(G_{\alpha}^{i} \times \mathbf{R}^{2})}^{2} \\ + \sum_{i \geq j, i > 0; N(G_{k}^{j}) \leq \epsilon_{3}^{1/2} 2^{i - 5}} \| P_{\xi(G_{k}^{j}), i - 2 \leq \cdot \leq i + 2} u \|_{U_{\Delta}^{2}(G_{k}^{j} \times \mathbf{R}^{2})}^{2}.
\endaligned
\end{equation}



\noindent Define $\| u \|_{\tilde{Y}_{k_{\ast}}([0, T] \times \mathbf{R}^{2})}$ using $Y(G_{k}^{j} \times \mathbf{R}^{2})$ in the same manner as $\| u \|_{\tilde{X}_{k_{\ast}}([0, T] \times \mathbf{R}^{2})}$ was defined using $X(G_{k}^{j} \times \mathbf{R}^{2})$.\vspace{5mm}
\end{definition}

\noindent Then by $(\ref{2.27})$ and $(\ref{2.28})$, for $i < j$, $(p, q)$ an admissible pair, by the definition of $\tilde{X}_{j}$,

\begin{equation}\label{2.37}
\aligned
\| P_{\xi(t), i} u \|_{L_{t}^{p} L_{x}^{q}(G_{k}^{j} \times \mathbf{R}^{2})} = (\sum_{G_{\alpha}^{i} \subset G_{k}^{j}} \| P_{\xi(t), i} u \|_{L_{t}^{p} L_{x}^{q}(G_{\alpha}^{i} \times \mathbf{R}^{2})}^{p})^{1/p} \\ \lesssim (\sum_{G_{\alpha}^{i} \subset G_{k}^{j}} \| P_{\xi(t), i} u \|_{L_{t}^{p} L_{x}^{q}(G_{\alpha}^{i} \times \mathbf{R}^{2})}^{2})^{1/p} (\sup_{G_{\alpha}^{i} \subset G_{k}^{j}} \| P_{\xi(t), i} u \|_{L_{t}^{p} L_{x}^{q}(G_{\alpha}^{i} \times \mathbf{R}^{2})})^{1 - \frac{2}{p}} \\ \lesssim_{p,q} (\sum_{G_{\alpha}^{i} \subset G_{k}^{j}} \| P_{\xi(G_{\alpha}^{i}), i - 2 \leq \cdot \leq i + 2} u \|_{U_{\Delta}^{2}(G_{\alpha}^{i} \times \mathbf{R}^{2})}^{2})^{1/p} (\sup_{G_{\alpha}^{i} \subset G_{k}^{j}} \| P_{\xi(G_{\alpha}^{i}), i - 2 \leq \cdot \leq i + 2} u \|_{U_{\Delta}^{2}(G_{\alpha}^{i} \times \mathbf{R}^{2})})^{1 - \frac{2}{p}} \\ \lesssim
2^{\frac{(j - i)}{p}} \| u \|_{X(G_{k}^{j} \times \mathbf{R}^{2})}^{2/p} \| u \|_{\tilde{X}_{j}(G_{k}^{j} \times \mathbf{R}^{2})}^{1 - 2/p} \lesssim 2^{\frac{(j - i)}{p}} \| u \|_{\tilde{X}_{j}(G_{k}^{j} \times \mathbf{R}^{2})}.
\endaligned
\end{equation}

\noindent Also, by $(\ref{2.5.1})$,

\begin{equation}\label{2.38}
\aligned
\| P_{\xi(t), \geq j} u \|_{L_{t}^{p} L_{x}^{q}(G_{k}^{j} \times \mathbf{R}^{2})} \sim_{q} \| (\sum_{l \geq j} |P_{\xi(t), l} u|^{2})^{1/2} \|_{L_{t}^{p} L_{x}^{q}(G_{k}^{j} \times \mathbf{R}^{2})} \\
\lesssim (\sum_{l \geq j} \| P_{\xi(t), l} u \|_{L_{t}^{p} L_{x}^{q}(G_{k}^{j} \times \mathbf{R}^{2})}^{2})^{1/2} \lesssim_{p,q} \| u \|_{X(G_{k}^{j} \times \mathbf{R}^{2})}.
\endaligned
\end{equation}

\section{Long time Strichartz estimate}

\begin{theorem}[Long time Strichartz estimate]\label{t3.1}
If $u$ is an almost periodic solution to $(\ref{0.1})$ then for any $M = 2^{k_{0}}$, $\epsilon_{1}$, $\epsilon_{2}$, $\epsilon_{3}$ satisfying $(\ref{2.20})$ - $(\ref{2.22})$, $\int_{0}^{T} N(t)^{3} dt = \epsilon_{3} M$, and $\int_{0}^{T} \int |u(t,x)|^{4} dx dt = M$,

\begin{equation}\label{3.1}
\| u \|_{\tilde{X}_{k_{0}}([0, T] \times \mathbf{R}^{2})} \lesssim 1.
\end{equation}
\end{theorem}

\noindent \textbf{Remark:} Throughout this section the implicit constant depends only on $u$, and not on $M$, or $\epsilon_{1}$, $\epsilon_{2}$, $\epsilon_{3}$.\vspace{5mm}

\noindent In this section we will prove the long time Strichartz estimate. This proof will occupy the bulk of the paper, encompassing all of section four as well as an appendix in section six. The idea of the proof is actually quite similar to proof of the long time Strichartz estimates in \cite{D2}, although much more technically complicated.\vspace{5mm}

\noindent In \cite{D2}, the long time Strichartz estimates appeared in theorem $1.24$, which stated that for an almost periodic solution to the mass - critical problem,

\begin{equation}\label{3.0.1}
\| P_{|\xi - \xi(t)| > N} u \|_{L_{t}^{2} L_{x}^{\frac{2d}{d - 2}}(J \times \mathbf{R}^{d})} \lesssim (\frac{K}{N})^{1/2} + 1,
\end{equation}

\noindent when $J$ is an interval satisfying

\begin{equation}\label{3.0.2}
\int_{J} N(t)^{3} dt = K.
\end{equation}

\noindent The corresponding result is unavailable in two dimensions, due to the failure of the endpoint Strichartz estimate $(2, \infty)$ to hold. The argument could not even be duplicated in the radial case, since the argument in \cite{D2} relies on the double endpoint Strichartz estimates.\vspace{5mm}

\noindent Instead, in section three, culminating with definitions $\ref{d2.9}$ and $\ref{d2.10}$, we constructed a function space which mimics the essential features of $(\ref{3.0.1})$. Observe that if $U_{\Delta}^{2}$ were replaced with $L_{t}^{2} L_{x}^{\infty}$, then $(\ref{2.32})$ and $(\ref{2.33})$ would be almost completely identical to $(\ref{3.0.1})$ and $(\ref{3.0.2})$ holding for any compact interval $J$ with $K \leq 2^{k_{0}}$.\vspace{5mm}

\noindent This lack of the double endpoint Strichartz estimate is also why the proof of theorem $\ref{t3.1}$ is so long. Suppose for a moment that the endpoint Strichartz estimate did hold in dimension $d = 2$ and that $\xi(t) \equiv 0$. Then for any $\eta > 0$, by concentration compactness,

\begin{equation}\label{3.0.3}
\aligned
\| P_{> N} u \|_{L_{t}^{2} L_{x}^{\infty}(J \times \mathbf{R}^{2})} \lesssim \| P_{> N} u \|_{L_{t}^{\infty} L_{x}^{2}(J \times \mathbf{R}^{2})} + \| P_{> N} F(u) \|_{L_{t}^{2} L_{x}^{1}(J \times \mathbf{R}^{2})}, \\
\lesssim 1 + \| P_{> \frac{N}{8}} u \|_{L_{t}^{2} L_{x}^{\infty}} \| P_{ > C(\eta) N(t)} u \|_{L_{t}^{\infty} L_{x}^{2}(J \times \mathbf{R}^{2})}^{2} + \| (P_{> \frac{N}{8}} u)(P_{ \leq C(\eta) N(t)} u)^{2} \|_{L_{t}^{2} L_{x}^{1}(J \times \mathbf{R}^{2})} \\
\lesssim 1 + \eta \| P_{> \frac{N}{8}} u \|_{L_{t}^{2} L_{x}^{\infty}(J \times \mathbf{R}^{2})} + \| (P_{> \frac{N}{8}} u)(P_{ \leq C(\eta) N(t)} u) \|_{L_{t,x}^{2}(J \times \mathbf{R}^{2})}.
\endaligned
\end{equation}

\noindent Then making a straightforward bilinear argument would prove the estimate corresponding to $(\ref{3.0.1})$ by induction on $N$, starting with $N = 1$.\vspace{5mm}

\noindent The reader should observe how the endpoint Strichartz estimate facilitated this computation; because we are projecting to frequencies higher than $N$, one of the terms in $F(u) = |u|^{2} u$ must be at a frequency higher than $\frac{N}{8}$. The other two terms may then be split into a piece at frequencies $\geq C(\eta) N(t)$, where the mass is small, and a piece at frequencies $\leq C(\eta) N(t)$. The pieces where mass is small give the gain $\eta \| P_{> \frac{N}{8}} u \|_{L_{t}^{2} L_{x}^{\infty}(J \times \mathbf{R}^{2})}$, while the pieces where mass is large may be estimated using the bilinear estimate.\vspace{5mm}

\noindent Replacing $L_{t}^{2} L_{x}^{\infty}$ with an admissible pair, say $L_{t}^{3} L_{x}^{6}$, is not quite good enough, because projecting $F(u) = |u|^{2} u$ to high frequencies means that only one component of the product lies at frequencies $> \frac{N}{8}$, not two. So instead, we substitute $L_{t}^{2} L_{x}^{\infty}$ with $U_{\Delta}^{2}$ and rely heavily on bilinear estimates.\vspace{5mm}

\noindent The proof is also complicated by the fact that $DU_{\Delta}^{2}$ is the dual of $V_{\Delta}^{2}$ and not $U_{\Delta}^{2} \subset V_{\Delta}^{2}$. If $DU_{\Delta}^{2}$ was the dual of $U_{\Delta}^{2}$, then

\begin{equation}\label{3.0.4}
\| P_{> N} F(u) \|_{DU_{\Delta}^{2}(J \times \mathbf{R}^{2})} = \sup_{\| v \|_{U_{\Delta}^{2}(J \times \mathbf{R}^{2})} = 1} \int_{J} \langle P_{> N} v, F(u) \rangle dt,
\end{equation}

\noindent so

\begin{equation}\label{3.0.5}
\| P_{> N} F(u) \|_{DU_{\Delta}^{2}(J \times \mathbf{R}^{2})} \leq \sup_{\| v \|_{U_{\Delta}^{2}(J \times \mathbf{R}^{2})} = 1} \| (P_{> N} v) u \|_{L_{t,x}^{2}(J \times \mathbf{R}^{2})} \| (P_{> \frac{N}{8}} u) u \|_{L_{t,x}^{2}(J \times \mathbf{R}^{2})}.
\end{equation}

\noindent Then estimate $(\ref{3.0.5})$ uses the bilinear estimates in theorems $\ref{t3.2.0}$ and $\ref{t3.2.1}$. These bilinear estimates use the interaction Morawetz estimates of \cite{PV}, and give a logarithmic improvement over a simple application of the bilinear estimates of \cite{B2}. Because of $(\ref{2.14})$, we utilize lemma $\ref{l2.6.1}$.\vspace{5mm}

\noindent There are a couple of other technical complications which are worked out in the proof. First, splitting up $u$ into a high frequency piece and a low frequency piece and then taking the $L_{t}^{\infty} L_{x}^{2}$ norm of the high frequency piece is relatively straightforward. However, computing the $U_{\Delta}^{2}$ norm is far more technically complicated. This is the reason for the introduction of the $\tilde{Y}_{k_{0}}$ spaces (definition $\ref{d2.10}$).\vspace{5mm}

\noindent Also, notice that $(\ref{2.32})$ with $U_{\Delta}^{2}$ replaced with $L_{t}^{2} L_{x}^{\infty}$ and $(\ref{3.0.1})$ are not exactly the same. $(\ref{2.32})$ is an estimate on the $l^{2}$ summation at different frequencies, but $(\ref{3.0.1})$ is only a $l^{\infty}$ estimate at different frequencies. The $l^{2}$ summation in $(\ref{2.32})$ is combined with theorem $\ref{t3.2.2}$, which gives an $l^{2}$ summation improvement over theorem $\ref{t3.2.1}$.\vspace{5mm}

\noindent Having sketched the proof, we move to the details.\vspace{5mm}





\noindent \emph{Proof of theorem $\ref{t3.1}$:} We wish to prove that for any $0 \leq j \leq k_{0}$ and $G_{k}^{j} \subset [0, T]$,

\begin{equation}\label{3.1.1}
\sum_{0 \leq i \leq j} 2^{i - j} \sum_{G_{\alpha}^{i} \subset G_{k}^{j}} \| P_{\xi(G_{\alpha}^{i}), i - 2 \leq \cdot \leq i + 2} u \|_{U_{\Delta}^{2}(G_{\alpha}^{i} \times \mathbf{R}^{2})}^{2} + \sum_{i > j} \| P_{\xi(G_{k}^{j}), i - 2 \leq \cdot \leq i + 2} u \|_{U_{\Delta}^{2}(G_{k}^{j} \times \mathbf{R}^{2})}^{2} \lesssim 1.
\end{equation}

\noindent First observe that since $V_{\Delta}^{2} \subset U_{\Delta}^{4}$, $(\ref{2.14})$, $(\ref{2.17.6})$, and $(\ref{2.24})$ imply that $\| u \|_{U_{\Delta}^{2}(J^{\alpha} \times \mathbf{R}^{2})} \lesssim 1$ for each $J^{\alpha} \subset G_{k}^{j}$. Therefore, by definition $\ref{d2.9}$,

\begin{equation}
 \| u \|_{\tilde{X}_{0}([0, T] \times \mathbf{R}^{2})} \leq C(u).
\end{equation}

\noindent Also by definition $\ref{d2.10}$, $(\ref{2.20})$, $(\ref{2.21})$, and Duhamel's principle,

\begin{equation}\label{3.8}
\aligned
(\sum_{i > 0 : N(J^{\alpha}) \leq \epsilon_{3}^{1/2} 2^{i - 5}} &\| P_{\xi(J^{\alpha}), i - 2 \leq \cdot \leq i + 2} u \|_{U_{\Delta}^{2}(J^{\alpha} \times \mathbf{R}^{2})}^{2})^{1/2} \\ &\lesssim \| P_{\xi(t), \geq 4 \epsilon_{3}^{-1/2} N(t)} u \|_{L_{t}^{\infty} L_{x}^{2}(J^{\alpha} \times \mathbf{R}^{2})} + \| P_{\xi(t), \geq 4 \epsilon_{3}^{-1/2} N(t)} F(u) \|_{L_{t}^{1} L_{x}^{2}(J^{\alpha} \times \mathbf{R}^{2})} \\ &\lesssim \| P_{\xi(t), \geq \epsilon_{3}^{-1/2} N(t)} u \|_{L_{t}^{\infty} L_{x}^{2}(J^{\alpha} \times \mathbf{R}^{2})}^{3/4} (\| u \|_{L_{t}^{\infty} L_{x}^{2}(J^{\alpha} \times \mathbf{R}^{2})}^{1/4} + \| u \|_{L_{t}^{9/4} L_{x}^{18}(J^{\alpha} \times \mathbf{R}^{2})}^{9/4}) \lesssim \epsilon_{2}^{3/4}.
\endaligned
\end{equation}

\noindent Therefore,

\begin{equation}\label{3.8.1}
\| u \|_{\tilde{Y}_{0}([0, T] \times \mathbf{R}^{2})} \leq C(u) \epsilon_{2}^{3/4}.
\end{equation}

\noindent Moreover it is clear from definitions $\ref{d2.9}$ and $\ref{d2.10}$ that for any $0 \leq k_{\ast} < k_{0}$,

\begin{equation}\label{3.9}
\aligned
 \| u \|_{\tilde{X}_{k_{\ast} + 1}([0, T] \times \mathbf{R}^{2})}^{2} \leq 2 \| u \|_{\tilde{X}_{k_{\ast}}([0, T] \times \mathbf{R}^{2})}^{2}, \\
\| u \|_{\tilde{Y}_{k_{\ast} + 1}([0, T] \times \mathbf{R}^{2})}^{2} \leq 2 \| u \|_{\tilde{Y}_{k_{\ast}}([0, T] \times \mathbf{R}^{2})}^{2},
\endaligned
\end{equation}

\noindent so

\begin{equation}\label{3.9.0}
\aligned
 \| u \|_{\tilde{X}_{11}([0, T] \times \mathbf{R}^{2})}^{2} \leq 2^{11} C(u), \\
\| u \|_{\tilde{Y}_{11}([0, T] \times \mathbf{R}^{2})}^{2} \leq 2^{11} C(u) \epsilon_{2}^{3/4}.
\endaligned
\end{equation}

\noindent $(\ref{3.9.0})$ also implies that for any $j > 11$ and $G_{k}^{j} \subset [0, T]$,

\begin{equation}\label{3.9.0.1}
\sum_{0 \leq i \leq 11} 2^{i - j} \sum_{G_{\alpha}^{i} \subset G_{k}^{j}} \| P_{\xi(G_{\alpha}^{i}), i - 2 \leq \cdot \leq i + 2} u \|_{U_{\Delta}^{2}(G_{\alpha}^{i} \times \mathbf{R}^{2})}^{2} \leq 2^{11} C(u),
\end{equation}

\noindent and

\begin{equation}\label{3.9.0.2}
\sum_{0 < i \leq 11} 2^{i - j} \sum_{G_{\alpha}^{i} \subset G_{k}^{j} ; N(G_{\alpha}^{i}) \leq \epsilon_{3}^{1/2} 2^{i - 5}} \| P_{\xi(G_{\alpha}^{i}), i - 2 \leq \cdot \leq i + 2} u \|_{U_{\Delta}^{2}(G_{\alpha}^{i} \times \mathbf{R}^{2})}^{2} \leq C(u) \epsilon_{2}^{3/2}.
\end{equation}

\noindent Fix $k_{0}$, $12 \leq j \leq k_{0}$, and $G_{k}^{j} \subset [0, T]$. For $11 \leq i < j$, Duhamel's principle implies

\begin{equation}\label{3.2}
\aligned
\| P_{\xi(G_{\alpha}^{i}), i - 2 \leq \cdot \leq i + 2} u \|_{U_{\Delta}^{2}(G_{\alpha}^{i} \times \mathbf{R}^{2})} &\lesssim \| P_{\xi(G_{\alpha}^{i}), i - 2 \leq \cdot \leq i + 2} u(t_{\alpha}^{i}) \|_{L^{2}(\mathbf{R}^{2})} \\ + &\| \int_{t_{\alpha}^{i}}^{t} e^{i(t - \tau) \Delta} P_{\xi(G_{\alpha}^{i}), i - 2 \leq \cdot \leq i + 2} F(u(\tau)) d\tau \|_{U_{\Delta}^{2}(G_{\alpha}^{i} \times \mathbf{R}^{2})}.
\endaligned
\end{equation}

\noindent Choose $t_{\alpha}^{i}$ satisfying

\begin{equation}\label{3.3}
\| P_{\xi(G_{\alpha}^{i}), i - 2 \leq \cdot \leq i + 2} u(t_{\alpha}^{i}) \|_{L^{2}(\mathbf{R}^{2})} = \inf_{t \in G_{\alpha}^{i}} \| P_{\xi(G_{\alpha}^{i}), i - 2 \leq \cdot \leq i + 2} u(t) \|_{L^{2}(\mathbf{R}^{2})}.
\end{equation}

\noindent Then by $(\ref{2.24})$, $(\ref{2.27})$, and $(\ref{2.28})$, 

\begin{equation}\label{3.4}
 \sum_{11 \leq i < j} 2^{i - j} \sum_{G_{\alpha}^{i} \subset G_{k}^{j}} \| P_{\xi(G_{\alpha}^{i}), i - 2 \leq \cdot \leq i + 2} u(t_{\alpha}^{i}) \|_{L^{2}(\mathbf{R}^{2})}^{2}
\end{equation}

\begin{equation}\label{3.5}
\aligned
 \lesssim 2^{-j} \epsilon_{3}^{-1} \int_{G_{k}^{j}} (\epsilon_{3} \| u(t) \|_{L^{4}(\mathbf{R}^{2})}^{4} + N(t)^{3}) \sum_{11 \leq i < j} \| P_{\xi(t), i - 3 \leq \cdot \leq i + 3} u(t) \|_{L^{2}(\mathbf{R}^{2})}^{2} dt \\ \lesssim 2^{-j} \epsilon_{3}^{-1} \| u \|_{L_{t}^{\infty} L^{2}([0, T] \times \mathbf{R}^{2})}^{2} \int_{G_{k}^{j}} (N(t)^{3} + \epsilon_{3} \| u(t) \|_{L^{4}(\mathbf{R}^{2})}^{4}) dt \lesssim 1.
\endaligned
\end{equation}

\noindent For $i \geq j$ simply take $t_{\alpha}^{i} = t_{0}$, where $t_{0}$ is a fixed element of $G_{k}^{j}$, say the left endpoint. Then

\begin{equation}\label{3.6}
\sum_{i \geq j} \| P_{\xi(G_{k}^{j}), i - 2 \leq \cdot \leq i + 2} u(t_{0}) \|_{L^{2}(\mathbf{R}^{2})}^{2} \lesssim \| u(t_{0}) \|_{L^{2}(\mathbf{R}^{2})}^{2} \lesssim 1.
\end{equation}

\noindent Therefore,

\begin{equation}\label{3.6.1}
\sum_{0 \leq i < j} 2^{i - j} \sum_{G_{\alpha}^{i} \subset G_{k}^{j}} \| P_{\xi(G_{\alpha}^{i}), i - 2 \leq \cdot \leq i + 2} u(t_{\alpha}^{i}) \|_{L_{x}^{2}(\mathbf{R}^{2})}^{2} + \sum_{i \geq j} \| P_{\xi(G_{k}^{j}), i - 2 \leq \cdot \leq i + 2} u(t_{0}) \|_{L^{2}(\mathbf{R}^{2})}^{2} \lesssim 1,
\end{equation}

\noindent so $(\ref{3.9.0.1})$ and $(\ref{3.6.1})$ imply

\begin{equation}\label{3.7}
\aligned
 \| u \|_{X(G_{k}^{j} \times \mathbf{R}^{2})}^{2} \lesssim 1 + \sum_{i \geq j; i \geq 11} \| \int_{t_{\alpha}^{i}}^{t} e^{i(t - \tau) \Delta} P_{\xi(G_{k}^{j}), i - 2 \leq \cdot \leq i + 2} F(u(\tau)) d\tau \|_{U_{\Delta}^{2}(G_{k}^{j} \times \mathbf{R}^{2})}^{2} \\ + \sum_{11 \leq i < j} 2^{i - j} \sum_{G_{\alpha}^{i} \subset G_{k}^{j}} \| \int_{t_{\alpha}^{i}}^{t} e^{i(t - \tau) \Delta} P_{\xi(G_{\alpha}^{i}), i - 2 \leq \cdot \leq i + 2} F(u(\tau)) d\tau \|_{U_{\Delta}^{2}(G_{\alpha}^{i} \times \mathbf{R}^{2})}^{2}.
\endaligned
\end{equation}

\noindent Similarly, by definition $\ref{d2.10}$, $(\ref{2.21})$, $(\ref{2.24})$, $(\ref{2.27})$, $(\ref{2.28})$, and $(\ref{3.9.0.2})$,

\begin{equation}\label{3.7.1}
\aligned
 \| u \|_{Y(G_{k}^{j} \times \mathbf{R}^{2})}^{2} \lesssim \epsilon_{2}^{3/2} + \sum_{i \geq j, i \geq 11; N(G_{k}^{j}) \leq 2^{i - 5} \epsilon_{3}^{1/2}} \| \int_{t_{\alpha}^{i}}^{t} e^{i(t - \tau) \Delta} P_{\xi(G_{k}^{j}), i - 2 \leq \cdot \leq i + 2} F(u(\tau)) d\tau \|_{U_{\Delta}^{2}(G_{k}^{j} \times \mathbf{R}^{2})}^{2}
\\ + \sum_{11 \leq i < j} 2^{i - j} \sum_{G_{\alpha}^{i} \subset G_{k}^{j}; N(G_{\alpha}^{i}) \leq 2^{i - 5} \epsilon_{3}^{1/2}} \| \int_{t_{\alpha}^{i}}^{t} e^{i(t - \tau) \Delta} P_{\xi(G_{\alpha}^{i}), i - 2 \leq \cdot \leq i + 2} F(u(\tau)) d\tau \|_{U_{\Delta}^{2}(G_{\alpha}^{i} \times \mathbf{R}^{2})}^{2}.
\endaligned
\end{equation}

\noindent First take the intervals $G_{\alpha}^{i} \subset G_{k}^{j}$ with $N(G_{\alpha}^{i}) \geq \epsilon_{3}^{1/2} 2^{i - 5}$. These intervals appear in $(\ref{3.7})$ but not $(\ref{3.7.1})$.\vspace{5mm}

\noindent There are at most two small intervals, call them $J_{1}$ and $J_{2}$, that intersect $G_{k}^{j}$ but are not contained in $G_{k}^{j}$. By theorem $\ref{t1.1}$ and conservation of mass, $\| u \|_{L_{t}^{3} L_{x}^{6}(J_{l} \times \mathbf{R}^{2})} \lesssim 1$. Therefore,

\begin{equation}\label{3.7.1.1}
\sum_{11 \leq i < j} 2^{i - j} \sum_{G_{\alpha}^{i} \subset G_{k}^{j}} \| F(u) \|_{L_{t}^{1} L_{x}^{2}(G_{\alpha}^{i} \cap (J_{1} \cup J_{2}) \times \mathbf{R}^{2})}^{2} \lesssim 1.
\end{equation}

\noindent Next observe that $(\ref{2.20})$ and $N(G_{\alpha}^{i}) \geq \epsilon_{3}^{1/2} 2^{i - 5}$ implies that $N(t) \geq \epsilon_{3}^{1/2} 2^{i - 6}$ for all $t \in G_{\alpha}^{i}$, so by $(\ref{3.7.1.1})$ and $l^{1} \subset l^{2}$,

\begin{equation}\label{3.7.2}
\aligned
\sum_{11 \leq i < j} 2^{i - j} \sum_{G_{\alpha}^{i} \subset G_{k}^{j}; N(G_{\alpha}^{i}) \geq \epsilon_{3}^{1/2} 2^{i - 5}} \| F(u) \|_{L_{t}^{1} L_{x}^{2}(G_{\alpha}^{i} \times \mathbf{R}^{2})}^{2} \\ \lesssim 1 + \sum_{11 \leq i < j} 2^{i - j} (\sum_{J_{l} \subset G_{k}^{j}: N(J_{l}) \geq 2^{i - 6} \epsilon_{3}^{1/2}} \| F(u) \|_{L_{t}^{1} L_{x}^{2}(J_{l} \times \mathbf{R}^{2})}^{2})
\endaligned
\end{equation}

\begin{equation}\label{3.7.3}
\lesssim 1 + \sum_{J_{l} \subset G_{k}^{j}} \sum_{11 \leq i < j; 2^{i} \leq 2^{6} \epsilon_{3}^{-1/2} N(J_{l})} 2^{i - j} \lesssim 2^{-j} \epsilon_{3}^{-1/2} \sum_{J_{l} \subset G_{k}^{j}} N(J_{l}) \lesssim 1.
\end{equation}

\noindent The last inequality follows from $(\ref{1.28})$ and $(\ref{2.24})$.\vspace{5mm}

\noindent Similarly, if $N(G_{k}^{j}) \geq 2^{j - 5} \epsilon_{3}^{1/2}$, $(\ref{2.20})$ implies that $N(t) \geq \epsilon_{3}^{1/2} 2^{j - 6}$ for all $t \in G_{k}^{j}$. Therefore by $(\ref{1.28})$ (which implies $\int_{G_{k}^{j}} N(t)^{3} dt \lesssim N(G_{k}^{j})$), $\int_{G_{k}^{j}} N(t)^{2} dt \lesssim 1$. Therefore, $(\ref{1.26})$ implies that

\begin{equation}\label{3.7.4}
\sum_{i \geq j; i \geq 11} \| P_{\xi(G_{k}^{j}), i - 2 \leq \cdot \leq i + 2} F(u) \|_{L_{t}^{1} L_{x}^{2}(G_{k}^{j} \times \mathbf{R}^{2})}^{2} \lesssim \| u \|_{L_{t}^{3} L_{x}^{6}(G_{k}^{j} \times \mathbf{R}^{2})}^{6} \lesssim 1.
\end{equation}

\noindent Next, suppose that $2^{i - 10} \epsilon_{3}^{1/2} \leq N(G_{\alpha}^{i}) \leq 2^{i - 5} \epsilon_{3}^{1/2}$. Once again let $J_{1}$ and $J_{2}$ be the two small intervals that intersect $G_{k}^{j}$ but are not contained in $G_{k}^{j}$. Then by $(\ref{2.20})$ and $(\ref{2.21})$,

\begin{equation}\label{3.7.5}
\aligned
\sum_{11 \leq i < j} 2^{i - j} \sum_{N(G_{\alpha}^{i}) \leq \epsilon_{3}^{1/2} 2^{i - 5}} \| P_{\xi(G_{\alpha}^{i}), i - 2 \leq \cdot \leq i + 2} F(u) \|_{L_{t,x}^{4/3}(G_{\alpha}^{i} \cap (J_{1} \cup J_{2}))}^{2} \\ \lesssim \sum_{11 \leq i < j} 2^{i - j} \| P_{\xi(t), \geq 4 \epsilon_{3}^{-1/2} N(t)} F(u) \|_{L_{t, x}^{4/3}(J_{1} \cup J_{2} \times \mathbf{R}^{2})}^{2} \\ \lesssim \| P_{\xi(t), \geq \epsilon_{3}^{-1/2} N(t)} u \|_{L_{t}^{\infty} L_{x}^{2}}^{2} \| u \|_{L_{t}^{8/3} L_{x}^{8}(J_{1} \cup J_{2} \times \mathbf{R}^{2})}^{4} \lesssim \epsilon_{2}^{2}.
\endaligned
\end{equation}

\noindent Next, as in $(\ref{3.7.3})$, $(\ref{2.20})$ implies $N(t) \geq 2^{i - 11} \epsilon_{3}^{1/2}$ on $G_{\alpha}^{i}$, so

\begin{equation}
\sum_{11 \leq i < j} 2^{i - j} (\sum_{J_{l} \subset G_{k}^{j}: N(J_{l}) \geq 2^{i - 11} \epsilon_{3}^{1/2}} \| F(u) \|_{L_{t}^{1} L_{x}^{2}(J_{l} \times \mathbf{R}^{2})}^{2})
\end{equation}

\begin{equation}\label{3.7.6}
\lesssim \sum_{J_{l} \subset G_{k}^{j}} \sum_{11 \leq i < j; 2^{i} \leq 2^{11} \epsilon_{3}^{-1/2} N(J_{l})} 2^{i - j} \lesssim 2^{-j} \epsilon_{3}^{-1/2} \sum_{J_{l} \subset G_{k}^{j}} N(J_{l}) \lesssim \epsilon_{3}^{1/2}.
\end{equation}

\noindent Finally, if $\epsilon_{3}^{1/2} 2^{j - 10} \leq N(G_{k}^{j}) \leq \epsilon_{3}^{1/2} 2^{j - 5}$, $N(t) \geq \epsilon_{2}^{1/2} 2^{j - 11}$, $\int_{G_{k}^{j}} N(t)^{2} dt \lesssim 1$, so

\begin{equation}
\sum_{i \geq j; i \geq 11} \| P_{\xi(G_{k}^{j}), i - 2 \leq \cdot \leq i + 2} F(u) \|_{L_{t}^{1} L_{x}^{2}(G_{k}^{j} \times \mathbf{R}^{2})}^{2} \lesssim \| P_{\xi(G_{k}^{j}), \geq j - 5} u \|_{L_{t}^{\infty} L_{x}^{2}(G_{k}^{j} \times \mathbf{R}^{2})}^{3/2} \| u \|_{L_{t}^{9/4} L_{x}^{18}(G_{k}^{j} \times \mathbf{R}^{2})}^{9/2} \lesssim \epsilon_{2}^{3/2}.
\end{equation}



\noindent Therefore,

\begin{equation}
 \aligned
(\ref{3.7}) \lesssim 1 + \sum_{11 \leq i \leq j} 2^{i - j} \sum_{G_{\alpha}^{i} \subset G_{k}^{j}; N(G_{\alpha}^{i}) \leq 2^{i - 10} \epsilon_{3}^{1/2}} \| \int_{t_{\alpha}^{i}}^{t} e^{i(t - \tau) \Delta} P_{\xi(G_{\alpha}^{i}), i - 2 \leq \cdot \leq i + 2} F(u(\tau)) d\tau \|_{U_{\Delta}^{2}(G_{\alpha}^{i} \times \mathbf{R}^{2})}^{2} \\
+ \sum_{i \geq j; i \geq 11; N(G_{k}^{j}) \leq 2^{i - 10} \epsilon_{3}^{1/2}} \| \int_{t_{\alpha}^{i}}^{t} e^{i(t - \tau) \Delta} P_{\xi(G_{k}^{j}), i - 2 \leq \cdot \leq i + 2} F(u(\tau)) d\tau \|_{U_{\Delta}^{2}(G_{k}^{j} \times \mathbf{R}^{2})}^{2},
\endaligned
\end{equation}

\noindent and

\begin{equation}
 \aligned
(\ref{3.7.1}) \lesssim \epsilon_{2}^{3/2} + \sum_{11 \leq i \leq j} 2^{i - j} \sum_{G_{\alpha}^{i} \subset G_{k}^{j}; N(G_{\alpha}^{i}) \leq 2^{i - 10} \epsilon_{3}^{1/2}} \| \int_{t_{\alpha}^{i}}^{t} e^{i(t - \tau) \Delta} P_{\xi(G_{\alpha}^{i}), i - 2 \leq \cdot \leq i + 2} F(u(\tau)) d\tau \|_{U_{\Delta}^{2}(G_{\alpha}^{i} \times \mathbf{R}^{2})}^{2} \\
+ \sum_{i \geq j; i \geq 11; N(G_{k}^{j}) \leq 2^{i - 10} \epsilon_{3}^{1/2}} \| \int_{t_{\alpha}^{i}}^{t} e^{i(t - \tau) \Delta} P_{\xi(G_{k}^{j}), i - 2 \leq \cdot \leq i + 2} F(u(\tau)) d\tau \|_{U_{\Delta}^{2}(G_{k}^{j} \times \mathbf{R}^{2})}^{2}.
\endaligned
\end{equation}

\noindent To estimate this term it suffices to prove

\begin{theorem}\label{t3.0}
\begin{equation}\label{3.9.1}
 \aligned
\sum_{i \geq j; N(G_{k}^{j}) \leq 2^{i - 5} \epsilon_{3}^{1/2}} \| \int_{t_{k}^{j}}^{t} e^{i(t - \tau) \Delta} P_{\xi(G_{k}^{j}), i - 2 \leq \cdot \leq i + 2} F(u(\tau)) d\tau \|_{U_{\Delta}^{2}(G_{k}^{j} \times \mathbf{R}^{2})}^{2}
\\ + \sum_{0 \leq i \leq j} 2^{i - j} \sum_{G_{\alpha}^{i} \subset G_{k}^{j}; N(G_{\alpha}^{i}) \leq 2^{i - 5} \epsilon_{3}^{1/2}} \| \int_{t_{\alpha}^{i}}^{t} e^{i(t - \tau) \Delta} P_{\xi(G_{\alpha}^{i}), i - 2 \leq \cdot \leq i + 2} F(u(\tau)) d\tau \|_{U_{\Delta}^{2}(G_{\alpha}^{i} \times \mathbf{R}^{2})}^{2} \\
\lesssim \epsilon_{2}^{1/3} \| u \|_{\tilde{X}_{j}([0, T] \times \mathbf{R}^{2})}^{5/3} \| u \|_{\tilde{Y}_{j}([0, T] \times \mathbf{R}^{2})}^{2} + \epsilon_{2}^{2} \| u \|_{\tilde{Y}_{j}([0, T] \times \mathbf{R}^{2})}^{2} + \| u \|_{\tilde{Y}_{j}([0, T] \times \mathbf{R}^{2})}^{4} (1 + \| u \|_{\tilde{X}_{j}([0, T] \times \mathbf{R}^{2})}^{8}).
\endaligned
\end{equation}
\end{theorem}

\noindent Indeed, make a bootstrap argument. Suppose

\begin{equation}\label{3.9.1.1}
\| u \|_{\tilde{X}_{k_{\ast}}([0, T] \times \mathbf{R}^{2})}^{2} \leq C_{0},
\end{equation}

\noindent and

\begin{equation}\label{3.9.1.2}
\| u \|_{\tilde{Y}_{k_{\ast}}([0, T] \times \mathbf{R}^{2})}^{2} \leq C(u) \epsilon_{2}^{3/2} \leq \epsilon_{2}.
\end{equation}

\noindent Then

\begin{equation}
\| u \|_{\tilde{X}_{k_{\ast} + 1}([0, T] \times \mathbf{R}^{2})}^{2} \leq 2 C_{0},
\end{equation}

\noindent and

\begin{equation}
\| u \|_{\tilde{Y}_{k_{\ast} + 1}([0, T] \times \mathbf{R}^{2})}^{2} \leq 2 \epsilon_{2}.
\end{equation}

\noindent Then by $(\ref{3.7})$, $(\ref{3.7.1})$, $(\ref{3.7.3})$, $(\ref{3.7.4})$, and $(\ref{3.9.1})$,

\begin{equation}
\| u \|_{\tilde{X}_{k_{\ast} + 1}([0, T] \times \mathbf{R}^{2})} \leq C(u)(1 + \epsilon_{2}^{2/3} (2C_{0})^{5/6} + \epsilon_{2}^{3/2} + \epsilon_{2} (1 + 2 C_{0})^{8}),
\end{equation}

\noindent and

\begin{equation}
\| u \|_{\tilde{Y}_{k_{\ast} + 1}([0, T] \times \mathbf{R}^{2})} \leq C(u)(\epsilon_{2}^{3/4} + \epsilon_{2}^{2/3} (2C_{0})^{5/6} + \epsilon_{2}^{3/2} + \epsilon_{2} (1 + 2 C_{0})^{8}).
\end{equation}

\noindent Taking $C_{0} = 2^{6} C(u)$, $\epsilon_{2} > 0$ sufficiently small implies that $(\ref{3.9.1.1})$ and $(\ref{3.9.1.2})$ hold for $k_{\ast} = 11$, and also closes the bootstrap, implying

\begin{equation}
\aligned
\| u \|_{\tilde{X}_{k_{\ast} + 1}([0, T] \times \mathbf{R}^{2})} \leq C_{0}, \\
\| u \|_{\tilde{Y}_{k_{\ast} + 1}([0, T] \times \mathbf{R}^{2})} \leq \epsilon_{2}^{1/2}.
\endaligned
\end{equation}

\noindent Theorem $\ref{t3.1}$ then follows by $(\ref{3.8})$, $(\ref{3.9})$, and induction on $k_{\ast}$. $\Box$\vspace{5mm}

\noindent \emph{Proof of theorem $\ref{t3.0}$:} By $(\ref{2.6})$,

\begin{equation}\label{3.9.2}
\aligned
P_{\xi(G_{\alpha}^{i}), i - 2 \leq \cdot \leq i + 2} F(u) = e^{ix \cdot \xi(G_{\alpha}^{i})} P_{i - 2 \leq \cdot \leq i + 2} e^{-ix \cdot \xi(G_{\alpha}^{i})} F(u) 
 = e^{ix \cdot \xi(G_{\alpha}^{i})} P_{i - 2 \leq \cdot \leq i + 2} F(e^{-ix \cdot \xi(G_{\alpha}^{i})} u).
 \endaligned
 \end{equation}

\noindent Since the nonlinearity is an algebraic, power - type nonlinearity,

\begin{equation}\label{3.9.3}
\aligned
e^{ix \cdot \xi(G_{\alpha}^{i})} P_{i - 2 \leq \cdot \leq i + 2} F(e^{-ix \cdot \xi(G_{\alpha}^{i})} u) = P_{\xi(G_{\alpha}^{i}), i - 2 \leq \cdot \leq i + 2} O((P_{\xi(G_{\alpha}^{i}), \geq i + 4} u)^{2} u) \\ 
+ P_{\xi(G_{\alpha}^{i}), i - 2 \leq \cdot \leq i + 2} O((P_{\xi(G_{\alpha}^{i}), i - 5 \leq \cdot \leq i + 5} u) u^{2}).
\endaligned
\end{equation}

\noindent Then by $(\ref{2.26})$, $(\ref{3.9.2})$, and $(\ref{3.9.3})$,

\begin{equation}\label{3.10}
 \| \int_{t_{\alpha}^{i}}^{t} e^{i(t - \tau) \Delta} P_{\xi(G_{\alpha}^{i}), i - 2 \leq \cdot \leq i + 2} F(u(\tau)) d\tau \|_{U_{\Delta}^{2}(G_{\alpha}^{i} \times \mathbf{R}^{2})} 
\end{equation}

\begin{equation}\label{3.11}
\lesssim \| \int_{t_{\alpha}^{i}}^{t} e^{i(t - \tau) \Delta} P_{\xi(G_{\alpha}^{i}), i - 2 \leq \cdot \leq i + 2} ((P_{\xi(G_{\alpha}^{i}), \geq i - 5} u)(P_{\xi(\tau), \geq i - 10} u) u) d\tau \|_{U_{\Delta}^{2}(G_{\alpha}^{i} \times \mathbf{R}^{2})}
\end{equation}

\begin{equation}\label{3.12}
 + \| \int_{t_{\alpha}^{i}}^{t} e^{i(t - \tau) \Delta} P_{\xi(G_{\alpha}^{i}), i - 2 \leq \cdot \leq i + 2} ((P_{\xi(G_{\alpha}^{i}), i - 5 \leq \cdot \leq i + 5} u)(P_{\xi(\tau), \leq i - 10} u)^{2}) d\tau \|_{U_{\Delta}^{2}(G_{\alpha}^{i} \times \mathbf{R}^{2})}.
\end{equation}

\noindent First take $(\ref{3.11})$. By $(\ref{2.21})$, $(\ref{2.37})$, $(\ref{2.38})$, and $N(G_{\alpha}^{i}) \leq \epsilon_{3}^{1/2} 2^{i - 5}$,

\begin{equation}\label{3.13}
\| P_{\xi(\tau), \geq i - 10} u \|_{L_{t}^{5/2} L_{x}^{10}(G_{\alpha}^{i} \times \mathbf{R}^{2})} \lesssim \| P_{\xi(\tau), \geq i - 10} u \|_{L_{t}^{\infty} L_{x}^{2}}^{1/6} \| P_{\xi(\tau), \geq i - 10} u \|_{L_{t}^{25/12} L_{x}^{50}}^{5/6} \lesssim \epsilon_{2}^{1/6} \| u \|_{\tilde{X}_{i}([0, T] \times \mathbf{R}^{2})}^{5/6}.
\end{equation}

\noindent Therefore,

\begin{theorem}\label{t3.2}
For a fixed $G_{k}^{j} \subset [0, T]$, $j > 11$,

\begin{equation}\label{3.14.1}
\aligned
\sum_{11 \leq i < j} 2^{i - j} \sum_{G_{\alpha}^{i} \subset G_{k}^{j}} \| \int_{t_{\alpha}^{i}}^{t} e^{i(t - \tau) \Delta} P_{\xi(G_{\alpha}^{i}), i - 2 \leq \cdot \leq i + 2} ((P_{\xi(G_{\alpha}^{i}), \geq i - 5} u)(P_{\xi(\tau), \geq i - 10} u) u) d\tau \|_{U_{\Delta}^{2}(G_{\alpha}^{i} \times \mathbf{R}^{2})}^{2} \\ + \sum_{i \geq j} \| \int_{t_{\alpha}^{i}}^{t} e^{i(t - \tau) \Delta} P_{\xi(G_{k}^{j}), i - 2 \leq \cdot \leq i + 2} ((P_{\xi(G_{k}^{j}), \geq i - 5} u)(P_{\xi(\tau), \geq i - 10} u) u) d\tau \|_{U_{\Delta}^{2}(G_{k}^{j} \times \mathbf{R}^{2})}^{2} \\ \lesssim \epsilon_{2}^{1/3} \| u \|_{\tilde{X}_{j}([0, T] \times \mathbf{R}^{2})}^{5/3} \| u \|_{\tilde{Y}_{j}([0, T] \times \mathbf{R}^{2})}^{2}.
\endaligned
\end{equation}
\end{theorem}

\noindent \emph{Proof:} By $(\ref{2.14})$ let $\hat{v}(t, \xi)$ be supported on $2^{i - 2} \leq |\xi - \xi(G_{\alpha}^{i})| \leq 2^{i + 2}$, $\| v \|_{V_{\Delta}^{2}(G_{\alpha}^{i} \times \mathbf{R}^{2})} = 1$. By $(\ref{2.37})$ and $(\ref{2.38})$,

\begin{equation}\label{3.15}
\int_{G_{\alpha}^{i}} \langle v, (P_{\xi(G_{\alpha}^{i}), \geq i - 5} u)(P_{\xi(t), \geq i - 10} u) u \rangle dt 
\end{equation}

\begin{equation}\label{3.16}
\lesssim \sum_{l \geq i - 5} \| v (P_{\xi(G_{\alpha}^{i}), l} u) \|_{L_{t}^{5/2} L_{x}^{5/3}}^{1/2} \| v \|_{L_{t}^{5/2} L_{x}^{10}}^{1/2} \| P_{\xi(G_{\alpha}^{i}), l} u \|_{L_{t}^{5/2} L_{x}^{10}}^{1/2} \| P_{\xi(t), \geq i - 10} u \|_{L_{t}^{5/2} L_{x}^{10}} \| u \|_{L_{t}^{\infty} L_{x}^{2}(G_{\alpha}^{i} \times \mathbf{R}^{2})}
\end{equation}

\begin{equation}\label{3.17}
 \lesssim \epsilon_{2}^{1/6} \| u \|_{\tilde{X}_{i}([0, T] \times \mathbf{R}^{2})}^{5/6} \sum_{l \geq i - 5} 2^{(i - l)/5} \| P_{\xi(G_{\alpha}^{i}), l} u \|_{U_{\Delta}^{2}(G_{\alpha}^{i} \times \mathbf{R}^{2})}.
\end{equation}

\noindent The last inequality follows from $V_{\Delta}^{2} \subset U_{\Delta}^{5/2}$ and lemma $\ref{l2.6.2}$.\vspace{5mm}

\noindent Now for any $0 \leq l \leq j$, $G_{k}^{j}$ overlaps $2^{j - l}$ intervals $G_{\beta}^{l}$ and for $0 \leq i \leq l$, each $G_{\beta}^{l}$ overlaps $2^{l - i}$ intervals $G_{\alpha}^{i}$. Additionally, each $G_{\alpha}^{i}$ is the subset of one $G_{\beta}^{l}$. Therefore,

\begin{equation}\label{3.18}
 \sum_{11 \leq i \leq j} 2^{i - j} \sum_{G_{\alpha}^{i} \subset G_{k}^{j}; N(G_{\alpha}^{i}) \leq \epsilon_{3}^{1/2} 2^{i - 10}} (\sum_{l \geq i - 5} 2^{(i - l)/5} \| P_{\xi(G_{\alpha}^{i}), l} u \|_{U_{\Delta}^{2}(G_{\alpha}^{i} \times \mathbf{R}^{2})}^{2})
\end{equation}

\begin{equation}\label{3.19}
 \lesssim \sum_{11 \leq i \leq j} 2^{i - j} \sum_{G_{\beta}^{l} \subset G_{k}^{j}; N(G_{\beta}^{l}) \leq \epsilon_{3}^{1/2} 2^{l -5}} (\sum_{i - 5 \leq l \leq i} 2^{(i - l)/5} \| P_{\xi(G_{\beta}^{l}), l - 2 \leq \cdot \leq l + 2} u \|_{U_{\Delta}^{2}(G_{\beta}^{l} \times \mathbf{R}^{2})}^{2})
\end{equation}

\begin{equation}\label{3.20}
 + \sum_{0 \leq l \leq j} 2^{l - j} \sum_{G_{\beta}^{l} \subset G_{k}^{j}; N(G_{\beta}^{l}) \leq \epsilon_{3}^{1/2} 2^{l - 5}} (\| P_{\xi(G_{\beta}^{l}), l - 2 \leq \cdot \leq l + 2} u \|_{U_{\Delta}^{2}(G_{\beta}^{l} \times \mathbf{R}^{2})}^{2} (\sum_{11 \leq i \leq l} 2^{(i - l)/5}))
\end{equation}

\begin{equation}\label{3.21}
 + \sum_{11 \leq i \leq j} \sum_{l > j; N(G_{k}^{j}) \leq \epsilon_{3}^{1/2} 2^{j - 5}} 2^{(i - l)/5} \| P_{\xi(G_{k}^{j}), l - 2 \leq \cdot \leq l + 2} u \|_{U_{\Delta}^{2}(G_{k}^{j} \times \mathbf{R}^{2})}^{2} \lesssim \| u \|_{\tilde{Y}_{j}([0, T] \times \mathbf{R}^{2})}^{2}.
\end{equation}

\noindent Also,

\begin{equation}\label{3.22}
 \sum_{i \geq j; N(G_{k}^{j}) \leq \epsilon_{3}^{1/2} 2^{i - 10}}  (\sum_{l \geq i - 5} 2^{(i - l)/5} \| P_{\xi(G_{k}^{j}), l} u \|_{U_{\Delta}^{2}(G_{k}^{j} \times \mathbf{R}^{2})}^{2})
\end{equation}

\begin{equation}\label{3.23}
 \lesssim \sum_{i \geq j} \sum_{l \geq i - 5, l < j} (\sum_{G_{\beta}^{l} \subset G_{k}^{j}; N(G_{\beta}^{l}) \leq \epsilon_{3}^{1/2} 2^{l - 5}} 2^{(i - l)/5} \| P_{\xi(G_{\beta}^{l}), l - 2 \leq \cdot \leq l + 2} u \|_{U_{\Delta}^{2}(G_{\beta}^{l} \times \mathbf{R}^{2})}^{2})
\end{equation}

\begin{equation}\label{3.24}
 + \sum_{i \geq j} (\sum_{l \geq i - 5, l \geq j; N(G_{k}^{j}) \leq \epsilon_{3}^{1/2} 2^{l - 5}} 2^{(i - l)/5} \| P_{\xi(G_{k}^{j}), l - 2 \leq \cdot \leq l + 2} u \|_{U_{\Delta}^{2}(G_{k}^{j} \times \mathbf{R}^{2})}^{2}) \lesssim \| u \|_{\tilde{Y}_{j}([0, T] \times \mathbf{R}^{2})}^{2}.
\end{equation}

\noindent $(\ref{3.17})$, together with $(\ref{3.18}) - (\ref{3.21})$, and $(\ref{3.22}) - (\ref{3.24})$ implies theorem $\ref{t3.2}$. $\Box$\vspace{5mm}

\noindent Now take $(\ref{3.12})$.

\begin{theorem}\label{tbigtheorem}
For any $11 \leq i < j$, $G_{\alpha}^{i} \subset G_{k}^{j}$,

\begin{equation}\label{big1}
\aligned
\| \int_{t_{\alpha}^{i}}^{t} e^{i(t - \tau) \Delta} P_{\xi(G_{\alpha}^{i}), i - 2 \leq \cdot \leq i + 2} ((P_{\xi(G_{\alpha}^{i}), i - 5 \leq \cdot \leq i + 5} u)(P_{\xi(\tau), \leq i - 10} u)^{2}) d\tau \|_{U_{\Delta}^{2}(G_{\alpha}^{i} \times \mathbf{R}^{2})} \\ \lesssim \| P_{\xi(G_{\alpha}^{i}), i - 5 \leq \cdot \leq i + 5} u \|_{U_{\Delta}^{2}(G_{\alpha}^{i} \times \mathbf{R}^{2})} (\epsilon_{2} + \| u \|_{\tilde{Y}_{i}(G_{\alpha}^{i} \times \mathbf{R}^{2})} (1 + \| u \|_{\tilde{X}_{i}(G_{\alpha}^{i} \times \mathbf{R}^{2})})^{8}).
\endaligned
\end{equation}

\noindent Also, for $i \geq j$, $j > 11$,

\begin{equation}\label{big2}
\aligned
\| \int_{t_{\alpha}^{i}}^{t} e^{i(t - \tau) \Delta} P_{\xi(G_{k}^{j}), i - 2 \leq \cdot \leq i + 2} ((P_{\xi(G_{k}^{j}), i - 5 \leq \cdot \leq i + 5} u)(P_{\xi(\tau), \leq i - 10} u)^{2}) d\tau \|_{U_{\Delta}^{2}(G_{\alpha}^{i} \times \mathbf{R}^{2})} \\ \lesssim 2^{3(j - i)/4} \| P_{\xi(G_{k}^{j}), i - 5 \leq \cdot \leq i + 5} u \|_{U_{\Delta}^{2}(G_{\alpha}^{i} \times \mathbf{R}^{2})} (\epsilon_{2} + \| u \|_{\tilde{Y}_{i}(G_{\alpha}^{i} \times \mathbf{R}^{2})} (1 + \| u \|_{\tilde{X}_{j}(G_{k}^{j} \times \mathbf{R}^{2})})^{8}).
\endaligned
\end{equation}
\end{theorem}

\noindent By theorem $\ref{t3.2}$ and $(\ref{3.18})$ - $(\ref{3.24})$, theorem $\ref{t3.0}$ follows from theorem $\ref{tbigtheorem}$. $\Box$\vspace{5mm}

\noindent \emph{Proof of theorem $\ref{tbigtheorem}$:} The proof of theorem $\ref{tbigtheorem}$ is considerably involved, and will occupy the remainder of this section, as well as the appendix. The proof will focus on $(\ref{big1})$, as the proof of $(\ref{big2})$ is nearly identical.\vspace{5mm}

\noindent For a given $G_{\alpha}^{i}$ there are at most two small intervals $J_{1}$, $J_{2}$ that overlap $G_{\alpha}^{i}$ but are not contained in $G_{\alpha}^{i}$. Let $\tilde{G}_{\alpha}^{i} = G_{\alpha}^{i} \setminus (J_{1} \cup J_{2})$. Then by lemma $\ref{l2.6}$, Duhamel's principle, $V_{\Delta}^{2} \subset U_{\Delta}^{4}$, and $U_{\Delta}^{2} \subset L_{t}^{\infty} L_{x}^{2}$,

\begin{equation}\label{3.28}
\| \int_{t_{\alpha}^{i}}^{t} e^{i(t - \tau) \Delta} P_{\xi(G_{\alpha}^{i}), i - 2 \leq \cdot \leq i + 2} ((P_{\xi(G_{\alpha}^{i}), i - 5 \leq \cdot \leq i + 5} u)(P_{\xi(\tau), \leq i - 10} u)^{2}) d\tau \|_{U_{\Delta}^{2}(G_{\alpha}^{i} \times \mathbf{R}^{2})}
\end{equation}

\begin{equation}\label{3.28.1}
\aligned
 \lesssim \| \int_{t_{\alpha}^{i}}^{t} e^{i(t - \tau) \Delta} P_{\xi(G_{\alpha}^{i}), i - 2 \leq \cdot \leq i + 2} ((P_{\xi(G_{\alpha}^{i}), i - 5 \leq \cdot \leq i + 5} u)(P_{\xi(\tau), \leq i - 10} u)^{2}) d\tau \|_{U_{\Delta}^{2}(\tilde{G}_{\alpha}^{i} \times \mathbf{R}^{2})} \\
+ \| (P_{\xi(G_{\alpha}^{i}), i - 5 \leq \cdot \leq i + 5} u)(P_{\xi(\tau), \leq i - 10} u)^{2} \|_{L_{t,x}^{4/3}(J_{1} \cap G_{\alpha}^{i} \times \mathbf{R}^{2})} \\ + \| (P_{\xi(G_{\alpha}^{i}), i - 5 \leq \cdot \leq i + 5} u)(P_{\xi(\tau), \leq i - 10} u)^{2} \|_{L_{t,x}^{4/3}(J_{2} \cap G_{\alpha}^{i} \times \mathbf{R}^{2})}.
\endaligned
\end{equation}

\noindent Now $t_{\alpha}^{i}$ need not lie in $\tilde{G}_{\alpha}^{i}$. However, if $t_{\alpha}^{i} \notin \tilde{G}_{\alpha}^{i}$ then we can move $t_{\alpha}^{i}$ into $\tilde{G}_{\alpha}^{i}$ at a cost of

\begin{equation}\label{3.28.2}
 \| (P_{\xi(G_{\alpha}^{i}), i - 5 \leq \cdot \leq i + 5} u)(P_{\xi(\tau), \leq i - 10} u)^{2} \|_{L_{t,x}^{4/3}(J_{1} \cap G_{\alpha}^{i} \times \mathbf{R}^{2})} + \| (P_{\xi(G_{\alpha}^{i}), i - 5 \leq \cdot \leq i + 5} u)(P_{\xi(\tau), \leq i - 10} u)^{2} \|_{L_{t,x}^{4/3}(J_{2} \cap G_{\alpha}^{i} \times \mathbf{R}^{2})}
\end{equation}

\noindent Indeed suppose without loss of generality that $t_{\alpha}^{i} \in J_{1}$ and let $\tilde{t}_{\alpha}^{i}$ be the left endpoint of $\tilde{G}_{\alpha}^{i}$. Then for any $t \in \tilde{G}_{\alpha}^{i}$

\begin{equation}
\aligned
\int_{t_{\alpha}^{i}}^{t} e^{i(t - \tau) \Delta} P_{\xi(G_{\alpha}^{i}), i - 2 \leq \cdot \leq i + 2} ((P_{\xi(G_{\alpha}^{i}), i - 5 \leq \cdot \leq i + 5} u)(P_{\xi(\tau), \leq i - 10} u)^{2}) d\tau \\ - \int_{\tilde{t}_{\alpha}^{i}}^{t} e^{i(t - \tau) \Delta} P_{\xi(G_{\alpha}^{i}), i - 2 \leq \cdot \leq i + 2} ((P_{\xi(G_{\alpha}^{i}), i - 5 \leq \cdot \leq i + 5} u)(P_{\xi(\tau), \leq i - 10} u)^{2}) d\tau \\ = e^{i(t - \tilde{t}_{\alpha}^{i}) \Delta} \int_{t_{\alpha}^{i}}^{\tilde{t}_{\alpha}^{i}} e^{i(\tilde{t}_{\alpha}^{i} - \tau) \Delta} P_{\xi(G_{\alpha}^{i}), i - 2 \leq \cdot \leq i + 2} ((P_{\xi(G_{\alpha}^{i}), i - 5 \leq \cdot \leq i + 5} u)(P_{\xi(\tau), \leq i - 10} u)^{2}) d\tau.
\endaligned
\end{equation}

\noindent By theorem $\ref{t1.1}$,

\begin{equation}
\aligned
\| \int_{t_{\alpha}^{i}}^{\tilde{t}_{\alpha}^{i}} e^{i(\tilde{t}_{\alpha}^{i} - \tau) \Delta} P_{\xi(G_{\alpha}^{i}), i - 2 \leq \cdot \leq i + 2} ((P_{\xi(G_{\alpha}^{i}), i - 5 \leq \cdot \leq i + 5} u)(P_{\xi(\tau), \leq i - 10} u)^{2}) d\tau \|_{L_{x}^{2}(\mathbf{R}^{2})} \\ \lesssim \| (P_{\xi(G_{\alpha}^{i}), i - 5 \leq \cdot \leq i + 5} u)(P_{\xi(\tau), \leq i - 10} u)^{2} \|_{L_{t,x}^{4/3}(J_{1} \cap G_{\alpha}^{i} \times \mathbf{R}^{2})}.
\endaligned
\end{equation}

\noindent Now then, by the bilinear estimate $(\ref{1.3})$, $N(t) \leq 2^{i - 5} \epsilon_{3}^{1/2}$ on $G_{\alpha}^{i}$, $(\ref{2.20})$, and $(\ref{2.21})$,

\begin{equation}\label{3.28.3}
\aligned
\| (P_{\xi(G_{\alpha}^{i}), i - 5 \leq \cdot \leq i + 5} u)(P_{\xi(\tau), \leq i - 10} u)^{2} \|_{L_{t,x}^{4/3}(J_{1} \cap G_{\alpha}^{i} \times \mathbf{R}^{2})} \\ \lesssim \| (P_{\xi(G_{\alpha}^{i}), i - 5 \leq \cdot \leq i + 5} u)(P_{\xi(J_{1}), \leq \epsilon_{3}^{-1/4} N(J_{1})} u) \|_{L_{t,x}^{2}(J_{1} \cap G_{\alpha}^{i} \times \mathbf{R}^{2})} \| u \|_{L_{t,x}^{4}(J_{1} \cap G_{\alpha}^{i} \times \mathbf{R}^{2})} \\
+ \| (P_{\xi(G_{\alpha}^{i}), i - 5 \leq \cdot \leq i + 5} u) \|_{L_{t}^{8/3} L_{x}^{8}(J_{1} \cap G_{\alpha}^{i} \times \mathbf{R}^{2})} \| u \|_{L_{t}^{8/3} L_{x}^{8}(J_{1} \cap G_{\alpha}^{i} \times \mathbf{R}^{2})} \| P_{\xi(J_{1}), \geq \epsilon_{3}^{-1/4} N(J_{1})} u \|_{L_{t}^{\infty} L_{x}^{2}(J_{1} \cap G_{\alpha}^{i} \times \mathbf{R}^{2})} \\
\lesssim \epsilon_{2} \| P_{\xi(G_{\alpha}^{i}), i - 5 \leq \cdot \leq i + 5} u \|_{U_{\Delta}^{2}(G_{\alpha}^{i} \times \mathbf{R}^{2})}.
\endaligned
\end{equation}

\noindent Therefore, at the price of $(\ref{3.28.3})$ we have now simplified to a situation in which $G_{\alpha}^{i}$ is the union of a bunch of small intervals.\vspace{5mm}

\noindent Now by lemma $\ref{l2.6.1}$,

\begin{equation}\label{3.28.4}
\| \int_{t_{\alpha}^{i}}^{t} e^{i(t - \tau) \Delta} P_{\xi(G_{\alpha}^{i}), i - 2 \leq \cdot \leq i + 2} ((P_{\xi(G_{\alpha}^{i}), i - 5 \leq \cdot \leq i + 5} u)(P_{\xi(\tau), \leq i - 10} u)^{2}) d\tau \|_{U_{\Delta}^{2}(G_{\alpha}^{i} \times \mathbf{R}^{2})}
\end{equation}

\begin{equation}\label{3.29.2}
\aligned
\lesssim \sum_{0 \leq l_{2} \leq i - 10} (\sum_{J_{l} \subset G_{\alpha}^{i} ; N(J_{l}) \geq \epsilon_{3}^{1/2} 2^{l_{2} - 5}} &\| P_{\xi(G_{\alpha}^{i}), i - 2 \leq \cdot \leq i + 2} \\ &((P_{\xi(G_{\alpha}^{i}), i - 5 \leq \cdot \leq i + 5} u)(P_{\xi(t), l_{2}} u)(P_{\xi(t), \leq l_{2}} u)) \|_{V_{\Delta}^{2 \ast}(J_{l} \times \mathbf{R}^{2})}^{2})^{1/2}
\endaligned
\end{equation}

\begin{equation}\label{3.29.1}
\aligned
+ \sum_{0 \leq l_{2} \leq i - 10} \sum_{J_{l} \subset G_{\alpha}^{i} ; N(J_{l}) \geq \epsilon_{3}^{1/2} 2^{l_{2} - 5}} &\| \int_{J_{l}} e^{-it \Delta} P_{\xi(G_{\alpha}^{i}), i - 2 \leq \cdot \leq i + 2} \\ &((P_{\xi(G_{\alpha}^{i}), i - 5 \leq \cdot \leq i + 5} u)(P_{\xi(t), l_{2}} u)(P_{\xi(t), \leq l_{2}} u)) dt \|_{L_{x}^{2}(\mathbf{R}^{2})}
\endaligned
\end{equation}

\begin{equation}\label{3.30}
\aligned
+ \sum_{0 \leq l_{2} \leq i - 10} &(\sum_{G_{\beta}^{l_{2}} \subset G_{\alpha}^{i}; N(G_{\beta}^{l_{2}}) \leq \epsilon_{3}^{1/2} 2^{l_{2} - 5}} \\ &\| P_{\xi(G_{\alpha}^{i}), i - 2 \leq \cdot \leq i + 2} ((P_{\xi(G_{\alpha}^{i}), i - 5 \leq \cdot \leq i + 5} u)(P_{\xi(t), l_{2}} u)(P_{\xi(t), \leq l_{2}} u)) \|_{V_{\Delta}^{2 \ast}(G_{\beta}^{l_{2}} \times \mathbf{R}^{2})}^{2})^{1/2}
\endaligned
\end{equation}

\begin{equation}\label{3.29}
\aligned
+ \sum_{0 \leq l_{2} \leq i - 10} &\sum_{G_{\beta}^{l_{2}} \subset G_{\alpha}^{i}; N(G_{\beta}^{l_{2}}) \leq \epsilon_{3}^{1/2} 2^{l_{2} - 5}} \\ &\| \int_{G_{\beta}^{l_{2}}} e^{-it \Delta} P_{\xi(G_{\alpha}^{i}), i - 2 \leq \cdot \leq i + 2} ((P_{\xi(G_{\alpha}^{i}), i - 5 \leq \cdot \leq i + 5} u)(P_{\xi(t), l_{2}} u)(P_{\xi(t), \leq l_{2}} u)) dt \|_{L_{x}^{2}(\mathbf{R}^{2})}.
\endaligned
\end{equation}

\noindent \textbf{Remark:} It is possible that an interval $J_{l}$ with $N(J_{l}) \geq \epsilon_{3}^{1/2} 2^{l_{2} - 5}$ may intersect an interval $G_{\beta}^{l_{2}}$ satisfying $N(G_{\beta}^{l_{2}}) \leq \epsilon_{3}^{1/2} 2^{i - 5}$. However, since all our computations utilize Lebesgue norms $L_{t}^{p} L_{x}^{q}$, this technicality may safely be ignored.\vspace{5mm}

\noindent First take $(\ref{3.29.2})$. By proposition $\ref{p1.4}$ and $(\ref{2.20})$ (which implies $|\xi(t) - \xi(G_{\alpha}^{i})| << 2^{i}$),

\begin{equation}\label{11.1}
\| (P_{\xi(G_{\alpha}^{i}), i - 5 \leq \cdot \leq i + 5} u)(P_{\xi(t), \leq l_{2}} u) \|_{L_{t}^{2} L_{x}^{5/3}(J_{l} \times \mathbf{R}^{2})} \lesssim 2^{-i/5} \| P_{\xi(G_{\alpha}^{i}), i - 5 \leq \cdot \leq i + 5} u \|_{U_{\Delta}^{2}(G_{\alpha}^{i} \times \mathbf{R}^{2})}.
\end{equation}

\noindent For $\hat{v}(t, \xi)$ supported on $2^{i - 2} \leq |\xi - \xi(G_{\alpha}^{i})| \leq 2^{i + 2}$, $\| v \|_{V_{\Delta}^{2}(J_{l} \times \mathbf{R}^{2})} = 1$, the Sobolev embedding theorem and $V_{\Delta}^{2} \subset U_{\Delta}^{2 + \epsilon}$ imply,

\begin{equation}\label{11.2}
\aligned
\| (P_{\xi(G_{\alpha}^{i}), i - 5 \leq \cdot \leq i + 5} u)(P_{\xi(t), \leq l_{2}} u)(P_{\xi(t), l_{2}} u) v \|_{L_{t,x}^{1}(J_{l} \times \mathbf{R}^{2})} \\ \lesssim 2^{-i/5} \| P_{\xi(G_{\alpha}^{i}), i - 5 \leq \cdot \leq i + 5} u \|_{U_{\Delta}^{2}(G_{\alpha}^{i} \times \mathbf{R}^{2})} \| P_{\xi(t), l_{2}} u \|_{L_{t}^{4} L_{x}^{20/3}(J_{l} \times \mathbf{R}^{2})} \| v \|_{L_{t,x}^{4}(J_{l} \times \mathbf{R}^{2})} \\ \lesssim 2^{(l_{2} - i)/5} \| P_{\xi(G_{\alpha}^{i}), i - 5 \leq \cdot \leq i + 5} u \|_{U_{\Delta}^{2}(G_{\alpha}^{i} \times \mathbf{R}^{2})}.
\endaligned
\end{equation}

\noindent Also, by $(\ref{1.3})$, $(\ref{2.20})$, $(\ref{1.12.1})$ - $(\ref{1.14})$, and $N(J_{l}) \geq \epsilon_{3}^{1/2} 2^{l_{2} - 5}$,

\begin{equation}\label{11.3}
\aligned
\| (P_{\xi(G_{\alpha}^{i}), i - 5 \leq \cdot \leq i + 5} u)(P_{\xi(t), \leq l_{2}} u) \|_{L_{t,x}^{2}(J_{l} \times \mathbf{R}^{2})} \lesssim \| (P_{\xi(G_{\alpha}^{i}), i - 5 \leq \cdot \leq i + 5} u)(P_{\xi(J_{l}), \leq 32 \epsilon_{3}^{-1/2} N(J_{l})} u) \|_{L_{t,x}^{2}(J_{l} \times \mathbf{R}^{2})} \\ \lesssim \epsilon_{3}^{-1/4} 2^{-i/2} N(J_{l})^{1/2} \| P_{\xi(G_{\alpha}^{i}), i - 5 \leq \cdot \leq i + 5} u \|_{U_{\Delta}^{2}(G_{\alpha}^{i} \times \mathbf{R}^{2})},
\endaligned
\end{equation}

\noindent so by corollary $\ref{c1.3}$,

\begin{equation}\label{11.4}
\aligned
\| (P_{\xi(G_{\alpha}^{i}), i - 5 \leq \cdot \leq i + 5} u)(P_{\xi(t), \leq l_{2}} u)(P_{\xi(t), l_{2}} u)v \|_{L_{t,x}^{1}(J_{l} \times \mathbf{R}^{2})} \\ \lesssim \epsilon_{3}^{-1/4} 2^{-i/2} N(J_{l})^{1/2} \| P_{\xi(G_{\alpha}^{i}), i - 5 \leq \cdot \leq i + 5} u \|_{U_{\Delta}^{2}(G_{\alpha}^{i} \times \mathbf{R}^{2})} \| v(P_{\xi(t), l_{2}} u) \|_{L_{t}^{20/9} L_{x}^{20/11}(J_{l} \times \mathbf{R}^{2})}^{5/6} \\ \times \| v \|_{L_{t}^{8/3} L_{x}^{8}(J_{l} \times \mathbf{R}^{2})}^{1/6} \| P_{\xi(t), l_{2}} u \|_{L_{t}^{8/3} L_{x}^{8}(J_{l} \times \mathbf{R}^{2})}^{1/6} \\ \lesssim \epsilon_{3}^{-7/16} 2^{-7i/8} N(J_{l})^{7/8} \| P_{\xi(G_{\alpha}^{i}), i - 5 \leq \cdot \leq i + 5} u \|_{U_{\Delta}^{2}(G_{\alpha}^{i} \times \mathbf{R}^{2})}.
\endaligned
\end{equation} 

\noindent Therefore,

\begin{equation}\label{11.5}
\aligned
\| (P_{\xi(G_{\alpha}^{i}), i - 5 \leq \cdot \leq i + 5} u)(P_{\xi(t), \leq l_{2}} u)(P_{\xi(t), l_{2}} u)v \|_{L_{t,x}^{1}(J_{l} \times \mathbf{R}^{2})}^{2} \\ \lesssim 2^{-i} N(J_{l}) 2^{(l_{2} - i) \cdot \frac{6}{35}} \| P_{\xi(G_{\alpha}^{i}), i - 5 \leq \cdot \leq i + 5} u \|_{U_{\Delta}^{2}(G_{\alpha}^{i} \times \mathbf{R}^{2})}.
\endaligned
\end{equation}

\noindent Summing up,

\begin{equation}\label{11.6}
\aligned
(\ref{3.29.2}) \lesssim \sum_{0 \leq l_{2} \leq i - 10} (\sum_{J_{l} \subset G_{\alpha}^{i}} \epsilon_{3}^{-1/2} 2^{-i} N(J_{l}) 2^{\frac{6}{35} (l_{2} - i)} \| P_{\xi(G_{\alpha}^{i}), i - 2 \leq \cdot \leq i + 2} u \|_{U_{\Delta}^{2}(G_{\alpha}^{i} \times \mathbf{R}^{2})})^{1/2} \\ \lesssim \epsilon_{3}^{1/4} \| P_{\xi(G_{\alpha}^{i}), i - 5 \leq \cdot \leq i + 5} u \|_{U_{\Delta}^{2}(G_{\alpha}^{i} \times \mathbf{R}^{2})} \lesssim \epsilon_{2}^{2}  \| P_{\xi(G_{\alpha}^{i}), i - 5 \leq \cdot \leq i + 5} u \|_{U_{\Delta}^{2}(G_{\alpha}^{i} \times \mathbf{R}^{2})}
\endaligned
\end{equation}

\noindent For $(\ref{3.29.1}) - (\ref{3.29})$ we will prove three bilinear estimates that rely on the interaction Morawetz estimates of \cite{PV}. Such estimates will give a logarithmic improvement over what would be obtained from $(\ref{1.3})$ directly. This improvement is quite helpful to the proof.

\begin{theorem}[First bilinear Strichartz estimate]\label{t3.2.0}
Suppose $v_{0} \in L^{2}(\mathbf{R}^{2})$ has Fourier transform $\hat{v}_{0}(\xi)$ supported on $2^{i - 5} \leq |\xi - \xi(G_{\alpha}^{i})| \leq 2^{i + 5}$. Also suppose $J_{l} \subset G_{\alpha}^{i}$ is a small interval and $|\xi(t) - \xi(G_{\alpha}^{i})| \leq 2^{i - 10}$ for all $t \in G_{\alpha}^{i}$. Then for any $0 \leq l_{2} \leq i - 10$,

\begin{equation}\label{12.1}
\aligned
\| (e^{it \Delta} v_{0})(P_{\xi(t), \leq l_{2}} u) \|_{L_{t,x}^{2}(J_{l} \times \mathbf{R}^{2})}^{2}  \lesssim 2^{l_{2} - i} \| v_{0} \|_{L^{2}(\mathbf{R}^{2})}  \\ + 2^{-i} \| v_{0} \|_{L^{2}(\mathbf{R}^{2})}^{2} (\int_{J_{l}} |\xi'(t)| \sum_{l_{1} \leq l_{2}} 2^{(l_{1} - l_{2})} \| P_{\xi(t), l_{1}} u(t) \|_{L_{x}^{2}(\mathbf{R}^{2})} \| P_{\xi(t), l_{2} - 3 \leq \cdot \leq l_{2} + 3} u(t) \|_{L_{x}^{2}(\mathbf{R}^{2})} dt).
\endaligned
\end{equation}

\noindent The same estimate also holds when $P_{\xi(t), \leq l_{2}}$ is replaced by $P_{\xi(t), l_{2}}$.
\end{theorem}

\noindent \emph{Proof:} Let

\begin{equation}\label{12.2}
v = e^{it \Delta} v_{0},
\end{equation}

\noindent and

\begin{equation}\label{12.3}
w = P_{\xi(t), \leq l_{2}} u.
\end{equation}

\noindent Then $v$ and $w$ solve the equations

\begin{equation}\label{12.4}
i v_{t} + \Delta v = 0,
\end{equation}

\noindent and

\begin{equation}\label{12.5}
i w_{t} + \Delta w = F(w) + \mathcal N_{1} + \mathcal N_{2} = F(w) + \mathcal N,
\end{equation}

\noindent where

\begin{equation}\label{12.6}
\mathcal N_{1} = P_{\xi(t), \leq l_{2}} F(u) - F(w),
\end{equation}

\noindent and

\begin{equation}\label{12.7}
\mathcal N_{2} = (\frac{d}{dt} P_{\xi(t), \leq l_{2}}) u.
\end{equation}

\noindent $\frac{d}{dt} P_{\xi(t), \leq l_{2}}$ is the Fourier multiplier

\begin{equation}\label{12.8}
-\nabla \phi(\frac{\xi - \xi(t)}{2^{l_{2}}}) \cdot \frac{\xi'(t)}{2^{l_{2}}}.
\end{equation}

\noindent Following \cite{PV}, let

\begin{equation}\label{12.10}
M_{\omega}(t) = \int \int |w(t,y)|^{2} \frac{(x - y)_{\omega}}{|(x - y)_{\omega}|} Im[\bar{v} \partial_{\omega} v](t,x) dx dy + \int \int |v(t,y)|^{2} \frac{(x - y)_{\omega}}{|(x - y)_{\omega}|} Im[\bar{w} \partial_{\omega} w](t,x) dx dy,
\end{equation}

\noindent and

\begin{equation}\label{12.9}
M(t) = \int_{\omega \in S^{1}} M_{\omega}(t) d\omega.
\end{equation}

\noindent Integrating by parts,

\begin{equation}\label{12.11}
\frac{d}{dt} M_{\omega}(t) = 2 \int \int_{x_{\omega} = y_{\omega}} |w(t,y)|^{2} |\partial_{\omega} v(t,x)|^{2} +  \int \int_{x_{\omega} = y_{\omega}} \partial_{\omega}(|w(t,y)|^{2}) \partial_{\omega}(|v(t,x)|^{2}) dx dy
\end{equation}

\begin{equation}\label{12.12}
+ 2 \int \int_{x_{\omega} = y_{\omega}} |v(t,y)|^{2} |\partial_{\omega} w(t,x)|^{2} dx dy - 4 \int \int_{x_{\omega} = y_{\omega}} Im[\bar{w} \partial_{\omega} w](t,y) Im[\bar{v} \partial_{\omega} v](t,x) dx dy
\end{equation}

\begin{equation}\label{12.14}
+ \frac{1}{2} \int \int_{x_{\omega} = y_{\omega}} |v(t,y)|^{2} |w(t,x)|^{4} dx dy + 2 \int \int Im[\bar{w} \mathcal N](t,y) \frac{(x - y)_{\omega}}{|(x - y)_{\omega}|} Im[\bar{v} \partial_{\omega} v](t,x) dx dy
\end{equation}

\begin{equation}\label{12.13}
+ \int \int |v(t,y)|^{2} \frac{(x - y)_{\omega}}{|(x - y)_{\omega}|} Im[\bar{\mathcal N} \partial_{\omega} w](t,x) dx dy + \int \int |v(t,y)|^{2} \frac{(x - y)_{\omega}}{|(x - y)_{\omega}|} Im[\bar{w} \partial_{\omega} \mathcal N](t,x) dx dy,
\end{equation}

\noindent which is equal to

\begin{equation}\label{12.15}
2 \int \int_{x_{\omega} = y_{\omega}} |\partial_{\omega}(\overline{w(t,y)} v(t,x))|^{2} dx dy + \frac{1}{2} \int \int_{x_{\omega} = y_{\omega}} |v(t,y)|^{2} |w(t,x)|^{4} dx dy
\end{equation}

\begin{equation}\label{12.16}
+ 2 \int \int Im[\bar{w} \mathcal N](t,y) \frac{(x - y)_{\omega}}{|(x - y)_{\omega}|} Im[\bar{v} \partial_{\omega} v](t,x) dx dy
\end{equation}

\begin{equation}\label{12.17}
+ \int \int |v(t,y)|^{2} \frac{(x - y)_{\omega}}{|(x - y)_{\omega}|} Im[\bar{\mathcal N} \partial_{\omega} w](t,x) dx dy + \int \int |v(t,y)|^{2} \frac{(x - y)_{\omega}}{|(x - y)_{\omega}|} Im[\bar{w} \partial_{\omega} \mathcal N](t,x) dx dy.
\end{equation}

\noindent For a moment suppose $\omega = (1, 0)$. $|w(t,x)|^{2}$ has Fourier support on $|\xi| \leq 2^{l_{2} + 5}$, so if $\phi(\xi) = 1$ on $|\xi| \leq 1$, $\phi \in C_{0}^{\infty}(\mathbf{R})$ is a real - valued, even function, then as in $(\ref{2.3})$,



\begin{equation}\label{12.18}
|w(t,x_{1}, x_{2})|^{2} = 2^{l_{2} + 5} \int \check{\phi}(2^{l_{2} + 5}(x_{2} - y_{2})) |w(t, x_{1}, y_{2})|^{2} dy_{2},
\end{equation}

\noindent where $\check{\phi}$ is the inverse Fourier transform of $\phi$. Similarly,

\begin{equation}\label{12.19}
|\partial_{1} w(t, x_{1}, x_{2})|^{2} = 2^{l_{2} + 5} \int \check{\phi}(2^{l_{2} + 5}(x_{2} - y_{2})) |w(t, x_{1}, y_{2})|^{2} dy_{2},
\end{equation}

\noindent and

\begin{equation}\label{12.20}
Re[(\partial_{1} w(t, x_{1}, x_{2})) \overline{w(t, x_{1}, x_{2})}] = 2^{l_{2} + 5} \int \check{\phi}(2^{l_{2} + 5}(x_{2} - y_{2})) Re[(\partial_{1} w(t, x_{1}, y_{2})) \overline{w(t, x_{1}, y_{2})}] dy_{2}.
\end{equation}

\noindent Therefore,

\begin{equation}\label{12.21}
\int |\partial_{1}(\overline{w(t,x)} v(t,x))|^{2} dx \lesssim 2^{l_{2}} \int \int_{x_{1} = y_{1}} |\partial_{1}(\overline{w(t,y)} v(t,x))|^{2} dx dy,
\end{equation}

\noindent and

\begin{equation}\label{12.22}
\int |\nabla (\overline{w(t,x)} v(t,x))|^{2} dx \lesssim 2^{l_{2}} \int_{\omega \in S^{1}} \int \int_{x_{\omega} = y_{\omega}} |\partial_{\omega} (\overline{w(t,y)} v(t,x))|^{2} dx dy d\omega.
\end{equation}

\noindent By the fundamental theorem of calculus and $(\ref{12.15})$ - $(\ref{12.17})$,

\begin{equation}\label{12.23}
2^{-l_{2}} \int_{J_{l}} \int |\nabla (\overline{w(t,x)} v(t,x))|^{2} dx dt \lesssim \sup_{t \in J_{l}} |M(t)| 
\end{equation}

\begin{equation}\label{12.24}
- \int_{S^{1}} \int_{J_{l}} \int \int_{x_{\omega} = y_{\omega}} |v(t,y)|^{2} |w(t,x)|^{4} dx dy dt d\omega
\end{equation}

\begin{equation}\label{12.26}
-  \int_{S^{1}} \int_{J_{l}} \int \int |v(t,y)|^{2} \frac{(x - y)_{\omega}}{|(x - y)_{\omega}|} Im[\bar{\mathcal N} \partial_{\omega} w](t,x) dx dy dt d\omega
\end{equation}

\begin{equation}\label{12.27}
- \int_{S^{1}} \int_{J_{l}} \int \int |v(t,y)|^{2} \frac{(x - y)_{\omega}}{|(x - y)_{\omega}|} Im[\bar{w} \partial_{\omega} \mathcal N](t,x) dx dy dt d\omega.
\end{equation}

\begin{equation}\label{12.25}
- 2  \int_{S^{1}} \int_{J_{l}} \int \int Im[\bar{w} \mathcal N](t,y) \frac{(x - y)_{\omega}}{|(x - y)_{\omega}|} Im[\bar{v} \partial_{\omega} v](t,x) dx dy dt d\omega.
\end{equation}

\noindent Since $\partial_{\omega} = \nabla \cdot \omega = \cos(\omega) \partial_{1} + \sin(\omega) \partial_{2}$, there exists a constant $C$ such that

\begin{equation}\label{12.28}
\int_{S^{1}} \frac{x_{\omega}}{|x_{\omega}|} (\nabla \cdot \omega) d\omega = C \frac{x}{|x|} \cdot \nabla.
\end{equation}

\noindent Therefore,

\begin{equation}\label{12.29}
M(t) = C \int \int |w(t,y)|^{2} \frac{(x - y)}{|x - y|} \cdot Im[\bar{v} \nabla v](t,x) dx dy + C \int \int |v(t,y)|^{2} \frac{(x - y)}{|x - y|} \cdot Im[\bar{w} \nabla w](t,x) dx dy.
\end{equation}

\noindent Furthermore, since $\frac{(x - y)}{|x - y|}$ is an odd function of $x - y$, for any $t \in J_{l}$,

\begin{equation}\label{12.30}
\aligned
M(t) =  \int \int |w(t,y)|^{2} \frac{(x - y)}{|x - y|} \cdot Im[\bar{v} (\nabla - i \xi(t)) v](t,x) dx dy \\ + \int \int |v(t,y)|^{2} \frac{(x - y)}{|x - y|} \cdot Im[\bar{w} (\nabla - i \xi(t)) w](t,x) dx dy \lesssim 2^{i} \| w \|_{L_{t}^{\infty} L_{x}^{2}}^{2} \| v_{0} \|_{L_{x}^{2}}^{2} \lesssim 2^{i} \| v_{0} \|_{L^{2}}^{2}.
\endaligned
\end{equation}

\noindent \textbf{Remark:} The Galilean invariance of the interaction Morawetz estimate was observed and utilized by \cite{PV}.\vspace{5mm}

\noindent Next take $(\ref{12.24})$. Since this term is positive definite, we could simply ignore this term when bounding $2^{-l_{2}} \int \int |\nabla (v(t,x) \bar{w}(t,x))|^{2} dx dt$. However, the results of this paper will also be used in an upcoming paper (\cite{D5}), addressing the focusing equation

\begin{equation}\label{12.33}
(i \partial_{t} + \Delta) u = -|u|^{2} u.
\end{equation}

\noindent In the focusing problem the signs of $(\ref{10.24})$ - $(\ref{10.25})$ are positive instead of negative. Therefore, we will go ahead and do the computation here.

\begin{equation}\label{12.31}
\int_{S^{1}} \delta(\cos(\theta) r) d\theta = \frac{C'}{r},
\end{equation}

\noindent so

\begin{equation}\label{12.32}
\int \int \int_{x_{\omega} = y_{\omega}} |v(t,y)|^{2} |w(t,x)|^{4} dx dy d\omega = \int \int \frac{C'}{|x - y|} |v(t,y)|^{2} |w(t,x)|^{4} dx dy.
\end{equation}

\noindent By the Hardy - Littlewood - Sobolev inequality,

\begin{equation}\label{12.34}
\aligned
\int \int \int \frac{1}{|x - y|} |v(t,y)|^{2} |w(t,x)|^{4} dx dy dt \\ \lesssim \| v \|_{L_{t}^{6} L_{x}^{3}(J_{l} \times \mathbf{R}^{2})}^{2} \| w \|_{L_{t}^{6} L_{x}^{24/5}(J_{l} \times \mathbf{R}^{2})}^{4} \lesssim \| v_{0} \|_{L^{2}(\mathbf{R}^{2})}^{2} \| w \|_{L_{t}^{6} L_{x}^{24/5}(J_{l} \times \mathbf{R}^{2})}^{4}.
\endaligned
\end{equation}

\noindent The Sobolev embedding theorem and conservation of mass implies

\begin{equation}\label{12.35}
\| w \|_{L_{t}^{6} L_{x}^{24/5}(J_{l} \times \mathbf{R}^{2})} \lesssim 2^{l_{2}/4} \| P_{\xi(t), \leq l_{2}} u \|_{L_{t}^{6} L_{x}^{3}(J_{l} \times \mathbf{R}^{2})} \lesssim 2^{l_{2}/4} \| u \|_{L_{t,x}^{4}(J_{l} \times \mathbf{R}^{2})}^{2/3} \| u \|_{L_{t}^{\infty} L_{x}^{2}(J_{l} \times \mathbf{R}^{2})}^{1/3} \lesssim 2^{l_{2}/4}.
\end{equation}

\noindent Therefore, $(\ref{12.34}) \lesssim 2^{l_{2}/4}$.\vspace{5mm}

\noindent Next, applying $(\ref{12.28})$ and $(\ref{12.29})$ to $(\ref{12.26})$ - $(\ref{12.25})$, by direct computation,

\begin{equation}\label{12.36}
(\ref{12.26}) + (\ref{12.24}) + (\ref{12.25}) = C \int \int \int |v(t,y)|^{2} \frac{(x - y)}{|x - y|} \cdot Im[\bar{\mathcal N} (\nabla - i \xi(t)) w](t,x) dx dy dt,
\end{equation}

\begin{equation}\label{12.37}
+ C \int \int \int |v(t,y)|^{2} \frac{(x - y)}{|x - y|} \cdot Im[\bar{w} (\nabla - i \xi(t)) \mathcal N](t,x) dx dy dt,
\end{equation}

\begin{equation}\label{12.38}
+ 2C \int \int \int Im[\bar{w} \mathcal N](t,y) \frac{(x - y)}{|x - y|} \cdot Im[\bar{v} (\nabla - i \xi(t)) v](t,x) dx dy dt.
\end{equation}

\noindent First by $(\ref{12.6})$, $(\ref{12.7})$, $(\ref{12.8})$, $\| u \|_{L_{t,x}^{4}(J_{l} \times \mathbf{R}^{2})} = 1$, and the definition of $w$,

\begin{equation}\label{12.39}
\aligned
(\ref{12.36}) \lesssim \| \mathcal N_{1} \|_{L_{t,x}^{4/3}(J_{l} \times \mathbf{R}^{2})} \| (\nabla - i \xi(t)) w \|_{L_{t,x}^{4}(J_{l} \times \mathbf{R}^{2})} \| v \|_{L_{t}^{\infty} L_{x}^{2}(J_{l} \times \mathbf{R}^{2})}^{2} \\ + (2^{-l_{2}} \int_{J_{l}} |\xi'(t)| \| P_{\xi(t), l_{2} - 3 \leq \cdot \leq l_{2} + 3} u(t) \|_{L_{x}^{2}(\mathbf{R}^{2})} \| (\nabla - i \xi(t)) P_{\xi(t), \leq l_{2}} u(t) \|_{L_{x}^{2}(\mathbf{R}^{2})} dt) \| v \|_{L_{t}^{\infty} L_{x}^{2}(J_{l} \times \mathbf{R}^{2})}^{2} \\ \lesssim 2^{l_{2}} \| v_{0} \|_{L^{2}}^{2} + \| v_{0} \|_{L^{2}}^{2} (2^{-l_{2}} \int_{J_{l}} |\xi'(t)| \| P_{\xi(t), l_{2} - 3 \leq \cdot \leq l_{2} + 3} u(t) \|_{L_{x}^{2}(\mathbf{R}^{2})} \| (\nabla - i \xi(t)) P_{\xi(t), \leq l_{2}} u(t) \|_{L_{x}^{2}(\mathbf{R}^{2})} dt).
\endaligned
\end{equation}

\noindent Next, integrating $(\ref{12.37})$ by parts,

\begin{equation}\label{12.40}
(\ref{12.37}) = (\ref{12.36}) + C \int_{J_{l}} \int \int |v(t,y)|^{2} \frac{1}{|x - y|} Re[\bar{w} \mathcal N](t,x) dx dy dt.
\end{equation}

\noindent By the Hardy - Littlewood - Sobolev inequality and Sobolev embedding theorem,

\begin{equation}\label{12.41}
\aligned
\int_{J_{l}} \int \int |v(t,y)|^{2} \frac{1}{|x - y|} |w(t,x)| |\mathcal N_{1}(t,x)| dx dy dt \\ \lesssim \| v \|_{L_{t}^{6} L_{x}^{3}(J_{l} \times \mathbf{R}^{2})}^{2} \| u \|_{L_{t}^{9/2} L_{x}^{18/5}(J_{l} \times \mathbf{R}^{2})}^{3} \| w \|_{L_{t,x}^{\infty}(J_{l} \times \mathbf{R}^{2})} \lesssim 2^{l_{2}} \| v_{0} \|_{L_{x}^{2}(\mathbf{R}^{2})}^{2},
\endaligned
\end{equation}

\noindent and by Hardy's inequality,

\begin{equation}\label{12.42}
\aligned
\| \frac{1}{|x - y|^{1/2}} P_{\xi(t), \leq l_{2}} u \|_{L_{x}^{2}(\mathbf{R}^{2})} \lesssim \| |\nabla - i \xi(t)|^{1/2} P_{\xi(t), \leq l_{2}} u \|_{L_{x}^{2}(\mathbf{R}^{2})}, \\
\| \frac{1}{|x - y|^{1/2}} P_{\xi(t), l_{2} - 3 \leq \cdot \leq l_{2} + 3} u(t) \|_{L^{2}} \lesssim 2^{l_{2}/2} \| P_{\xi(t), l_{2} - 3 \leq \cdot \leq l_{2} + 3} u(t) \|_{L^{2}},
\endaligned
\end{equation}

\noindent so by $(\ref{12.7})$ and $(\ref{12.8})$,

\begin{equation}\label{12.43}
\aligned
\int_{J_{l}} \int \int |v(t,y)|^{2} \frac{1}{|x - y|} |w(t,x)| |\mathcal N_{2}(t,x)| dx dy dt \\ \lesssim 2^{-l_{2}/2} \| v_{0} \|_{L_{x}^{2}(\mathbf{R}^{2})}^{2} (\int_{J_{l}} |\xi'(t)| \| P_{\xi(t), l_{2} - 3 \leq \cdot \leq l_{2} + 3} u(t) \|_{L_{x}^{2}(\mathbf{R}^{2})} \| |\nabla - i \xi(t)|^{1/2} P_{\xi(t), \leq l_{2}} u(t) \|_{L_{x}^{2}(\mathbf{R}^{2})} dt).
\endaligned
\end{equation}

\noindent Finally take $(\ref{12.38})$.

\begin{equation}\label{12.44}
\int_{J_{l}} \int \int Im[\bar{\mathcal N}_{1} w](t,y) \frac{(x - y)}{|x - y|} \cdot Im[\bar{v} (\nabla - i \xi(t)) v](t,x) dx dy dt \lesssim 2^{i} \| v_{0} \|_{L_{x}^{2}(\mathbf{R}^{2})}^{2} \| u \|_{L_{t,x}^{4}(J_{l} \times \mathbf{R}^{2})}^{4}.
\end{equation}

\noindent Also, by $(\ref{12.7})$ and $(\ref{12.8})$,

\begin{equation}\label{12.45}
\aligned
\int_{J_{l}} \int \int Im[\bar{\mathcal N_{2}} P_{\xi(t), l_{2} - 5 \leq \cdot \leq l_{2}} u](t,y) \frac{(x - y)}{|x - y|} \cdot Im[\bar{v} (\nabla - i \xi(t)) v](t,x) dx dy dt \\ \lesssim 2^{i - 2l_{2}} \| v_{0} \|_{L^{2}(\mathbf{R}^{2})}^{2} (\int_{J_{l}} |\xi'(t)| \| P_{\xi(t), l_{2} - 3 \leq \cdot \leq l_{2} + 3} u(t) \|_{L_{x}^{2}(\mathbf{R}^{2})} \| (\nabla - i \xi(t)) P_{\xi(t), \leq l_{2}} u \|_{L_{x}^{2}(\mathbf{R}^{2})} dt).
\endaligned
\end{equation}

\noindent $(\ref{12.8})$ implies $\bar{\mathcal N_{2}} P_{\xi(t), \leq l_{2} - 5} u$ is supported on $|\xi| \geq 2^{l_{2} - 5}$, so integrating by parts,

\begin{equation}\label{12.46}
\aligned
\int_{J_{l}} \int \int (\frac{\Delta_{y}}{\Delta_{y}} Im[\bar{\mathcal N}_{2} P_{\xi(t), \leq l_{2} - 5} u])(t,y) \frac{(x - y)}{|x - y|} \cdot Im[\bar{v} (\nabla - i \xi(t)) v](t,x) dx dy dt \\
= \int_{J_{l}} \int \int \frac{\partial_{k}}{\Delta_{y}} Im[\bar{\mathcal N}_{2} P_{\xi(t), \leq l_{2} - 5} u](t,y) [\frac{\delta_{jk}}{|x - y|} - \frac{(x - y)_{j} (x - y)_{k}}{|x - y|^{3}}] Im[\bar{v} (\partial_{j} - i \xi_{j}(t)) v](t,x) dx dy dt,
\endaligned
\end{equation}

\noindent which by $(\ref{12.42})$, Bernstein's inequality, and the Fourier support of $\bar{\mathcal N_{2}} P_{\xi(t), l_{2} - 5} u$,

\begin{equation}\label{12.47}
\lesssim 2^{i - 2l_{2}} \| v_{0} \|_{L^{2}(\mathbf{R}^{2})}^{2} (\int_{J_{l}} |\xi'(t)| \| P_{\xi(t), l_{2} - 3 \leq \cdot \leq l_{2} + 3} u(t) \|_{L_{x}^{2}(\mathbf{R}^{2})} \| (\nabla - i \xi(t)) P_{\xi(t), \leq l_{2}} u(t) \|_{L_{x}^{2}(\mathbf{R}^{2})} dt).
\end{equation}

\noindent Collecting $(\ref{12.23})$ - $(\ref{12.25})$, $(\ref{12.39})$, $(\ref{12.41})$, $(\ref{12.43})$, $(\ref{12.44})$, $(\ref{12.45})$, and $(\ref{12.47})$,

\begin{equation}\label{12.48}
\aligned
\| \nabla (v \bar{w}) \|_{L_{t,x}^{2}(G_{\alpha}^{i} \times \mathbf{R}^{2})}^{2} \lesssim 2^{i + l_{2}} \| v_{0} \|_{L^{2}(\mathbf{R}^{2})}^{2} \\ + 2^{i - 2l_{2}} \| v_{0} \|_{L^{2}}^{2} (\int_{J_{l}} |\xi'(t)| \| P_{\xi(t), l_{2} - 3 \leq \cdot \leq l_{2} + 3} u(t) \|_{L_{x}^{2}} \| (\nabla - i \xi(t)) P_{\xi(t), \leq l_{2}} u(t) \|_{L_{x}^{2}} dt)).
\endaligned
\end{equation}

\noindent Now by $(\ref{2.20})$, since $\hat{v}_{0}(\xi)$ is supported on $2^{i - 5} \leq |\xi - \xi(G_{\alpha}^{i})| \leq 2^{i + 5}$, $\hat{w}(t,\xi)$ is supported on $|\xi - \xi(t)| \leq 2^{i - 10}$, and $|\xi(t) - \xi(G_{\alpha}^{i})| << 2^{i}$, $v \bar{w}$ is supported on $|\xi| \sim 2^{i}$. Therefore, by Bernstein's inequality,

\begin{equation}\label{12.49}
\aligned
\| v \bar{w} \|_{L_{t,x}^{2}(G_{\alpha}^{i} \times \mathbf{R}^{2})}^{2} \lesssim 2^{l_{2} - i} \| v_{0} \|_{L^{2}(\mathbf{R}^{2})}^{2} \\ + 2^{-l_{2} - i} \| v_{0} \|_{L^{2}}^{2} (\int_{J_{l}} |\xi'(t)| \| P_{\xi(t), l_{2} - 3 \leq \cdot \leq l_{2} + 3} u(t) \|_{L_{x}^{2}} \| (\nabla - i \xi(t)) P_{\xi(t), \leq l_{2}} u(t) \|_{L_{x}^{2}} dt)).
\endaligned
\end{equation}

\noindent Since 

\begin{equation}\label{12.50}
\| \bar{v} w \|_{L_{t,x}^{2}}^{2} = \| |v|^{2} |w|^{2} \|_{L_{t,x}^{1}} = \| vw \|_{L_{t,x}^{2}}^{2},
\end{equation}

\noindent the proof of theorem $\ref{t3.2.0}$ is complete. $\Box$\vspace{5mm}

\noindent \textbf{Remark:} Notice that if $\xi(t) \equiv 0$, such as when $u$ has radial symmetry, theorem $\ref{t3.2.0}$ merely gives the bilinear estimate from $(\ref{1.3})$. However, in the case when $\xi(t)$ is free to move around, theorem $\ref{t3.2.0}$ represents an improvement over $(\ref{1.3})$.\vspace{5mm}

\noindent By H{\"o}lder's lemma,

\begin{equation}\label{12.51}
\aligned
\| \int_{J_{l}} e^{-it \Delta} P_{\xi(G_{\alpha}^{i}), i - 2 \leq \cdot \leq i + 2} ((P_{\xi(G_{\alpha}^{i}), i - 5 \leq \cdot \leq i + 5} u)(P_{\xi(t), l_{2}} u)(P_{\xi(t), \leq l_{2}} u)) dt \|_{L^{2}(\mathbf{R}^{2})} \\
\leq \sup_{\| v_{0} \|_{L^{2}} = 1} \| (e^{it \Delta} v_{0}) (P_{\xi(G_{\alpha}^{i}), i - 5 \leq \cdot \leq i + 5} u)(P_{\xi(t), l_{2}} u)(P_{\xi(t), \leq l_{2}} u) \|_{L_{t,x}^{1}(J_{l} \times \mathbf{R}^{2})},
\endaligned
\end{equation}

\noindent for $\hat{v}_{0}(\xi)$ supported on $2^{i - 2} \leq |\xi - \xi(G_{\alpha}^{i})| \leq 2^{i + 2}$. Applying theorem $\ref{t3.2.0}$ to atoms of $P_{\xi(G_{\alpha}^{i}), i - 5 \leq \cdot \leq i + 5} u$ and $e^{it \Delta} v_{0}$ (see lemma $\ref{l2.6.2}$) implies

\begin{equation}\label{12.52}
\aligned
(\ref{12.51}) \lesssim \| (e^{it \Delta} v_{0})(P_{\xi(t), l_{2}} u) \|_{L_{t,x}^{2}(J_{l} \times \mathbf{R}^{2})} \| (P_{\xi(G_{\alpha}^{i}), i - 5 \leq \cdot \leq i + 5} u)(P_{\xi(t), \leq l_{2}} u) \|_{L_{t,x}^{2}(J_{l} \times \mathbf{R}^{2})} \\
\lesssim 2^{l_{2} - i} \| P_{\xi(G_{\alpha}^{i}), i - 5 \leq \cdot \leq i + 5} u \|_{U_{\Delta}^{2}(G_{\alpha}^{i} \times \mathbf{R}^{2})} \\ + 2^{-l_{2} - i} \| P_{\xi(G_{\alpha}^{i}), i - 5 \leq \cdot \leq i + 5} u \|_{U_{\Delta}^{2}(J_{l} \times \mathbf{R}^{2})} (\int_{J_{l}} |\xi'(t)| \| P_{\xi(t), l_{2} - 3 \leq \cdot \leq l_{2} + 3} u(t) \|_{L^{2}} \| (\nabla - i \xi(t)) P_{\xi(t), \leq l_{2}} u(t) \|_{L^{2}} dt).
\endaligned
\end{equation}

\noindent Rearranging the order of summation,

\begin{equation}\label{12.53}
\sum_{0 \leq l_{2} \leq i - 10} \sum_{J_{l} \subset G_{\alpha}^{i}; N(J_{l}) \geq 2^{i - 5} \epsilon_{3}^{1/2}} 2^{l_{2} - i} \leq 2^{-i} \epsilon_{3}^{-1/2} \sum_{J_{l} \subset G_{\alpha}^{i}} N(J_{l}) \lesssim 2^{-i} \epsilon_{3}^{-1/2} \sum_{J_{l} \subset G_{\alpha}^{i}} N(J_{l}) \lesssim \epsilon_{3}^{1/2}.
\end{equation}

\noindent Also by $(\ref{2.20})$, and conservation of mass,

\begin{equation}\label{12.54}
\sum_{0 \leq l_{2} \leq i - 10} \sum_{J_{l} \subset G_{\alpha}^{i}} 2^{-l_{2} - i} (\int_{J_{l}} |\xi'(t)| \| P_{\xi(t), l_{2} - 3 \leq \cdot \leq l_{2} + 3} u(t) \|_{L^{2}} \| (\nabla - i \xi(t)) P_{\xi(t), \leq l_{2}} u(t) \|_{L^{2}} dt) \lesssim \epsilon_{3} \epsilon_{1}^{-1}.
\end{equation}

\noindent Combining $(\ref{12.53})$ and $(\ref{12.54})$,

\begin{equation}\label{12.55}
(\ref{3.29.1}) \lesssim \epsilon_{3}^{1/2} \| P_{\xi(G_{\alpha}^{i}), i - 5 \leq \cdot \leq i + 5} u \|_{U_{\Delta}^{2}(G_{\alpha}^{i} \times \mathbf{R}^{2})}.
\end{equation}

\noindent Next take $(\ref{3.30})$.


\begin{theorem}[Second bilinear Strichartz estimate]\label{t3.2.1}
Suppose $v_{0} \in L^{2}(\mathbf{R}^{2})$ has Fourier transform $\hat{v}_{0}(\xi)$ supported on $2^{i - 5} \leq |\xi - \xi(G_{\alpha}^{i})| \leq 2^{i + 5}$. Then for $0 \leq l_{2} \leq i - 10$, $G_{\beta}^{l_{2}} \subset G_{\alpha}^{i}$,

\begin{equation}\label{10.1}
\| (e^{it \Delta} v_{0})(P_{\xi(t), \leq l_{2}} u) \|_{L_{t,x}^{2}(G_{\beta}^{l_{2}} \times \mathbf{R}^{2})}^{2} \lesssim \| v_{0} \|_{L^{2}(\mathbf{R}^{2})}^{2} (1 + \| u \|_{\tilde{X}_{i}(G_{\alpha}^{i} \times \mathbf{R}^{2})}^{4}).
\end{equation}
\end{theorem}

\noindent \emph{Proof:} Applying $(\ref{12.2})$ - $(\ref{12.22})$, Galilean invariance of the potential, $(\ref{12.31})$, $(\ref{12.32})$, ($(\ref{12.36})$ - $(\ref{12.38})$), $(\ref{2.20})$ (which implies $|\xi(t) - \xi(G_{\alpha}^{i})| << 2^{i}$, Bernstein's inequality (see $(\ref{12.48})$ - $(\ref{12.50})$), and the fundamental theorem of calculus,

\begin{equation}\label{10.23}
\int_{G_{\beta}^{l_{2}}} \int |\overline{w(t,x)} v(t,x)|^{2} dx dt \lesssim 2^{l_{2} - 2i} \sup_{t \in G_{\beta}^{l_{2}}} |M(t)| 
\end{equation}

\begin{equation}\label{10.24}
+ 2^{l_{2} - 2i} \int_{G_{\beta}^{l_{2}}} \int \int \frac{1}{|x - y|} |v(t,y)|^{2} |w(t,x)|^{4} dx dy dt
\end{equation}

\begin{equation}\label{10.26}
+ 2^{l_{2} - 2i} \int_{G_{\beta}^{l_{2}}} \int \int |v(t,y)|^{2} \frac{(x - y)}{|x - y|} \cdot Re[\bar{\mathcal N} (\nabla - i \xi(t)) w](t,x) dx dy dt
\end{equation}

\begin{equation}\label{10.27}
+ 2^{l_{2} - 2i} \int_{G_{\beta}^{l_{2}}} \int \int |v(t,y)|^{2} \frac{x - y}{|x - y|} \cdot Re[\bar{w} (\nabla - i \xi(t)) \mathcal N](t,x) dx dy dt
\end{equation}

\begin{equation}\label{10.25}
+ 2^{l_{2} - 2i} \int_{G_{\beta}^{l_{2}}} \int \int Im[\bar{w} \mathcal N](t,y) \frac{(x - y)}{|x - y|} Im[\bar{v} (\nabla - i \xi(t)) v](t,x) dx dy dt.
\end{equation}

\noindent Once again

\begin{equation}\label{10.30}
\aligned
M(t) =  \int \int |w(t,y)|^{2} \frac{(x - y)}{|x - y|} \cdot Im[\bar{v} (\nabla - i \xi(t)) v](t,x) dx dy \\ + \int \int |v(t,y)|^{2} \frac{(x - y)}{|x - y|} \cdot Im[\bar{w} (\nabla - i \xi(t)) w](t,x) dx dy \lesssim 2^{i} \| w \|_{L_{t}^{\infty} L_{x}^{2}}^{2} \| v_{0} \|_{L_{x}^{2}}^{2} \lesssim 2^{i} \| v_{0} \|_{L^{2}}^{2},
\endaligned
\end{equation}

\noindent so $2^{l_{2} - 2i} \sup_{t \in G_{\alpha}^{i}} |M(t)|$ is bounded by the right hand side of $(\ref{10.1})$.\vspace{5mm}

\noindent Again by the Hardy - Littlewood - Sobolev inequality and Strichartz estimates,

\begin{equation}\label{10.34}
\aligned
2^{l_{2} - 2i} \int_{G_{\beta}^{l_{2}}} \int \int \frac{1}{|x - y|} |v(t,y)|^{2} |w(t,x)|^{4} dx dy dt \\ \lesssim \| v \|_{L_{t}^{6} L_{x}^{3}(G_{\alpha}^{i} \times \mathbf{R}^{2})}^{2} \| w \|_{L_{t}^{6} L_{x}^{24/5}(G_{\alpha}^{i} \times \mathbf{R}^{2})}^{4} \lesssim \| v_{0} \|_{L^{2}(\mathbf{R}^{2})}^{2} \| w \|_{L_{t}^{6} L_{x}^{24/5}(G_{\alpha}^{i} \times \mathbf{R}^{2})}^{4}.
\endaligned
\end{equation}

\noindent By $(\ref{2.37})$, $G_{\beta}^{l_{2}} \subset G_{\alpha}^{i}$, and $l_{2} \leq i - 10$,

\begin{equation}\label{10.35}
\| w \|_{L_{t}^{6} L_{x}^{24/5}(G_{\beta}^{l_{2}} \times \mathbf{R}^{2})} \lesssim \sum_{0 \leq l_{3} \leq l_{2}} 2^{l_{3}/4} \| P_{\xi(t), l_{3}} u \|_{L_{t}^{6} L_{x}^{3}(G_{\beta}^{l_{2}} \times \mathbf{R}^{2})} \lesssim 2^{l_{2}/4} \| u \|_{\tilde{X}_{i}(G_{\beta}^{l_{2}} \times \mathbf{R}^{2})},
\end{equation}

\noindent so

\begin{equation}\label{10.36}
(\ref{10.34}) \lesssim 2^{2 l_{2} - 2i} \| v_{0} \|_{L^{2}}^{2} \| u \|_{\tilde{X}_{i}(G_{\alpha}^{i} \times \mathbf{R}^{2})}^{4}.
\end{equation}

\noindent Therefore, $(\ref{10.24})$ is bounded by the right hand side of $(\ref{10.1})$.\vspace{5mm}

\noindent Now for $(\ref{10.26})$ - $(\ref{10.25})$, to simplify notation let every $L_{t}^{p} L_{x}^{q}$ norm will be taken over $G_{\beta}^{l_{2}} \times \mathbf{R}^{2}$ unless otherwise stated. By $(\ref{12.39})$,

\begin{equation}\label{10.39.1}
(\ref{10.26}) \lesssim 2^{-2i} \| v \|_{L_{t}^{\infty} L_{x}^{2}}^{2} (\int_{G_{\beta}^{l_{2}}} |\xi'(t)| \| P_{\xi(t), l_{2} - 3 \leq \cdot \leq l_{2} + 3} u(t) \|_{L^{2}} \| (\nabla - i \xi(t)) P_{\xi(t), \leq l_{2}} u(t) \|_{L^{2}} dt)
\end{equation}

\begin{equation}\label{10.39}
+ 2^{l_{2} - 2i} \| v \|_{L_{t}^{\infty} L_{x}^{2}}^{2} \| \mathcal N_{1} \|_{L_{t}^{3/2} L_{x}^{6/5}} \| (\nabla - i \xi(t)) w \|_{L_{t}^{3} L_{x}^{6}}.
\end{equation}

\noindent By $(\ref{2.20})$ and conservation of mass,

\begin{equation}\label{5.31.1}
(\ref{10.39.1}) \lesssim 2^{2 l_{2} - 2i} \| v_{0} \|_{L_{x}^{2}}^{2}.
\end{equation}

\noindent By $(\ref{2.37})$,

\begin{equation}\label{10.40}
(\ref{10.39}) \lesssim 2^{l_{2} - 2i} \| v_{0} \|_{L^{2}}^{2} \| u \|_{\tilde{X}_{i}(G_{\beta}^{l_{2}} \times \mathbf{R}^{2})} \| P_{\xi(t), \leq l_{2}} F(u) - F(P_{\xi(t), \leq l_{2}} u) \|_{L_{t}^{3/2} L_{x}^{6/5}}.
\end{equation}

\noindent Split $u = u_{h} + u_{l}$, where $u_{l} = P_{\xi(t), \leq l_{2} - 5} u$. Following the same analysis as in $(\ref{3.9.2})$ for each $t \in G_{\beta}^{l_{2}}$ separately,

\begin{equation}\label{10.41}
P_{\xi(t), \leq l_{2}} F(u_{l}) - F(P_{\xi(t), \leq l_{2}} u_{l}) = 0.
\end{equation}

\noindent Next,

\begin{equation}\label{10.42}
P_{\xi(t), \leq l_{2}} (|u_{l}|^{2} u_{h}) - |P_{\xi(t), \leq l_{2}} u_{l}|^{2} (P_{\xi(t), \leq i - 10} u_{h}) = P_{\xi(t), \leq l_{2}} (|u_{l}|^{2} u_{h}) - |u_{l}|^{2} (P_{\xi(t), \leq l_{2}} u_{h}).
\end{equation}

\noindent By the fundamental theorem of calculus,

\begin{equation}\label{5.29}
|\phi(\frac{\xi_{2} + \xi_{1} - \xi(t)}{2^{l_{2}}}) - \phi(\frac{\xi_{1} - \xi(t)}{2^{l_{2}}})| \lesssim 2^{-l_{2}} |\xi_{2}|.
\end{equation}

\noindent Therefore,

\begin{equation}\label{5.29.1}
\| P_{\xi(t), \leq l_{2}} (|u_{l}|^{2} u_{h}) - |u_{l}|^{2} (P_{\xi(t), \leq l_{2}} u_{h}) \|_{L_{t}^{3/2} L_{x}^{6/5}} \lesssim 2^{-l_{2}} \| u_{h} \|_{L_{t}^{3} L_{x}^{6}} \| \nabla |u_{l}|^{2} \|_{L_{t}^{3} L_{x}^{3/2}}.
\end{equation}

\noindent By $(\ref{2.37})$ and the product rule,

\begin{equation}\label{5.29.2}
\| \nabla |u_{l}|^{2} \|_{L_{t}^{3} L_{x}^{3/2}} = \| \nabla ((e^{-ix \cdot \xi(t)} u_{l})(e^{ix \cdot \xi(t)} \bar{u}_{l})) \|_{L_{t}^{3} L_{x}^{3/2}}
\end{equation}

\begin{equation}\label{5.29.3}
\lesssim \| \nabla (e^{-ix \cdot \xi(t)} u_{l}) \|_{L_{t}^{3} L_{x}^{6}} \| u_{l} \|_{L_{t}^{\infty} L_{x}^{2}} \lesssim 2^{l_{2}} \| u \|_{\tilde{X}_{l_{2}}(G_{\beta}^{l_{2}} \times \mathbf{R}^{2})} \lesssim 2^{l_{2}} \| u \|_{\tilde{X}_{i}(G_{\alpha}^{i} \times \mathbf{R}^{2})}.
\end{equation}

\noindent Therefore,

\begin{equation}\label{5.30}
\| P_{\xi(t), \leq l_{2}}(|u_{l}|^{2} u_{h}) - |u_{l}|^{2} (P_{\xi(t), \leq l_{2}} u_{h}) \|_{L_{t}^{3/2} L_{x}^{6/5}} \lesssim \| u \|_{\tilde{X}_{i}(G_{\alpha}^{i} \times \mathbf{R}^{2})}^{2}.
\end{equation}

\noindent Similarly,

\begin{equation}\label{5.30.1}
\| P_{\xi(t), \leq l_{2}} (u_{l}^{2} \bar{u}_{h}) - (u_{l})^{2} (P_{\xi(t), \leq l_{2}} \bar{u}_{h}) \|_{L_{t}^{3/2} L_{x}^{6/5}} \lesssim \| u \|_{\tilde{X}_{i}(G_{\alpha}^{i} \times \mathbf{R}^{2})}^{2}.
\end{equation}

\noindent Finally, by $(\ref{2.37})$,

\begin{equation}\label{5.30.2}
\| u u_{h}^{2} \|_{L_{t}^{3/2} L_{x}^{6/5}} \lesssim \| u_{h} \|_{L_{t}^{3} L_{x}^{6}}^{2} \| u \|_{L_{t}^{\infty} L_{x}^{2}} \lesssim \| u \|_{\tilde{X}_{i}(G_{\alpha}^{i} \times \mathbf{R}^{2})}^{2}.
\end{equation}

\noindent Therefore,

\begin{equation}\label{5.31}
(\ref{10.39}) \lesssim 2^{l_{2}} \| v(0) \|_{L^{2}(\mathbf{R}^{2})}^{2} \| u \|_{\tilde{X}_{i}(G_{\alpha}^{i} \times \mathbf{R}^{2})}^{3},
\end{equation}

\noindent and

\begin{equation}\label{5.31.2}
(\ref{10.26}) \lesssim 2^{2l_{2} - 2i} \| v_{0} \|_{L^{2}}^{2} (1 + \| u \|_{\tilde{X}_{i}(G_{\alpha}^{i} \times \mathbf{R}^{2})}^{3}).
\end{equation}

\noindent The right hand side of this term is clearly bounded by the right hand side of $(\ref{10.1})$.\vspace{5mm}

\noindent Again integrating by parts in space,

\begin{equation}\label{5.32}
2^{l_{2} - 2i} \int_{G_{\beta}^{l_{2}}} \int \int |v(t,y)|^{2} \frac{(x - y)}{|x - y|} \cdot Im[\bar{w} (\nabla - i \xi(t)) \mathcal N](t,x) dx dy dt 
\end{equation}

\begin{equation}\label{5.33}
= (\ref{10.26}) - 2^{l_{2} - 2i} \int_{G_{\beta}^{l_{2}}} \int \int \frac{1}{|x - y|} |v(t,y)|^{2}  Re[\bar{w} \mathcal N](t,x) dx dy dt.
\end{equation}

\noindent Then by $(\ref{12.43})$, $(\ref{2.20})$, $(\ref{5.30})$ - $(\ref{5.30.2})$, $(\ref{2.37})$, the Hardy - Littlewood - Sobolev lemma, and the Sobolev embedding theorem

\begin{equation}\label{5.33.1}
\aligned
(\ref{5.33}) \lesssim 2^{l_{2} - 2i} \| v \|_{L_{t}^{6} L_{x}^{3}}^{2} \| \mathcal N_{1} \|_{L_{t}^{3/2} L_{x}^{6/5}} \| w \|_{L_{t,x}^{\infty}} + 2^{2l_{2} - 2i} \| v_{0} \|_{L^{2}}^{2}(1 + \| u \|_{\tilde{X}_{i}(G_{\alpha}^{i} \times \mathbf{R}^{2})}^{3}) \\ \lesssim 2^{2l_{2} - 2i} \| v_{0} \|_{L^{2}}^{2} (1 + \| u \|_{\tilde{X}_{i}(G_{\alpha}^{i} \times \mathbf{R}^{2})}^{2}).
\endaligned
\end{equation}

\noindent Therefore, by $(\ref{5.31.2})$ and $(\ref{5.33.1})$,

\begin{equation}\label{5.33.2}
(\ref{10.27}) \lesssim 2^{2l_{2} - 2i} \| v_{0} \|_{L^{2}}^{2} (1 + \| u \|_{\tilde{X}_{i}(G_{\alpha}^{i} \times \mathbf{R}^{2})}^{3}).
\end{equation}


\noindent The term $(\ref{10.25})$ is much more difficult to estimate, and so the computations will be postponed to the appendix.

\begin{lemma}\label{appendix}
\begin{equation}
(\ref{10.25}) \lesssim  \| v_{0} \|_{L^{2}}^{2} (1 + \| u \|_{\tilde{X}_{i}(G_{\alpha}^{i} \times \mathbf{R}^{2})}^{4}).
\end{equation}
\end{lemma}

\noindent \emph{Proof:} See the appendix.\vspace{5mm}

\noindent Assuming lemma $\ref{appendix}$ is true, the proof of theorem $\ref{t3.2.1}$ is complete. $\Box$\vspace{5mm}

\noindent So by lemma $\ref{l2.6.2}$, for any $G_{\beta}^{i - 10} \subset G_{\alpha}^{i}$,

\begin{equation}
\| (P_{\xi(G_{\alpha}^{i}), i - 5 \leq \cdot \leq i + 5} u)(P_{\xi(t), \leq i - 10} u) \|_{L_{t,x}^{2}(G_{\beta}^{i - 10} \times \mathbf{R}^{2})} \lesssim \| P_{\xi(G_{\alpha}^{i}), i - 5 \leq \cdot \leq i + 5} u \|_{U_{\Delta}^{2}(G_{\alpha}^{i} \times \mathbf{R}^{2})} (1 + \| u \|_{\tilde{X}_{i}(G_{\alpha}^{i} \times \mathbf{R}^{2})}^{2}),
\end{equation}

\noindent and so by $(\ref{1.12.1})$ - $(\ref{1.14})$, for any $l_{2} \leq i - 10$,

\begin{equation} 
\| (P_{\xi(G_{\alpha}^{i}), i - 5 \leq \cdot \leq i + 5} u)(P_{\xi(t), \leq l_{2}} u) \|_{L_{t,x}^{2}(G_{\beta}^{i - 10} \times \mathbf{R}^{2})} \lesssim \| P_{\xi(G_{\alpha}^{i}), i - 5 \leq \cdot \leq i + 5} u \|_{U_{\Delta}^{2}(G_{\alpha}^{i} \times \mathbf{R}^{2})} (1 + \| u \|_{\tilde{X}_{i}(G_{\alpha}^{i} \times \mathbf{R}^{2})}^{2}).
\end{equation}

\noindent Since $G_{\alpha}^{i}$ is the union of $2^{10}$ subintervals $G_{\beta}^{i - 10}$, this means that

\begin{equation} 
\| (P_{\xi(G_{\alpha}^{i}), i - 5 \leq \cdot \leq i + 5} u)(P_{\xi(t), \leq l_{2}} u) \|_{L_{t,x}^{2}(G_{\alpha}^{i} \times \mathbf{R}^{2})} \lesssim \| P_{\xi(G_{\alpha}^{i}), i - 5 \leq \cdot \leq i + 5} u \|_{U_{\Delta}^{2}(G_{\alpha}^{i} \times \mathbf{R}^{2})} (1 + \| u \|_{\tilde{X}_{i}(G_{\alpha}^{i} \times \mathbf{R}^{2})}^{2}).
\end{equation}

\noindent Now suppose that for each $G_{\beta}^{l_{2}} \subset G_{\alpha}^{i}$, $\| v_{\beta}^{l_{2}} \|_{V_{\Delta}^{2}(G_{\beta}^{l_{2}} \times \mathbf{R}^{2})} = 1$ and $\hat{v}_{\beta}^{l_{2}}(t, \xi)$ supported on $2^{i - 2} \leq |\xi - \xi(G_{\alpha}^{i})| \leq 2^{i + 2}$. Then,

\begin{equation}\label{3.30.2}
(\sum_{G_{\beta}^{l_{2}} \subset G_{\alpha}^{i}; N(G_{\beta}^{l_{2}}) \leq 2^{l_{2} - 5} \epsilon_{3}^{1/2}} \| v_{\beta}^{l_{2}} ((P_{\xi(G_{\alpha}^{i}), i - 5 \leq \cdot \leq i + 5} u)(P_{\xi(t), l_{2}} u)(P_{\xi(t), \leq l_{2}} u))  \|_{L_{t,x}^{1}(G_{\beta}^{l_{2}} \times \mathbf{R}^{2})}^{2})^{1/2}
\end{equation}

\begin{equation}\label{3.30.3}
\lesssim (\sup_{G_{\beta}^{l_{2}} \subset G_{\alpha}^{i}; N(G_{\beta}^{l_{2}}) \leq 2^{l_{2} - 5} \epsilon_{3}^{1/2}} \| v_{\beta}^{l_{2}} (P_{\xi(t), l_{2}} u) \|_{L_{t,x}^{2}(G_{\beta}^{l_{2}} \times \mathbf{R}^{2})}) \| P_{\xi(G_{\alpha}^{i}), i - 5 \leq \cdot \leq i + 5} u \|_{U_{\Delta}^{2}(G_{\alpha}^{i} \times \mathbf{R}^{2})} (1 + \| u \|_{\tilde{X}_{i}(G_{\alpha}^{i} \times \mathbf{R}^{2})}^{2}).
\end{equation}

\noindent Then by $V_{\Delta}^{2} \subset U_{\Delta}^{2 + \epsilon}$, $N(G_{\beta}^{l_{2}}) \leq \epsilon_{3}^{1/2} 2^{l_{2} - 5}$, and corollary $\ref{c1.3}$,

\begin{equation}\label{3.30.4}
\| v_{\beta}^{l_{2}} (P_{\xi(t), l_{2}} u) \|_{L_{t,x}^{2}(G_{\beta}^{l_{2}} \times \mathbf{R}^{2})} \lesssim \| v_{\beta}^{l_{2}} (P_{\xi(t), l_{2}} u) \|_{L_{t}^{3} L_{x}^{3/2}}^{1/2} \| v_{\beta}^{l_{2}} \|_{L_{t}^{3} L_{x}^{6}}^{1/2} \| P_{\xi(t), l_{2}} u \|_{L_{t}^{3} L_{x}^{6}}^{1/2} \lesssim 2^{(l_{2} - i)/6} \| u \|_{\tilde{Y}_{i}(G_{\alpha}^{i} \times \mathbf{R}^{2})}.
\end{equation}

\noindent Therefore, 

\begin{equation}
(\ref{3.30}) \lesssim \| P_{\xi(G_{\alpha}^{i}), i - 5 \leq \cdot \leq i + 5} u \|_{U_{\Delta}^{2}(G_{\alpha}^{i} \times \mathbf{R}^{2})} \| u \|_{\tilde{Y}_{i}(G_{\alpha}^{i} \times \mathbf{R}^{2})} (1 + \| u \|_{\tilde{X}_{i}(G_{\alpha}^{i} \times \mathbf{R}^{2})}^{2}).
\end{equation}

\noindent Finally, we estimate $(\ref{3.29})$. Again by H{\"o}lder's lemma

\begin{equation}\label{13.1}
\aligned
\| \int_{G_{\beta}^{l_{2}}} e^{-it \Delta} P_{\xi(G_{\alpha}^{i}), i - 2 \leq \cdot \leq i + 2} ((P_{\xi(G_{\alpha}^{i}), i - 5 \leq \cdot \leq i + 5} u)(P_{\xi(t), l_{2}} u)(P_{\xi(t), \leq l_{2}} u)) dt \|_{L^{2}} \\ = \sup_{\| v_{0} \|_{L^{2}} = 1} \| (e^{it \Delta} P_{\xi(G_{\alpha}^{i}), i - 2 \leq \cdot \leq i + 2} v_{0})(P_{\xi(G_{\alpha}^{i}), i - 5 \leq \cdot \leq i + 5} u)(P_{\xi(t), l_{2}} u)(P_{\xi(t), \leq l_{2}} u) \|_{L_{t,x}^{1}(G_{\beta}^{l_{2}} \times \mathbf{R}^{2})},
\endaligned
\end{equation}

\noindent so by $(\ref{1.3})$,

\begin{equation}\label{13.2}
(\ref{13.1}) \lesssim 2^{(l_{2} - i)/2} \| P_{\xi(G_{\beta}^{l_{2}}), l_{2}} u \|_{U_{\Delta}^{2}(G_{\beta}^{l_{2}} \times \mathbf{R}^{2})} \| (P_{\xi(G_{\alpha}^{i}), i - 5 \leq \cdot \leq i + 5} u)(P_{\xi(t), \leq l_{2}} u) \|_{L_{t,x}^{2}(G_{\beta}^{l_{2}} \times \mathbf{R}^{2})}.
\end{equation}

\noindent Therefore, by the Cauchy - Schwartz inequality,

\begin{equation}\label{13.4}
(\ref{3.29}) \lesssim \| u \|_{\tilde{Y}_{i}(G_{\alpha}^{i} \times \mathbf{R}^{2})} (\sum_{0 \leq l_{2} \leq i - 10} \| (P_{\xi(G_{\alpha}^{i}), i - 5 \leq \cdot \leq i + 5} u)(P_{\xi(t), \leq l_{2}} u) \|_{L_{t,x}^{2}(G_{\alpha}^{i} \times \mathbf{R}^{2})}^{2})^{1/2}.
\end{equation}

\noindent Theorem $\ref{tbigtheorem}$ will then finally follow from the final bilinear estimate,

\begin{theorem}[Third bilinear Strichartz estimate]\label{t3.2.2}
\begin{equation}\label{13.5}
\sum_{0 \leq l_{2} \leq i - 10} \| (e^{it \Delta} v(0))(P_{\xi(t), \leq l_{2}} u) \|_{L_{t,x}^{2}}^{2} \lesssim \| v(0) \|_{L^{2}}^{2} (1 + \| u \|_{\tilde{X}_{i}(G_{\alpha}^{i} \times \mathbf{R}^{2})}^{6}).
\end{equation}
\end{theorem}

\noindent Indeed, if $(\ref{13.5})$ holds, lemma $\ref{l2.6.2}$ implies

\begin{equation}\label{13.6}
(\ref{13.4}) \lesssim \| u \|_{\tilde{Y}_{i}(G_{\alpha}^{i} \times \mathbf{R}^{2})} \| P_{\xi(G_{\alpha}^{i}), i - 5 \leq \cdot \leq i + 5} u \|_{U_{\Delta}^{2}(G_{\alpha}^{i} \times \mathbf{R}^{2})} (1 + \| u \|_{\tilde{X}_{i}(G_{\alpha}^{i} \times \mathbf{R}^{2})}^{3}).
\end{equation}

\noindent \emph{Proof of theorem $\ref{t3.2.2}$:} This time let $v = e^{it \Delta} v_{0}$ and $w_{l_{2}} = P_{\xi(t), \leq l_{2}} u$. Then, as in $(\ref{10.23})$ - $(\ref{10.25})$,

\begin{equation}\label{13.7}
\sum_{0 \leq l_{2} \leq i - 10} \int_{G_{\alpha}^{i}} \int |\overline{w_{l_{2}}(t,x)} v(t,x)|^{2} dx dt \lesssim \sum_{0 \leq l_{2} \leq i - 10} 2^{l_{2} - 2i} \sup_{t \in G_{\alpha}^{i}} |M_{l_{2}}(t)| 
\end{equation}

\begin{equation}\label{13.8}
+ \sum_{0 \leq l_{2} \leq i - 10} 2^{l_{2} - 2i} \int_{G_{\beta}^{l_{2}}} \int \int \frac{1}{|x - y|} |v(t,y)|^{2} |w_{l_{2}}(t,x)|^{4} dx dy dt
\end{equation}

\begin{equation}\label{13.9}
+ \sum_{0 \leq l_{2} \leq i - 10} 2^{l_{2} - 2i} \int_{G_{\beta}^{l_{2}}} \int \int |v(t,y)|^{2} \frac{(x - y)}{|x - y|} \cdot Re[\bar{\mathcal N} (\nabla - i \xi(t)) w_{l_{2}}](t,x) dx dy dt
\end{equation}

\begin{equation}\label{13.10}
+ \sum_{0 \leq l_{2} \leq i - 10} 2^{l_{2} - 2i} \int_{G_{\beta}^{l_{2}}} \int \int |v(t,y)|^{2} \frac{x - y}{|x - y|} \cdot Re[\bar{w}_{l_{2}} (\nabla - i \xi(t)) \mathcal N](t,x) dx dy dt
\end{equation}

\begin{equation}\label{13.11}
+ \sum_{0 \leq l_{2} \leq i - 10} 2^{l_{2} - 2i}  \int_{G_{\beta}^{l_{2}}} \int \int Im[\bar{w}_{l_{2}} \mathcal N](t,y) \frac{(x - y)}{|x - y|} Im[\bar{v} (\nabla - i \xi(t)) v](t,x) dx dy dt,
\end{equation}

\noindent where 

\begin{equation}\label{13.12}
\aligned
M_{l_{2}}(t) =  \int \int |w_{l_{2}}(t,y)|^{2} \frac{(x - y)}{|x - y|} \cdot Im[\bar{v} (\nabla - i \xi(t)) v](t,x) dx dy \\ + \int \int |v(t,y)|^{2} \frac{(x - y)}{|x - y|} \cdot Im[\bar{w}_{l_{2}} (\nabla - i \xi(t)) w_{l_{2}}](t,x) dx dy \lesssim 2^{i} \| w_{l_{2}} \|_{L_{t}^{\infty} L_{x}^{2}}^{2} \| v_{0} \|_{L_{x}^{2}}^{2} \lesssim 2^{i} \| v_{0} \|_{L^{2}}^{2}.
\endaligned
\end{equation}

\noindent Therefore,

\begin{equation}\label{13.13}
\sum_{0 \leq l_{2} \leq i - 10} 2^{l_{2} - 2i} \sup_{t \in G_{\alpha}^{i}} |M_{l_{2}}(t)| \lesssim \sum_{0 \leq l_{2} \leq i - 10} 2^{l_{2} - i} \| v_{0} \|_{L^{2}}^{2} \lesssim \| v_{0} \|_{L^{2}}^{2}.
\end{equation}

\noindent Next, by $(\ref{10.36})$, $(\ref{5.31.2})$, and $(\ref{5.33.2})$, since there are $2^{i - l_{2}}$ intervals $G_{\beta}^{l_{2}} \subset G_{\alpha}^{i}$,

\begin{equation}\label{13.14}
(\ref{13.8}) + (\ref{13.9}) + (\ref{13.10}) \lesssim \sum_{0 \leq l_{2} \leq i - 10} 2^{i - l_{2}} 2^{2l_{2} - 2i} (1 + \| u \|_{\tilde{X}_{i}(G_{\alpha}^{i} \times \mathbf{R}^{2})}^{3}) \lesssim (1 + \| u \|_{\tilde{X}_{i}(G_{\alpha}^{i} \times \mathbf{R}^{2})}^{3}).
\end{equation}

\noindent Also, combining $(\ref{12.45})$, $(\ref{12.47})$, and $(\ref{2.20})$,

\begin{equation}\label{13.15}
\aligned
\sum_{0 \leq l_{2} \leq i - 10} 2^{l_{2} - 2i} \int_{G_{\alpha}^{i}} \int \int Im[\bar{w}_{l_{2}} \mathcal N_{2}] \frac{(x - y)}{|x - y|} \cdot Im[\bar{v} (\nabla - i \xi(t)) v](t,x) dx dy dt \\
\lesssim \sum_{0 \leq l_{2} \leq i - 10} \sum_{0 \leq l_{3} \leq l_{2}} 2^{l_{2} - 2i} 2^{i - 2l_{2}} 2^{l_{3}} \| v_{0} \|_{L^{2}}^{2} \int_{G_{\alpha}^{i}} |\xi'(t)| \| P_{\xi(t), l_{2} \leq \cdot \leq l_{2} + 3} u(t) \|_{L_{x}^{2}(\mathbf{R}^{2})} \| P_{\xi(t), l_{3}} u(t) \|_{L_{x}^{2}(\mathbf{R}^{2})} \\
\lesssim 2^{-i} \| v_{0} \|_{L^{2}}^{2} \int_{G_{\alpha}^{i}} |\xi'(t)| \| u(t) \|_{L^{2}}^{2} dt \lesssim \epsilon_{3} \epsilon_{1}^{-1} \| v_{0} \|_{L^{2}}^{2} \lesssim \| v_{0} \|_{L^{2}}^{2}.
\endaligned
\end{equation}

\noindent Therefore it only remains to show

\begin{equation}\label{13.16}
\sum_{0 \leq l_{2} \leq i - 10} 2^{l_{2} - 2i} \int_{G_{\alpha}^{i}} \int \int Im[\bar{w}_{l_{2}} \mathcal N_{1}](t,y) \frac{(x - y)}{|x - y|} \cdot Im[\bar{v} (\nabla - i \xi(t)) v](t,x) dx dy dt \lesssim 1.
\end{equation}

\noindent Again split $u = u_{h} + u_{l}$ where $u_{l} = P_{\xi(t), \leq l_{2} - 5} u$, and decompose

\begin{equation}\label{13.17}
Im[\bar{w}_{l_{2}} \mathcal N_{1}](t,y) = F_{0}(t,y) + F_{1}(t, y) + F_{2}(t, y) + F_{3}(t, y) + F_{4}(t, y),
\end{equation}

\noindent where $F_{0}$ has four $u_{l}$ terms and zero $u_{h}$ terms, $F_{1}$ has three $u_{l}$ terms and one $u_{h}$ term, and so on.

\begin{equation}\label{13.18}
\aligned
\sum_{0 \leq l_{2} \leq i - 10} 2^{l_{2} - 2i} \int_{G_{\alpha}^{i}} \int \int &[F_{2}(t,y) + F_{3}(t,y) + F_{4}(t, y)] \frac{(x - y)}{|x - y|} Im[\bar{v} (\nabla - i \xi(t)) v](t,x) dx dy dt \\ &\lesssim \| v_{0} \|_{L^{2}}^{2} \sum_{0 \leq l_{2} \leq i - 10} 2^{l_{2} - i}  \| F_{2} + F_{3} + F_{4} \|_{L_{t,x}^{1}(G_{\alpha}^{i} \times \mathbf{R}^{2})}.
\endaligned
\end{equation}

\noindent To simplify notation, $L_{t}^{p} L_{x}^{q}$ will refer to $L_{t}^{p} L_{x}^{q}(G_{\alpha}^{i} \times \mathbf{R}^{2})$ unless stated otherwise.


\begin{equation}\label{13.19}
\aligned
\sum_{0 \leq l_{2} \leq i - 10} 2^{l_{2} - i} \| F_{4} \|_{L_{t,x}^{1}} \lesssim \sum_{0 \leq l_{2} \leq i - 10} 2^{l_{2} - i} \| P_{\xi(t), l_{2} - 5 \leq \cdot \leq i} u \|_{L_{t,x}^{4}}^{4} +  \sum_{0 \leq l_{2} \leq i - 10} 2^{l_{2} - i} \| P_{\xi(t), \geq i} u \|_{L_{t,x}^{4}}^{4},
\endaligned
\end{equation}

\noindent By $(\ref{2.38})$ and conservation of mass,

\begin{equation}\label{13.19.1}
\sum_{0 \leq l_{2} \leq i - 10} 2^{l_{2} - i} \| P_{\xi(t), \geq i} u \|_{L_{t,x}^{4}(G_{\alpha}^{i} \times \mathbf{R}^{2})}^{4} \lesssim \| u \|_{\tilde{X}_{i}(G_{\alpha}^{i} \times \mathbf{R}^{2})}^{4}.
\end{equation}

\noindent Also, by conservation of mass and Young's inequality,

\begin{equation}
\aligned
\sum_{0 \leq l_{2} \leq i - 10} 2^{l_{2} - i} \| P_{\xi(t), l_{2} - 5 \leq \cdot \leq i} u \|_{L_{t,x}^{4}(G_{\alpha}^{i} \times \mathbf{R}^{2})}^{4} \\ \lesssim \sum_{0 \leq j_{1} \leq j_{2} \leq j_{3} \leq j_{4} \leq i} 2^{j_{1} - i} \| P_{\xi(t), j_{1}} u \|_{L_{t}^{\infty} L_{x}^{2}} \| P_{\xi(t), j_{2}} u \|_{L_{t}^{3} L_{x}^{6}} \| P_{\xi(t), j_{3}} u \|_{L_{t}^{3} L_{x}^{6}}  \| P_{\xi(t), j_{4}} u \|_{L_{t}^{3} L_{x}^{6}} \\
\lesssim \sum_{0 \leq j_{2} \leq j_{3} \leq j_{4} \leq i} 2^{j_{2} - i} \| P_{\xi(t), j_{2}} u \|_{L_{t}^{3} L_{x}^{6}} \| P_{\xi(t), j_{3}} u \|_{L_{t}^{3} L_{x}^{6}} \| P_{\xi(t), j_{4}} u \|_{L_{t}^{3} L_{x}^{6}} \\
\lesssim \sum_{0 \leq j_{2} \leq i} 2^{j_{2} - i} \| P_{\xi(t), j_{2}} u \|_{L_{t}^{3} L_{x}^{6}(G_{\alpha}^{i} \times \mathbf{R}^{2})}^{3}.
\endaligned
\end{equation}

\noindent By $(\ref{2.26})$ - $(\ref{2.28})$, the definition of $X(G_{\alpha}^{i} \times \mathbf{R}^{2})$, and $\tilde{X}_{i}(G_{\alpha}^{i} \times \mathbf{R}^{2})$,

\begin{equation}\label{13.20}
\aligned
\lesssim \sum_{0 \leq j_{2} \leq i} 2^{j_{2} - i} \sum_{G_{\beta}^{j_{2}} \subset G_{\alpha}^{i}} \| P_{\xi(G_{\beta}^{j_{2}}), j_{2} - 2 \leq \cdot \leq j_{2} + 2} u \|_{U_{\Delta}^{2}(G_{\beta}^{j_{2}} \times \mathbf{R}^{2})}^{3} \lesssim \| u \|_{\tilde{X}_{i}(G_{\alpha}^{i} \times \mathbf{R}^{2})}^{3}.
\endaligned
\end{equation}







\noindent Next take $F_{2}$.

\begin{equation}\label{13.23}
\aligned
\| F_{2} \|_{L_{t,x}^{1}(G_{\alpha}^{i} \times \mathbf{R}^{2})} \lesssim \| (\overline{P_{\xi(t), \geq l_{2} - 5} u})(P_{\xi(t), \leq l_{2} - 15} u) \|_{L_{t,x}^{2}(G_{\alpha}^{i} \times \mathbf{R}^{2})}^{2} \\ + \| (\overline{P_{\xi(t), \geq l_{2} - 5} u})(P_{\xi(t), l_{2} - 15 \leq \cdot \leq l_{2} - 5} u) \|_{L_{t,x}^{2}(G_{\alpha}^{i} \times \mathbf{R}^{2})}^{2}.
\endaligned
\end{equation}

\noindent Since $l_{2} \leq i - 10$, by Plancherel's theorem,

\begin{equation}\label{13.24}
\| (\overline{P_{\xi(t), \geq i} u})(P_{\xi(t), \leq l_{2} - 15} u) \|_{L_{t,x}^{2}(G_{\alpha}^{i} \times \mathbf{R}^{2})}^{2} \lesssim \sum_{l_{1} \geq i} \| (\overline{P_{\xi(t), l_{1}} u})(P_{\xi(t), \leq l_{2} - 15} u) \|_{L_{t,x}^{2}(G_{\alpha}^{i} \times \mathbf{R}^{2})}^{2}.
\end{equation}

\noindent Then by theorem $\ref{t3.2.1}$, $(\ref{1.12.1})$ - $(\ref{1.14})$, and $(\ref{2.26})$ - $(\ref{2.28})$,

\begin{equation}\label{13.25}
\aligned
\sum_{0 \leq l_{2} \leq i - 10} 2^{l_{2} - i} (\ref{13.24}) &\lesssim \sum_{0 \leq l_{2} \leq i - 10} 2^{l_{2} - i} \sum_{l_{1} \geq i} \| (\overline{P_{\xi(G_{\alpha}^{i}), l_{1} - 5 \leq \cdot \leq l_{1} + 5} u})(P_{\xi(t), \leq l_{2} - 5} u) \|_{L_{t,x}^{2}(G_{\alpha}^{i} \times \mathbf{R}^{2})} \\ &\lesssim \| u \|_{\tilde{X}_{i}(G_{\alpha}^{i} \times \mathbf{R}^{2})}^{2} (1 + \| u \|_{\tilde{X}_{i}(G_{\alpha}^{i} \times \mathbf{R}^{2})}^{4}).
\endaligned
\end{equation}

\noindent Also by $(\ref{1.12.1})$ - $(\ref{1.14})$, and $(\ref{2.26})$ - $(\ref{2.28})$,

\begin{equation}\label{13.26}
\sum_{0 \leq l_{2} \leq i - 10} 2^{l_{2} - i} \sum_{l_{2} - 5 \leq j_{1} \leq j_{2} \leq i} \| (\overline{P_{\xi(t), j_{1}} u})(P_{\xi(t), \leq l_{2} - 5} u) \|_{L_{t,x}^{2}(G_{\alpha}^{i} \times \mathbf{R}^{2})} \| (\overline{P_{\xi(t), j_{2}} u})(P_{\xi(t), \leq l_{2} - 5} u) \|_{L_{t,x}^{2}(G_{\alpha}^{i} \times \mathbf{R}^{2})} 
\end{equation}

\begin{equation}\label{13.27}
\aligned
\lesssim \sum_{0 \leq l_{2} \leq i - 10} 2^{l_{2} - i} \sum_{l_{2} - 5 \leq j_{1} \leq j_{2} \leq i} (\sum_{G_{\beta}^{j_{1}} \subset G_{\alpha}^{i}} \| (\overline{P_{\xi(G_{\beta}^{j_{1}}), j_{1} - 2 \leq \cdot \leq j_{1} + 2} u})(P_{\xi(t), \leq l_{2} - 15} u) \|_{L_{t,x}^{2}(G_{\beta}^{j_{1}} \times \mathbf{R}^{2})}^{2})^{1/2} \\
\times (\sum_{G_{\beta'}^{j_{2}} \subset G_{\alpha}^{i}} \| (\overline{P_{\xi(G_{\beta'}^{j_{2}}), j_{2} - 2 \leq \cdot \leq j_{2} + 2} u})(P_{\xi(t), \leq l_{2} - 15} u) \|_{L_{t,x}^{2}(G_{\beta'}^{j_{2}} \times \mathbf{R}^{2})}^{2})^{1/2},
\endaligned
\end{equation}

\noindent so by theorem $\ref{t3.2.1}$,

\begin{equation}\label{13.28}
\aligned
\lesssim \sum_{0 \leq l_{2} \leq i - 10} 2^{l_{2} - i} \sum_{l_{2} - 5 \leq j_{1} \leq j_{2} \leq i} (\sum_{G_{\beta}^{j_{1}} \subset G_{\alpha}^{i}} \| P_{\xi(G_{\beta}^{j_{1}}), j_{1} - 2 \leq \cdot \leq j_{1} + 2} u \|_{U_{\Delta}^{2}(G_{\beta}^{j_{1}} \times \mathbf{R}^{2})}^{2})^{1/2} \\ \times (\sum_{G_{\beta'}^{j_{2}} \subset G_{\alpha}^{i}} \| P_{\xi(G_{\beta'}^{j_{2}}), j_{2} - 2 \leq \cdot \leq j_{2} + 2} u \|_{U_{\Delta}^{2}(G_{\beta'}^{j_{2}} \times \mathbf{R}^{2})}^{2})^{1/2} (1 + \| u \|_{\tilde{X}_{i}(G_{\alpha}^{i} \times \mathbf{R}^{2})}^{4}) \lesssim (1 + \| u \|_{\tilde{X}_{i}(G_{\alpha}^{i} \times \mathbf{R}^{2})}^{6}).
\endaligned
\end{equation}



\noindent Therefore,

\begin{equation}\label{13.30}
\sum_{0 \leq l_{2} \leq i - 10} 2^{l_{2} - i} \| F_{2} \|_{L_{t,x}^{1}(G_{\alpha}^{i} \times \mathbf{R}^{2})} \lesssim (1 + \| u \|_{\tilde{X}_{i}(G_{\alpha}^{i} \times \mathbf{R}^{2})})^{6}.
\end{equation}

\noindent Interpolating $(\ref{13.19.1})$, $(\ref{13.20})$, and $(\ref{13.30})$,

\begin{equation}\label{13.31}
\sum_{0 \leq l_{2} \leq i - 10} 2^{l_{2} - i} \| F_{3} \|_{L_{t,x}^{1}(G_{\alpha}^{i} \times \mathbf{R}^{2})} \lesssim (1 + \| u \|_{\tilde{X}_{i}(G_{\alpha}^{i} \times \mathbf{R}^{2})})^{6}.
\end{equation}

\noindent Therefore, collecting $(\ref{13.19.1})$, $(\ref{13.20})$, and $(\ref{13.30})$,

\begin{equation}\label{13.32}
\sum_{0 \leq l_{2} \leq i - 10} 2^{l_{2} - 2i} \int_{G_{\alpha}^{i}} [F_{2} + F_{3} + F_{4}](t,y) \frac{(x - y)}{|x - y|} \cdot Im[\bar{v} (\nabla - i \xi(t)) v](t,x) dx dy dt \lesssim \| v_{0} \|_{L^{2}}^{2} (1 + \| u \|_{\tilde{X}_{i}(G_{\alpha}^{i} \times \mathbf{R}^{2})}^{6}).
\end{equation}

\noindent Now

\begin{equation}\label{13.33}
F_{0}(t, y) = Im[(\overline{P_{\xi(t), \leq l_{2}} u_{l}}) P_{\xi(t), \leq l_{2}} F(u_{l})] = Im[|u_{l}|^{4}] = 0,
\end{equation}

\noindent and

\begin{equation}\label{13.34}
\aligned
F_{1}(t, y) = Im[(\overline{P_{\xi(t), \leq l_{2}} u_{l}}) P_{\xi(t), \leq l_{2}} (2 |u_{l}|^{2} u_{h} + u_{l}^{2} \bar{u}_{h})] + Im[(\overline{P_{\xi(t), \leq l_{2}} u_{h}}) P_{\xi(t), \leq l_{2}} F(u_{l})] \\
= Im[\overline{u_{l}} P_{\xi(t), \leq l_{2}} (2 |u_{l}|^{2} u_{h} + u_{l}^{2} \overline{u_{h}}) + (\overline{P_{\xi(t), \leq l_{2}} u_{h}}) F(u_{l})].
\endaligned
\end{equation}

\begin{equation}\label{13.35}
\aligned
Im[\overline{u_{l}} P_{\xi(t), \leq l_{2}} (2 |u_{l}|^{2} (P_{\xi(t), \leq l_{2} - 2} u_{h}) + u_{l}^{2} (\overline{P_{\xi(t), \leq l_{2} - 2} u_{h}})) + \overline{(P_{\xi(t), \leq l_{2}} P_{\xi(t), \leq l_{2} - 2} u_{h})} F(u_{l})] \\
= Im[2 F(u_{l}) (\overline{P_{\xi(t), \leq l_{2} - 2} u_{h}}) + 2 F(\bar{u}_{l}) (P_{\xi(t), \leq l_{2} - 2} u_{h})] = 0.
\endaligned
\end{equation}

\noindent Therefore,

\begin{equation}\label{13.36}
F_{0}(t, y) + F_{1}(t, y) = O((P_{\xi(t), \geq l_{2} - 2} u_{h}) |u_{l}|^{2} \overline{u_{l}} + (\overline{P_{\xi(t), \geq l_{2} - 2} u_{h}}) |u_{l}|^{2} u_{l}),
\end{equation}

\noindent and is supported on $|\xi| \geq 2^{l_{2} - 4}$. Integrating by parts,

\begin{equation}\label{13.37}
\aligned
2^{l_{2} - 2i} \int_{G_{\alpha}^{i}} \int \int (\frac{\Delta_{y}}{\Delta_{y}} F_{1}(t,y)) \frac{(x - y)}{|x - y|} \cdot Im[\bar{v} (\nabla - i \xi(t)) v](t,x) dx dy dt \\ = 2^{l_{2} - 2i} \int_{G_{\alpha}^{i}} \int \int (\frac{\partial_{k}}{\Delta_{y}} F_{1}(t,y)) [\frac{\delta_{jk}}{|x - y|} + \frac{(x - y)_{j} (x - y)_{k}}{|x - y|^{3}}] Im[\bar{v} (\partial_{j} - i \xi_{j}(t)) v](t,x) dx dy dt.
\endaligned
\end{equation}

\noindent Then by the Hardy - Littlewood - Sobolev inequality and the Fourier support of $F_{1}$,

\begin{equation}\label{13.38}
\aligned
\lesssim 2^{l_{2} - i} \| v \|_{L_{t}^{6} L_{x}^{3}(G_{\alpha}^{i} \times \mathbf{R}^{2})}^{2} \| |\nabla - i \xi(t)|^{-1} P_{\xi(t), \geq l_{2} - 2} u \|_{L_{t}^{3} L_{x}^{6}(G_{\alpha}^{i} \times \mathbf{R}^{2})} \| u_{l} \|_{L_{t}^{9} L_{x}^{9/2}(G_{\alpha}^{i} \times \mathbf{R}^{2})}^{3} \\
\lesssim 2^{l_{2} - i} \| v_{0} \|_{L^{2}}^{2} \| |\nabla - i \xi(t)|^{-1} P_{\xi(t), \geq l_{2} - 2} u \|_{L_{t}^{3} L_{x}^{6}(G_{\alpha}^{i} \times \mathbf{R}^{2})} \| P_{\xi(t), \leq l_{2}} u \|_{L_{t}^{9} L_{x}^{9/2}(G_{\alpha}^{i} \times \mathbf{R}^{2})}^{3}.
\endaligned
\end{equation}

\begin{equation}\label{13.39}
\sum_{0 \leq l_{2} \leq i - 10} 2^{l_{2} - i} \| v_{0} \|_{L^{2}}^{2} \| |\nabla - i \xi(t)|^{-1} P_{\xi(t), \geq l_{2} - 2} u \|_{L_{t}^{3} L_{x}^{6}(G_{\alpha}^{i} \times \mathbf{R}^{2})} \| P_{\xi(t), \leq l_{2} - 5} u \|_{L_{t}^{9} L_{x}^{9/2}(G_{\alpha}^{i} \times \mathbf{R}^{2})}^{3}
\end{equation}

\begin{equation}\label{13.40}
\lesssim \sum_{0 \leq l_{2} \leq i - 10} 2^{l_{2} - i}  \| |\nabla - i \xi(t)|^{-1} P_{\xi(t), \geq l_{2} - 2} u \|_{L_{t}^{3} L_{x}^{6}} \sum_{l_{5} \leq l_{4} \leq l_{3} \leq l_{2} - 5} \| P_{\xi(t), l_{5}} u \|_{L_{t,x}^{\infty}} \| P_{\xi(t), l_{4}} u \|_{L_{t}^{\infty} L_{x}^{2}} \| P_{\xi(t), l_{3}} u \|_{L_{t}^{3} L_{x}^{6}},
\end{equation}

\noindent so by the Sobolev embedding theorem,

\begin{equation}\label{13.41}
\aligned
\lesssim \sum_{0 \leq l_{2} \leq i - 10} \sum_{0 \leq l_{3} \leq l_{2} - 5} \sum_{l_{1} \geq l_{2} - 2} 2^{l_{3} - i} 2^{l_{2} - l_{1}} \| P_{\xi(t), l_{3}} u \|_{U_{\Delta}^{2}(G_{\alpha}^{i} \times \mathbf{R}^{2})} \| P_{\xi(t), l_{1}} u \|_{U_{\Delta}^{2}(G_{\alpha}^{i} \times \mathbf{R}^{2})} \\ \lesssim \sum_{0 \leq l_{3} \leq l_{1}; l_{3} \leq i - 15} \inf(2^{l_{1}}, 2^{i}) 2^{l_{3} - i} 2^{-l_{1}} \| P_{\xi(t), l_{3}} u \|_{U_{\Delta}^{2}(G_{\alpha}^{i} \times \mathbf{R}^{2})} \| P_{\xi(t), l_{1}} u \|_{U_{\Delta}^{2}(G_{\alpha}^{i} \times \mathbf{R}^{2})} \lesssim \| u \|_{\tilde{X}_{i}(G_{\alpha}^{i} \times \mathbf{R}^{2})}^{2}.
\endaligned
\end{equation}

\noindent This completes the proof of theorem $\ref{t3.2.2}$. $\Box$\vspace{5mm}

\noindent Theorem $\ref{t3.2.2}$ then implies theorem $\ref{tbigtheorem}$, which implies theorem $\ref{t3.0}$, which in turn implies theorem $\ref{t3.1}$. $\Box$

\section{Rigidity}
\noindent The proof of theorem $\ref{t0.9}$ follows directly from the long time Strichartz estimates in the previous section. We first show that if $u$ is an almost periodic solution to $(\ref{0.1})$ in the form of theorem $\ref{t1.9}$ with $\int_{0}^{\infty} N(t)^{3} dt = K < \infty$, then $u \equiv 0$. 

\begin{theorem}\label{t4.1}
 If $u$ is an almost periodic solution in the form of theorem $\ref{t1.9}$, and $\int_{0}^{\infty} N(t)^{3} dt = K < \infty$, then 

\begin{equation}\label{4.0}
 \| u(t) \|_{\dot{H}^{3}(\mathbf{R}^{2})} \lesssim K^{3}.
\end{equation}

\end{theorem}

\noindent \emph{Proof:} Let $[0, T]$ be an interval such that for some $k_{0} \in \mathbf{Z}_{+}$,

\begin{equation}\label{4.1}
 \int_{0}^{T} \int |u(t,x)|^{4} dx dt = 2^{k_{0}}.
\end{equation}

\noindent Let $u_{\lambda}(t,x) = \lambda u(\lambda^{2} t, \lambda x)$, with $\lambda = \frac{\epsilon_{3} 2^{k_{0}}}{\int_{0}^{T} N(t)^{3} dt}$. Then for $N_{\lambda}(t) = \lambda N(\lambda^{2} t)$,

\begin{equation}\label{4.2}
\int_{0}^{\frac{T}{\lambda^{2}}} N_{\lambda}(t)^{3} dt = \epsilon_{3} 2^{k_{0}}.
\end{equation}

\noindent By theorems $\ref{t3.1}$ and $\ref{t3.0}$,

\begin{equation}\label{4.3}
 \| u_{\lambda} \|_{\tilde{X}_{k_{0}}([0, \frac{T}{\lambda^{2}}] \times \mathbf{R}^{2})} \leq C_{0},
\end{equation}

\noindent and

\begin{equation}\label{4.4}
 \| u_{\lambda} \|_{\tilde{Y}_{k_{0}}([0, \frac{T}{\lambda^{2}}] \times \mathbf{R}^{2})} \leq \epsilon_{2}^{1/2}.
\end{equation}

\noindent By $(\ref{2.20})$, $|\xi_{\lambda}(t)| \leq 2^{-20} \epsilon_{3} \epsilon_{1}^{-1/2} 2^{k_{0}} << 2^{k_{0}}$ for all $t \in [0, \frac{T}{\lambda^{2}}]$.\vspace{5mm}

\noindent \textbf{Remark:} $N(0) = 1$ so we only need to be concerned about $K \geq 1$. Also since $N(t) \leq 1$ on $[0, \infty)$, after rescaling $N_{\lambda}(t) \leq \epsilon_{3} 2^{k_{0}}$ for $t \in [0, \infty)$.\vspace{5mm}



\noindent Duhamel's principle combined with theorems $\ref{t3.2}$ and $\ref{tbigtheorem}$, implies that

\begin{equation}\label{4.5}
 \| P_{> N} u_{\lambda} \|_{U_{\Delta}^{2}([0, \frac{T}{\lambda^{2}}] \times \mathbf{R}^{2})} \lesssim \inf_{t \in [0, \frac{T}{\lambda^{2}}]} \| P_{> N} u_{\lambda}(t) \|_{L_{x}^{2}(\mathbf{R}^{2})} + \epsilon_{2}^{1/3} C_{0}^{3} \| P_{> \frac{N}{64}} u_{\lambda} \|_{U_{\Delta}^{2}([0, \frac{T}{\lambda^{2}}] \times \mathbf{R}^{2})}.
\end{equation}

\noindent Since the $U_{\Delta}^{2}$ and $L^{2}$ norms are scale invariant, $(\ref{4.5})$ implies that for $N \geq \epsilon_{3}^{-1} K$,

\begin{equation}\label{4.6}
\| P_{> N} u(t) \|_{U_{\Delta}^{2}([0, T] \times \mathbf{R}^{2})} \lesssim \inf_{t \in [0, T]} \| P_{> N} u(t) \|_{L_{x}^{2}(\mathbf{R}^{2})} + \epsilon_{2}^{1/3} C_{0}^{3} \| P_{> \frac{N}{64}} u \|_{U_{\Delta}^{2}([0, T] \times \mathbf{R}^{2})}.
\end{equation}

\noindent Now by theorem $\ref{t3.1}$,

\begin{equation}\label{4.7}
 \| P_{> \epsilon_{3}^{-1} K} u(t) \|_{U_{\Delta}^{2}([0, \infty) \times \mathbf{R}^{2})} \leq C_{0}.
\end{equation}

\noindent Also, $\int_{0}^{\infty} N(t)^{3} dt = K < \infty$ combined with $(\ref{2.20})$ implies that $\lim_{t \rightarrow \infty} N(t) = 0$ and $|\xi(t)| \leq 2^{-20} \epsilon_{3}^{-1/2} K$, so for $N \geq \epsilon_{3}^{-1} K$,

\begin{equation}\label{4.8}
 \lim_{t \rightarrow \infty} \| P_{> N} u(t) \|_{L^{2}(\mathbf{R}^{2})} = 0.
\end{equation}

\noindent Choosing $\epsilon_{2}, \epsilon_{3} > 0$ sufficiently small and satisfying $(\ref{2.17.13})$, in particular $\epsilon_{2}^{1/3} C(u) C_{0}^{3} << 1$, $(\ref{4.6})$ implies that for $N \geq \epsilon_{3}^{-1} K$,

\begin{equation}\label{4.9}
\| P_{> N} u(t) \|_{U_{\Delta}^{2}([0, \infty) \times \mathbf{R}^{2})} \leq 2^{-20} \| P_{> \frac{N}{64}} u \|_{U_{\Delta}^{2}([0, \infty) \times \mathbf{R}^{2})}.
\end{equation}

\noindent Therefore, by conservation of mass, $(\ref{4.7})$, and induction on $N$,

\begin{equation}\label{4.10}
 \| u(t) \|_{\dot{H}^{3}(\mathbf{R}^{2})} \lesssim \epsilon_{3}^{-3} K^{3}.
\end{equation}

\noindent $\Box$\vspace{5mm}

\noindent Then by conservation of energy,

\begin{theorem}\label{t4.2}
If $u$ is an almost periodic solution to $(\ref{0.1})$ in the form of theorem $\ref{t1.9}$ and  $\int_{0}^{\infty} N(t)^{3} dt = K < \infty$, then $u \equiv 0$.
\end{theorem}

\noindent \emph{Proof:} By $(\ref{2.20})$ and $\xi(0) = 0$, the limit

\begin{equation}
\xi_{\infty} = \lim_{t \rightarrow +\infty} \xi(t)
\end{equation}

\noindent exists and moreover $|\xi(t)| \leq 2^{-20} \epsilon_{1}^{-1/2} K$ for all $t \in [0, \infty)$, thus $|\xi_{\infty}| \leq 2^{-20} \epsilon_{1}^{-1/2} K$. Then make a Galilean transformation mapping $\xi_{\infty}$ to the origin. After the Galilean transformation,

\begin{equation}\label{4.11}
\| u(t) \|_{\dot{H}^{3}(\mathbf{R}^{2})} \lesssim \epsilon_{3}^{-3} K^{3}.
\end{equation}

\noindent Then $\xi(t) \rightarrow 0$, $N(t) \rightarrow 0$, and $(\ref{1.20})$ implies that

\begin{equation}\label{4.12}
\lim_{t \rightarrow \infty} \| P_{\xi(t), \leq C(\eta) N(t)} u(t) \|_{\dot{H}^{3}(\mathbf{R}^{2})} = 0.
\end{equation}

\noindent Interpolating $(\ref{4.11})$ with $(\ref{1.20})$ also implies

\begin{equation}\label{4.13}
\| P_{\xi(t), \geq C(\eta) N(t)} u(t) \|_{\dot{H}^{1}(\mathbf{R}^{2})} \lesssim \eta^{1/3}.
\end{equation}

\noindent Since $\eta > 0$ can be arbitrarily small, $(\ref{4.12})$ and $(\ref{4.13})$ imply

\begin{equation}\label{4.14}
 \lim_{t \rightarrow +\infty} \| u(t) \|_{\dot{H}^{1}(\mathbf{R}^{2})} = 0.
\end{equation}

\noindent By Sobolev embedding this implies

\begin{equation}\label{4.15}
 \lim_{t \rightarrow +\infty} E(u(t)) = 0,
\end{equation}

\noindent which by conservation of energy proves that $E(u(t)) = 0$, which then implies $u \equiv 0$. $\Box$\vspace{5mm}

\noindent Now suppose that $u$ is as in theorem $\ref{t1.9}$ with $\int_{0}^{\infty} N(t)^{3} dt = \infty$. In this case we use the interaction Morawetz estimate of \cite{PV}. It should be observed that \cite{CGT1} also proved the interaction Morawetz estimate that we use here. We rely on the work of \cite{PV} throughout this paper due to their work on the bilinear estimates as well as their work with Galilean invariance.

\begin{theorem}[Frequency localized interaction Morawetz estimate]\label{t4.2}
If $u$ is an almost periodic solution to $(\ref{0.1})$ on $[0, T]$ with $\int_{0}^{T} N(t)^{3} dt = K$, then

\begin{equation}\label{4.16}
\| |\nabla|^{1/2} |P_{\leq 10 \epsilon_{1}^{-1} K} u(t,x)|^{2} \|_{L_{t,x}^{2}([0, T] \times \mathbf{R}^{2})} \lesssim o(K),
\end{equation}

\noindent where $o(K)$ is a quantity, $\frac{o(K)}{K} \rightarrow 0$ as $K \nearrow \infty$.
\end{theorem}

\noindent \emph{Proof:} Again suppose $[0, T]$ is an interval such that for some integer $k_{0}$,

\begin{equation}\label{4.17}
\int_{0}^{T} \int |u(t,x)|^{4} dx dt = 2^{k_{0}}.
\end{equation}

\noindent Again rescale, with $\lambda = \frac{\epsilon_{3} 2^{k_{0}}}{K}$. Again by theorem $\ref{t3.1}$,



\begin{equation}\label{4.18}
\| u_{\lambda} \|_{\tilde{X}_{k_{0}}([0, \frac{T}{\lambda^{2}}] \times \mathbf{R}^{2})} \lesssim 1.
\end{equation}

\noindent Then let $w = P_{\leq k_{0}} u$.

\begin{equation}\label{4.19}
i w_{t} + \Delta w = F(w) + \mathcal N = F(w) + P_{\leq k_{0}} F(u) - F(w).
\end{equation}

\noindent Let

\begin{equation}\label{4.20}
M_{\omega}(t) = \int \int |w(t,y)|^{2} \frac{(x - y)_{\omega}}{|(x - y)_{\omega}|} \cdot Im[\bar{w} \partial_{\omega} w](t,x) dx dy.
\end{equation}

\noindent Making computations identical to the computations that gave $(\ref{12.15})$ - $(\ref{12.17})$,

\begin{equation}\label{4.21}
\frac{d}{dt} M_{\omega}(t) = \int \int_{x_{\omega} = y_{\omega}} |\partial_{\omega} (\overline{w(t,y)} w(t,x))|^{2} dx dy + \int \int_{x_{\omega} = y_{\omega}} |w(t,y)|^{2} |w(t,x)|^{4} dx dy
\end{equation}

\begin{equation}\label{4.22}
+ 2 \int \int Im[\bar{w} \mathcal N](t,y) \frac{(x - y)_{\omega}}{|(x - y)_{\omega}|} Im[\bar{w} \partial_{\omega} w](t,x) dx dy
\end{equation}

\begin{equation}\label{4.23}
+ \int \int |w(t,y)|^{2} \frac{(x - y)_{\omega}}{|(x - y)_{\omega}|} Im[\bar{\mathcal N} \partial_{\omega} w](t,x) dx dy + \int \int |w(t,y)|^{2} \frac{(x - y)_{\omega}}{|(x - y)_{\omega}|} Im[\bar{w} \partial_{\omega} \mathcal N](t,x) dx dy.
\end{equation}

\noindent Following \cite{PV}, the properties of the Radon transform imply

\begin{equation}\label{4.24}
\int \int \int_{x_{\omega} = y_{\omega}} |\partial_{\omega} (\overline{w(t,y)} w(t,x))|^{2} dx dy d\omega + \int \int \int_{x_{\omega} = y_{\omega}} |w(t,y)|^{2} |w(t,x)|^{4} dx dy d\omega \gtrsim \| |\nabla|^{1/2} |w(t,x)|^{2} \|_{L_{x}^{2}(\mathbf{R}^{2})}^{2}.
\end{equation}

\noindent \textbf{Remark:} Unlike in the previous section, we only consider the defocusing problem and therefore simply ignore the second term in $(\ref{4.21})$. This issue is the main reason the results of this paper do not directly carry over to focusing problem. This difficulty will be addressed in an upcoming paper (\cite{D5}).\vspace{5mm}

\noindent Again by the fundamental theorem of calculus in time, if

\begin{equation}
M(t) = \int M_{\omega}(t) d\omega,
\end{equation}

\begin{equation}\label{4.25}
\| |\nabla|^{1/2} |w(t,x)|^{2} \|_{L_{t,x}^{2}([0, \frac{T}{\lambda^{2}}] \times \mathbf{R}^{2})}^{2} \lesssim \sup_{t \in [0, \frac{T}{\lambda^{2}}]} |M(t)| + |\int_{0}^{\frac{T}{\lambda^{2}}} \int (\ref{4.22}) + (\ref{4.23}) d\omega dt|.
\end{equation}

\noindent As in $(\ref{12.30})$ and $(\ref{12.36})$ - $(\ref{12.38})$,

\begin{equation}\label{4.26}
M(t) = \int \int |w(t,y)|^{2} \frac{(x - y)}{|x - y|} \cdot Im[\bar{w} (\nabla - i \xi(t)) w](t,x) dx dy,
\end{equation}

\noindent and

\begin{equation}\label{4.27}
\int (\ref{4.22}) + (\ref{4.23}) d\omega = 2 \int \int Im[\bar{w} \mathcal N](t,y) \frac{(x - y)}{|x - y|} Im[\bar{w} (\nabla - i \xi(t)) w](t,x) dx dy
\end{equation}

\begin{equation}\label{4.27.1}
+ \int \int |w(t,y)|^{2} \frac{(x - y)}{|x - y|} \cdot Im[\bar{\mathcal N} (\nabla - i \xi(t)) w](t,x) dx dy 
\end{equation}

\begin{equation}\label{4.27.2}
+ \int \int |w(t,y)|^{2} \frac{(x - y)}{|x - y|} \cdot Im[\bar{w} (\nabla - i \xi(t)) \mathcal N](t,x) dx dy.
\end{equation}

\noindent In fact, the Galilean invariance of $(\ref{4.26})$ and $(\ref{4.27})$ would also hold if $\frac{(x - y)_{\omega}}{|(x - y)_{\omega}|}$ were replaced with a time dependent function $a(t, x - y)$, provided that $a(t, x - y)$ was an odd function of $x - y$. This fact will be useful to the focusing problem (in \cite{D5}), we will estimate $(\ref{4.26})$ and $(\ref{4.27})$ when $\frac{x - y}{|x - y|}$ is replaced by $a(t, x - y)$ provided there exists a constant $C$ such that

\begin{enumerate}
\item $|a(t, x - y)| \leq C$, and

\item $|\nabla a(t,x - y)| \leq \frac{C}{|x - y|}$.
\end{enumerate}

\noindent First take $(\ref{4.26})$. Under the conditions described in theorem $\ref{t1.9}$, $N(t) \leq 1$, so $N_{\lambda}(t) \leq \frac{\epsilon_{3} 2^{k_{0}}}{K}$. Therefore, by the Arzela - Ascoli theorem (see $(\ref{1.19})$ and $(\ref{1.20})$), for any $\eta > 0$, if $K(\eta) \geq C(\eta)$, with $C(\eta)$ given in $(\ref{1.20})$, then

\begin{equation}\label{4.28}
\| (\nabla - i \xi(t)) w \|_{L_{t}^{\infty} L_{x}^{2}([0, T] \times \mathbf{R}^{2})} \lesssim \eta 2^{k_{0}}.
\end{equation}

\noindent Therefore, by conservation of mass,

\begin{equation}\label{4.29}
\aligned
\sup_{t \in [0, \frac{T}{\lambda^{2}}]} \int \int |w(t,y)|^{2} a(t, x - y) \cdot Im[\bar{w} \nabla w](t,x) dx dy \\ = \sup_{t \in [0, \frac{T}{\lambda^{2}}]} \int \int |w(t,y)|^{2} a(t, x - y) \cdot Im[\bar{w} (\nabla - i \xi(t)) w](t,x) dx dy \lesssim \eta 2^{k_{0}}.
\endaligned
\end{equation}

\noindent Next take $(\ref{4.27})$. As in $(\ref{13.17})$, let $u_{l} = P_{\leq k_{0} - 3} u$, $u = u_{l} + u_{h}$, and decompose

\begin{equation}\label{4.30}
Im[\bar{w} \mathcal N](t,y) = F_{0}(t,y) + F_{1}(t,y) + F_{2}(t,y) + F_{3}(t,y) + F_{4}(t,y).
\end{equation}

\noindent By $(\ref{13.19.1})$, $(\ref{13.20})$, $(\ref{13.30})$, and $(\ref{13.31})$,

\begin{equation}\label{4.31}
\| F_{2} + F_{3} + F_{4} \|_{L_{t,x}^{1}([0, \frac{T}{\lambda^{2}}] \times \mathbf{R}^{2})} \lesssim 1,
\end{equation}

\noindent so by $(\ref{4.28})$ and conservation of mass,

\begin{equation}\label{4.32}
\int_{0}^{\frac{T}{\lambda^{2}}} \int \int [F_{2} + F_{3} + F_{4}](t,y) a(t, x - y) \cdot Im[\bar{w} (\nabla - i \xi(t)) w](t,x) dx dy dt \lesssim \eta 2^{k_{0}}.
\end{equation}

\noindent As in $(\ref{13.33})$, $F_{0} = 0$. Finally, as in $(\ref{13.36})$, $F_{1}(t,y)$ is supported on $|\xi| \geq 2^{k_{0} - 4}$ so integrating by parts,

\begin{equation}\label{4.32}
\aligned
\int \int \int F_{1}(t,y) a(t, x - y) \cdot Im[\bar{w} (\nabla - i \xi(t)) w](t,x) dx dy dt \\ = \int \int \int (\frac{\Delta_{y}}{\Delta_{y}} F_{1}(t,y)) a(t, x - y) \cdot Im[\bar{w} (\nabla - i \xi(t)) w](t,x) dx dy dt \\
= \int \int \int (\frac{\partial_{k}}{\Delta_{y}} F_{1}(t,y)) (\partial_{k} a)(t, x - y)_{y} Im[\bar{w} (\partial_{j} - i \xi_{j}(t)) w](t,x) dx dy dt,
\endaligned
\end{equation}

\noindent so by the Hardy - Littlewood - Sobolev inequality and the properties of $a(t, x - y)$,

\begin{equation}\label{4.33}
\lesssim 2^{-k_{0}} \| (\nabla - i \xi(t)) w \|_{L_{t,x}^{4}([0, \frac{T}{\lambda^{2}}] \times \mathbf{R}^{2})} \| w \|_{L_{t}^{\infty} L_{x}^{2}([0, \frac{T}{\lambda^{2}}] \times \mathbf{R}^{2})} \| u_{h} \|_{L_{t,x}^{4}([0, \frac{T}{\lambda^{2}}] \times \mathbf{R}^{2})} \| u_{l} \|_{L_{t, x}^{6}([0, \frac{T}{\lambda^{2}}] \times \mathbf{R}^{2})}^{3}.
\end{equation}

\noindent Now by $(\ref{2.37})$ and $(\ref{4.18})$,

\begin{equation}\label{4.34}
\| (\nabla - i \xi(t)) w \|_{L_{t}^{5/2} L_{x}^{10}([0, \frac{T}{\lambda^{2}}] \times \mathbf{R}^{2})} \lesssim \sum_{0 \leq j \leq k_{0}} 2^{j} 2^{\frac{2}{5} (k_{0} - j)} \lesssim 2^{k_{0}}.
\end{equation}

\noindent Interpolating $(\ref{4.28})$ and $(\ref{4.34})$,

\begin{equation}\label{4.35}
\| (\nabla - i \xi(t)) w \|_{L_{t,x}^{4}([0, \frac{T}{\lambda^{2}}] \times \mathbf{R}^{2})} \lesssim \eta^{3/8} 2^{k_{0}}.
\end{equation}

\noindent Similarly, by the Sobolev embedding theorem and $(\ref{4.18})$,

\begin{equation}\label{4.36}
\| u_{l} \|_{L_{t,x}^{6}([0, \frac{T}{\lambda^{2}}] \times \mathbf{R}^{2})} \lesssim \sum_{0 \leq j \leq k_{0}} 2^{j/3} \| P_{\xi(t), j} u \|_{L_{t}^{6} L_{x}^{3}([0, \frac{T}{\lambda^{2}}] \times \mathbf{R}^{2})} \lesssim \sum_{0 \leq j \leq k_{0}} 2^{j/3} 2^{(k_{0} - j)/6} \lesssim 2^{k_{0}/3}.
\end{equation}

\noindent Therefore, (even with $\frac{x - y}{|x - y|}$ replaced with $a(t, x - y)$), 

\begin{equation}\label{4.36.1}
(\ref{4.27}) \lesssim \eta 2^{k_{0}}.
\end{equation}

\noindent Next, by conservation of mass,

\begin{equation}\label{4.37}
(\ref{4.27.1}) \lesssim \| \mathcal N \|_{L_{t}^{3/2} L_{x}^{6/5}([0, \frac{T}{\lambda^{2}}] \times \mathbf{R}^{2})} \| (\nabla - i \xi(t)) w \|_{L_{t}^{3} L_{x}^{6}([0, \frac{T}{\lambda^{2}}] \times \mathbf{R}^{2})}.
\end{equation}

\noindent As in $(\ref{5.30})$ - $(\ref{5.30.2})$, by $(\ref{4.18})$, $(\ref{4.28})$, and $(\ref{4.34})$,

\begin{equation}\label{4.38}
(\ref{4.37}) \lesssim \| (\nabla - i \xi(t)) w \|_{L_{t}^{3} L_{x}^{6}([0, \frac{T}{\lambda^{2}}] \times \mathbf{R}^{2})} \lesssim \eta^{1/6} 2^{k_{0}}.
\end{equation}

\noindent Finally, integrating by parts,

\begin{equation}\label{4.39}
(\ref{4.27.2}) = (\ref{4.27.1}) - \int \int \int |w(t,y)|^{2} \frac{1}{|x - y|} Re[\bar{w} \mathcal N](t,x) dx dy dt.
\end{equation}

\noindent By $(\ref{10.41})$,

\begin{equation}\label{4.40}
\int \int \int |w(t,y)|^{2} \frac{1}{|x - y|} Re[\bar{w} \mathcal N](t,x) dx dy dt \lesssim \int \int \int |w(t,y)|^{2} \frac{1}{|x - y|} |w(t,x)| |u_{h}(t,x)|^{3} dx dy dt
\end{equation}

\begin{equation}\label{4.41}
+ \int \int \int |w(t,y)|^{2} \frac{1}{|x - y|} |w(t,x)| |u_{h}(t,x)| |u_{l}(t,x)|^{2} dx dy dt.
\end{equation}

\noindent By the Hardy - Littlewood - Sobolev inequality, $(\ref{2.37})$, $(\ref{4.35})$, $(\ref{4.28})$, the Sobolev embedding theorem, conservation of mass, and interpolation,

\begin{equation}\label{4.42}
(\ref{4.40}) \lesssim \| u_{h} \|_{L_{t,x}^{4}}^{3} \| w \|_{L_{t}^{12} L_{x}^{4}}^{3} \lesssim \eta^{3/8} 2^{k_{0}}.
\end{equation}

\noindent By the Hardy - Littlewood - Sobolev inequality and the Sobolev embedding theorem,

\begin{equation}\label{4.43}
(\ref{4.41}) \lesssim \| u_{h} \|_{L_{t}^{3} L_{x}^{6}} \| w \|_{L_{t}^{9} L_{x}^{10/3}}^{3} \| u_{l} \|_{L_{t}^{6} L_{x}^{60/11}}^{2}.  
\end{equation}

\noindent By $(\ref{4.38})$, Bernstein's inequality, and the Fourier support of $u_{h}$,

\begin{equation}\label{4.44}
(\ref{4.43}) \lesssim \eta^{1/6} \| w \|_{L_{t}^{9} L_{x}^{10/3}}^{3} \| u_{l} \|_{L_{t}^{6} L_{x}^{60/11}}^{2}.
\end{equation}

\noindent Next, by $(\ref{2.37})$, the Sobolev embedding theorem, and $(\ref{4.18})$,

\begin{equation}\label{4.45}
\| w \|_{L_{t}^{9} L_{x}^{10/3}} \lesssim \sum_{0 \leq j \leq k_{0}} 2^{j/9 + j/45} \| P_{\xi(t), j} u \|_{L_{t}^{9} L_{x}^{18/7}} \lesssim 2^{2 k_{0}/15}.
\end{equation}

\noindent Also,

\begin{equation}\label{4.46}
\| u_{l} \|_{L_{t}^{6} L_{x}^{60/11}} \lesssim \sum_{0 \leq j \leq k_{0}} 2^{3j/10} \| P_{\xi(t), j} u \|_{L_{t}^{6} L_{x}^{3}} \lesssim 2^{3k_{0}/10}.
\end{equation}

\noindent Therefore, $(\ref{4.29})$, $(\ref{4.44}) \lesssim \eta^{1/6} 2^{k_{0}}$. Combining $(\ref{4.36.1})$, $(\ref{4.38})$, $(\ref{4.42})$, and $(\ref{4.43})$ - $(\ref{4.46})$ imply $(\ref{4.26}) + (\ref{4.27}) + (\ref{4.27.1}) + (\ref{4.27.2}) \lesssim \eta^{1/6} 2^{k_{0}}$. Therefore, by $(\ref{4.25})$,

\begin{equation}\label{4.47}
\| |\nabla|^{1/2} |w(t,x)|^{2} \|_{L_{t,x}^{2}([0, \frac{T}{\lambda^{2}}] \times \mathbf{R}^{2})}^{2} \lesssim \eta^{1/6} 2^{k_{0}}.
\end{equation}

\noindent Undoing the scaling $u(t,x) \mapsto \lambda u(\lambda^{2} t, \lambda x)$, $\lambda = \frac{\epsilon_{3} 2^{k_{0}}}{K}$,

\begin{equation}\label{4.48}
\| |\nabla|^{1/2} |P_{\leq 10 \epsilon_{1}^{-1} K} u(t,x)|^{2} \|_{L_{t,x}^{2}([0, T] \times \mathbf{R}^{2})}^{2} \lesssim \epsilon_{3}^{-1} \eta(K)^{1/6} K.
\end{equation}

\noindent This proves theorem $\ref{t4.2}$. $\Box$\vspace{5mm}

\noindent Then by the Sobolev embedding theorem,

\begin{equation}\label{4.49}
\| P_{\leq 10 \epsilon_{1}^{-1} K} u(t,x) \|_{L_{t}^{4} L_{x}^{8}([0, T] \times \mathbf{R}^{2})}^{4} \lesssim \eta(K)^{1/6} K.
\end{equation}

\noindent By $(\ref{2.20})$, $(\ref{1.19})$, $(\ref{1.20})$, and the fact that $N(t) \leq 1$ on $[0, \infty)$, if $\| u(t) \|_{L^{2}(\mathbf{R}^{2})} = m_{0}$,

\begin{equation}\label{4.50}
\lim_{K \nearrow \infty} \int_{|x - x(t)| \leq \frac{\eta(K)^{-1/100}}{N(t)}} |P_{\leq 10 \epsilon_{1}^{-1} K} u(t, x)|^{2} dx = m_{0}^{2},
\end{equation}

\noindent uniformly in $t \in [0, \infty)$. Therefore, by H{\"o}lder's inequality,

\begin{equation}\label{4.51}
\aligned
m_{0}^{4} = \lim_{K \nearrow \infty} \frac{1}{K} \int_{0}^{T} N(t)^{3} \| u(t) \|_{L^{2}(\mathbf{R}^{2})}^{4} dt \\ \lesssim \lim_{K \nearrow \infty} \frac{1}{K} \int_{0}^{T} N(t)^{3} (\int_{|x - x(t)| \leq \frac{\eta(K)^{-1/100}}{N(t)}} |P_{\leq 10 \epsilon_{1}^{-1} K} u(t,x)|^{2} dx)^{2} dt \\ \lesssim \lim_{K \nearrow \infty} \frac{\eta(K)^{-3/100}}{K} \int_{0}^{T} \| P_{\leq 10 \epsilon_{1}^{-1} K} u \|_{L_{x}^{8}(\mathbf{R}^{2})}^{4} dt \lesssim \lim_{K \nearrow \infty} \eta(K)^{97/600} = 0.
\endaligned
\end{equation}

\noindent Therefore, $u \equiv 0$. $\Box$

\section{Appendix}
\noindent In the appendix, we pay our debt from section four and prove lemma $\ref{appendix}$. Some of the arguments in this section are very similar to arguments that appear later in section four or in section five. In recognition of the fact that some readers may wish to see the proof of lemma $\ref{appendix}$ immediately, and thus skip to this section, we present the full detail again.\vspace{5mm}

\noindent \emph{Proof of lemma $\ref{appendix}$:} Recalling lemma $\ref{appendix}$, we wish to prove

\begin{equation}\label{14.1}
2^{l_{2} - 2i}  \int_{G_{\beta}^{l_{2}}} \int \int Im[\bar{w} \mathcal N](t,y) \frac{(x - y)}{|x - y|} Im[\bar{v} (\nabla - i \xi(t)) v](t,x) dx dy dt \lesssim \| v_{0} \|_{L_{x}^{2}(\mathbf{R}^{2})}^{2} (1 + \| u \|_{\tilde{X}_{i}(G_{\alpha}^{i} \times \mathbf{R}^{2})}^{4}).
\end{equation}

\noindent First observe from $(\ref{2.20})$ and $(\ref{12.45})$ - $(\ref{12.47})$, that

\begin{equation}\label{14.2}
\aligned
2^{l_{2} - 2i} \int_{G_{\beta}^{l_{2}}} \int \int Im[\bar{w} \mathcal N_{2}](t,y) \frac{(x - y)}{|x - y|} Im[\bar{v} (\nabla - i \xi(t)) v](t,x) dx dy dt \\ \lesssim 2^{i - 2l_{2}} \| v_{0} \|_{L^{2}}^{2} 2^{l_{2} - 2i} (\int_{G_{\beta}^{l_{2}}} |\xi'(t)| \| P_{\xi(t), l_{2} - 3 \leq \cdot \leq l_{2} + 3} u(t) \|_{L^{2}} \| (\nabla - i \xi(t)) P_{\xi(t), \leq l_{2}} u(t) \|_{L^{2}} dt)
\lesssim 2^{l_{2} - i} \| v_{0} \|_{L^{2}}^{2}.
\endaligned
\end{equation}

\noindent Next partition $u = u_{h} + u_{l}$, where $u_{l} = P_{\xi(t), \leq l_{2} - 5} u$. Decompose

\begin{equation}\label{5.34}
Im[\bar{w} \mathcal N_{1}](t,y) = F_{0}(t,y) + F_{1}(t,y) + F_{2}(t,y) + F_{3}(t,y) + F_{4}(t,y),
\end{equation}

\noindent where $F_{0}$ consists of the terms in $Im[\bar{w} \mathcal N_{1}]$ with zero $u_{h}$ terms and four $u_{l}$ terms, $F_{1}$ has one $u_{h}$ term and three $u_{l}$ terms, and so on.

\begin{equation}\label{5.35}
F_{0}(t,y) = Im[(\overline{P_{\xi(t), \leq l_{2}} P_{\xi(t), \leq l_{2} - 5} u}) P_{\xi(t), \leq l_{2}} F(P_{\xi(t), \leq l_{2} - 5} u)] = Im[|u_{l}|^{4}] = 0.
\end{equation}

\noindent Next, by $(\ref{2.37})$ and $(\ref{2.38})$, since $l_{2} \leq i - 10$ and $G_{\beta}^{l_{2}} \subset G_{\alpha}^{i}$,

\begin{equation}
\| F_{3} + F_{4} \|_{L_{t,x}^{1}(G_{\beta}^{l_{2}} \times \mathbf{R}^{2})}\lesssim \| u_{h} \|_{L_{t}^{3} L_{x}^{6}(G_{\beta}^{l_{2}} \times \mathbf{R}^{2})}^{3} \| u \|_{L_{t}^{\infty} L_{x}^{2}(G_{\beta}^{l_{2}} \times \mathbf{R}^{2})} \lesssim \| u \|_{\tilde{X}_{i}(G_{\alpha}^{i} \times \mathbf{R}^{2})}^{3},
\end{equation}

\noindent so

\begin{equation}
2^{l_{2} - 2i} \int_{G_{\beta}^{l_{2}}} \int \int (F_{3} + F_{4})(t,y) \frac{(x - y)}{|x - y|} \cdot Im[\bar{v} (\nabla - i \xi(t)) v](t,x) dx dy dt \lesssim 2^{l_{2} - i} \| v_{0} \|_{L^{2}}^{2} \| u \|_{\tilde{X}_{i}(G_{\alpha}^{i} \times \mathbf{R}^{2})}^{3}.
\end{equation}

\noindent Next,

\begin{equation}
\aligned
F_{1}(t,y) = Im[(\overline{P_{\xi(t), \leq l_{2}} u_{h}}) P_{\xi(t), \leq l_{2}} F(u_{l}) + \overline{u_{l}} P_{\xi(t), \leq l_{2}}(2 |u_{l}|^{2} u_{h} + u_{l}^{2} \overline{u_{h}})] \\
= Im[(\overline{P_{\xi(t), \leq l_{2}} P_{\xi(t), \leq l_{2} - 2} u_{h}}) F(u_{l}) + \overline{u_{l}} P_{\xi(t), \leq l_{2}} (2 |u_{l}|^{2} (P_{\xi(t), \leq l_{2} - 2} u_{h}) + u_{l}^{2} (\overline{P_{\xi(t), \leq l_{2} - 2} u_{h}}))] \\
+ Im[(\overline{P_{\xi(t), \leq l_{2}} P_{\xi(t), \geq l_{2} - 2} u_{h}}) F(u_{l}) + \overline{u_{l}} P_{\xi(t), \leq l_{2}} (2 |u_{l}|^{2} (P_{\xi(t), \geq l_{2} - 2} u_{h}) + u_{l}^{2} (\overline{P_{\xi(t), \geq l_{2} - 2} u_{h}}))] \\
= 2 Im[(\overline{P_{\xi(t), \leq l_{2}} u_{h}}) |u_{l}|^{2} u_{l} + (P_{\xi(t), \leq l_{2}} u_{h}) |u_{l}|^{2} \overline{u_{l}}] \\ 
+ Im[(\overline{P_{\xi(t), \leq l_{2}} P_{\xi(t), \geq l_{2} - 2} u_{h}}) F(u_{l}) + \overline{u_{l}} P_{\xi(t), \leq l_{2}} (2 |u_{l}|^{2} (P_{\xi(t), \geq l_{2} - 2} u_{h}) + u_{l}^{2} (\overline{P_{\xi(t), \geq l_{2} - 2} u_{h}}))] \\
= Im[(\overline{P_{\xi(t), \leq l_{2}} P_{\xi(t), \geq l_{2} - 2} u_{h}}) F(u_{l}) + \overline{u_{l}} P_{\xi(t), \leq l_{2}} (2 |u_{l}|^{2} (P_{\xi(t), \geq l_{2} - 2} u_{h}) + u_{l}^{2} (\overline{P_{\xi(t), \geq l_{2} - 2} u_{h}}))].
\endaligned
\end{equation}

\noindent This implies that $F_{1}(t,y)$ is supported on $|\xi| \geq 2^{l_{2} - 4}$. Integrating by parts in space,

\begin{equation}\label{5.36}
2^{l_{2} - 2i} \int_{G_{\beta}^{l_{2}}} \int \int (\frac{\Delta_{y}}{\Delta_{y}} F_{1}(t,y)) \frac{(x - y)}{|x - y|} \cdot Im[\bar{v} (\nabla - i \xi(t)) v](t,x) dx dy dt
\end{equation}

\begin{equation}\label{5.37}
 = 2^{l_{2} - 2i} \int_{G_{\beta}^{l_{2}}} \int \int \frac{\partial_{k}}{\Delta_{y}} F_{1}(t,y) [\frac{\delta_{jk}}{|x - y|} - \frac{(x - y)_{j} (x - y)_{k}}{|x - y|^{3}}] Im[\bar{v} (\partial_{j} - i \xi_{j}(t)) v](t,x) dx dy dt.
\end{equation}

\noindent Then by the Hardy - Littlewood - Sobolev inequality and the Fourier support of $F_{1}$,

\begin{equation}\label{5.37.1}
\aligned
(\ref{5.37}) \lesssim 2^{-i} \| v \|_{L_{t}^{3} L_{x}^{6}(G_{\beta}^{l_{2}} \times \mathbf{R}^{2})}^{2} \| u_{h} \|_{L_{t}^{3} L_{x}^{6}(G_{\beta}^{l_{2}} \times \mathbf{R}^{2})} \| u_{l} \|_{L_{t}^{9} L_{x}^{9/2}(G_{\beta}^{l_{2}} \times \mathbf{R}^{2})}^{3}.
\endaligned
\end{equation}

\noindent By $(\ref{2.37})$, $(\ref{2.38})$, and the Sobolev embedding theorem,

\begin{equation}\label{5.37.2}
\| u_{l} \|_{L_{t}^{9} L_{x}^{9/2}(G_{\beta}^{l_{2}} \times \mathbf{R}^{2})} \lesssim \sum_{0 \leq l_{3} \leq l_{2} - 5} 2^{l_{3}/3} \| P_{\xi(t), l_{3}} u \|_{L_{t}^{9} L_{x}^{18/7}(G_{\beta}^{l_{2}} \times \mathbf{R}^{2})} \lesssim 2^{l_{2}/3} \| u \|_{\tilde{X}_{i}(G_{\alpha}^{i} \times \mathbf{R}^{2})},
\end{equation}

\noindent which again by $(\ref{2.37})$ and $(\ref{2.38})$ implies

\begin{equation}\label{5.37.3}
(\ref{5.37.1}) \lesssim 2^{l_{2} - i} \| v_{0} \|_{L^{2}}^{2} \| u \|_{\tilde{X}_{i}(G_{\alpha}^{i} \times \mathbf{R}^{2})}^{4}.
\end{equation}

\noindent Finally, we want to estimate

\begin{equation}\label{5.38}
2^{l_{2} - 2i} \int_{G_{\alpha}^{i}} \int \int F_{2}(t,y) \frac{(x - y)}{|x - y|} \cdot Im[\bar{v} (\nabla - i \xi(t)) v](t,x) dx dy dt.
\end{equation}

\noindent As in the analysis of $F_{1}(t,y)$, integrating by parts,

\begin{equation}\label{5.38.1}
2^{l_{2} - 2i} \int_{G_{\beta}^{l_{2}}} \int \int (\frac{\Delta_{y}}{\Delta_{y}} P_{\geq l_{2} - 10} F_{2})(t,y) \frac{(x - y)}{|x - y|} \cdot Im[\bar{v} (\nabla - i \xi(t)) v](t,x) dx dy dt
\end{equation}

\begin{equation}\label{5.38.2}
= 2^{l_{2} - 2i} \int_{G_{\beta}^{l_{2}}} \int \int (\frac{\partial_{k}}{\Delta_{y}} P_{\geq l_{2} - 10} F_{2})(t,y) [\frac{\delta_{jk}}{|x - y|} - \frac{(x - y)_{j} (x - y)_{k}}{|x - y|^{3}}] Im[\bar{v} (\nabla - i \xi(t)) v](t,x) dx dy dt
\end{equation}

\begin{equation}\label{5.38.3}
\lesssim 2^{-i} \| v_{0} \|_{L^{2}}^{2} \| u_{h} \|_{L_{t}^{3} L_{x}^{6}(G_{\beta}^{l_{2}} \times \mathbf{R}^{2})}^{2} \| u_{l} \|_{L_{t}^{\infty} L_{x}^{4}(G_{\beta}^{l_{2}} \times \mathbf{R}^{2})}^{2} \lesssim 2^{l_{2} - i} \| v_{0} \|_{L^{2}}^{2} \| u \|_{\tilde{X}_{i}(G_{\alpha}^{i} \times \mathbf{R}^{2})}^{2}.
\end{equation}

\noindent Next, expanding out $F_{2}(t,y)$,

\begin{equation}\label{5.39}
\aligned
F_{2}(t,y) = Im[2 \overline{u_{l}} P_{\xi(t), \leq l_{2}} (|u_{h}|^{2} u_{l}) + \overline{u_{l}} P_{\xi(t), \leq l_{2}} (u_{h}^{2} \overline{u_{l}}) \\ + 2 (\overline{P_{\xi(t), \leq l_{2}} u_{h}}) P_{\xi(t), \leq l_{2}} (|u_{l}|^{2} u_{h}) + (\overline{P_{\xi(t), \leq l_{2}} u_{h}}) P_{\xi(t), \leq l_{2}} (u_{l}^{2} \overline{u_{h}})].
\endaligned
\end{equation}

\noindent Observe that

\begin{equation}\label{5.39.1}
2 Im[\overline{u_{l}} P_{\xi(t), \leq l_{2}} (|u_{h}|^{2} u_{l})] = 2 Im[\overline{u_{l}} P_{\xi(t), \leq l_{2}}(|u_{h}|^{2} u_{l}) - |u_{l}|^{2} P_{\leq l_{2}} (|u_{h}|^{2})],
\end{equation}

\noindent and

\begin{equation}\label{5.39.1.1}
Im[2 (\overline{P_{\xi(t), \leq l_{2}} u_{h}}) P_{\xi(t), \leq l_{2}} (|u_{l}|^{2} u_{h})] = Im[2 (\overline{P_{\xi(t), \leq l_{2}} u_{h}}) P_{\xi(t), \leq l_{2}} (|u_{l}|^{2} u_{h}) - 2 |P_{\xi(t), \leq l_{2}} u_{h}|^{2} |u_{l}|^{2}].
\end{equation}

\noindent By $(\ref{5.29})$, $(\ref{2.37})$, $(\ref{2.38})$, and the product rule,

\begin{equation}\label{5.39.1.2}
\aligned
2^{l_{2} - i} \| Im[2 (\overline{P_{\xi(t), \leq l_{2}} u_{h}}) P_{\xi(t), \leq l_{2}} (|u_{l}|^{2} u_{h})] \|_{L_{t,x}^{1}(G_{\beta}^{l_{2}} \times \mathbf{R}^{2})} = 2^{l_{2} - i} \| (\ref{5.39}) \|_{L_{t,x}^{1}(G_{\beta}^{l_{2}} \times \mathbf{R}^{2})}  \\
\lesssim 2^{-i} \| u_{h} \|_{L_{t}^{3} L_{x}^{6}(G_{\beta}^{l_{2}} \times \mathbf{R}^{2})}^{2} \| \nabla |u_{l}|^{2} \|_{L_{t}^{3} L_{x}^{3/2}(G_{\beta}^{l_{2}} \times \mathbf{R}^{2})} \\
\lesssim 2^{l_{2} - i} \| u \|_{\tilde{X}_{i}(G_{\beta}^{l_{2}} \times \mathbf{R}^{2})}^{2} \| \nabla (e^{-ix \cdot \xi(t)} u_{l})(\overline{e^{-ix \cdot \xi(t)} u_{l}}) \|_{L_{t}^{3} L_{x}^{3/2}(G_{\beta}^{l_{2}} \times \mathbf{R}^{2})} \\ \lesssim 2^{l_{2} - i} \| u \|_{\tilde{X}_{i}(G_{\beta}^{l_{2}} \times \mathbf{R}^{2})}^{2} \| (\nabla - i \xi(t)) u_{l} \|_{L_{t}^{3} L_{x}^{6}(G_{\beta}^{l_{2}} \times \mathbf{R}^{2})} \lesssim 2^{l_{2} - i} \| u \|_{\tilde{X}_{i}(G_{\beta}^{l_{2}} \times \mathbf{R}^{2})}^{3}.
\endaligned
\end{equation}

\noindent Next, by a computation similar to $(\ref{5.29})$,

\begin{equation}\label{5.39.2}
\phi(\frac{\xi_{1} + \xi_{2} - \xi(t)}{2^{l_{2}}}) - \phi(\frac{\xi_{1}}{2^{l_{2}}}) \lesssim \frac{|\xi_{2} - \xi(t)|}{2^{l_{2}}},
\end{equation}

\noindent so by $(\ref{2.37})$ and $(\ref{2.38})$,

\begin{equation}\label{5.39.3}
\aligned
2^{l_{2} - i}  \| 2 Im[\overline{u_{l}} P_{\xi(t), \leq l_{2}} (|u_{h}|^{2} u_{l})] \|_{L_{t,x}^{1}(G_{\beta}^{l_{2}} \times \mathbf{R}^{2})} \\ \lesssim 2^{-i}  \| u_{h} \|_{L_{t}^{3} L_{x}^{6}(G_{\beta}^{l_{2}} \times \mathbf{R}^{2})}^{2} \| (\nabla - i \xi(t)) u_{l} \|_{L_{t}^{3} L_{x}^{6}(G_{\beta}^{l_{2}} \times \mathbf{R}^{2})} \| u \|_{L_{t}^{\infty} L_{x}^{2}} \lesssim 2^{l_{2} - i} \| u \|_{\tilde{X}_{i}(G_{\alpha}^{i} \times \mathbf{R}^{2})}^{3}.
\endaligned
\end{equation}

\noindent Also, 

\begin{equation}
P_{\xi(t), \leq l_{2}} (u_{l}^{2} \overline{u_{h}}) - (\overline{P_{\xi(t), \leq l_{2}} u_{h}}) u_{l}^{2} = e^{ix \cdot \xi(t)} P_{\leq l_{2}} ((e^{-ix \cdot \xi(t)} u_{l})^{2} (\overline{e^{-ix \cdot \xi(t)} u_{h}})) - e^{ix \cdot \xi(t)} (\overline{P_{\leq l_{2}} e^{-ix \cdot \xi(t)} u_{h}}) (e^{-ix \cdot \xi(t)} u_{l})^{2},
\end{equation}

\noindent so by $(\ref{5.29})$, $(\ref{2.37})$, $(\ref{2.38})$, and the product rule,

\begin{equation}\label{5.39.1.3}
\aligned
2^{l_{2} - i} \| (\overline{P_{\xi(t), \leq l_{2}} u_{h}}) P_{\xi(t), \leq l_{2}}(u_{l}^{2} \overline{u_{h}}) - (\overline{P_{\xi(t), \leq l_{2}} u_{h}})^{2} u_{l}^{2} \|_{L_{t,x}^{1}(G_{\beta}^{l_{2}} \times \mathbf{R}^{2})} \\ \lesssim 2^{-i} \| u_{h} \|_{L_{t}^{3} L_{x}^{6}(G_{\beta}^{l_{2}} \times \mathbf{R}^{2})}^{2} \| (\nabla - i \xi(t)) u_{l} \|_{L_{t}^{3} L_{x}^{6}(G_{\beta}^{l_{2}} \times \mathbf{R}^{2})} \| u \|_{L_{t}^{\infty} L_{x}^{2}(G_{\beta}^{l_{2}} \times \mathbf{R}^{2})} \lesssim 2^{l_{2} - i} \| u \|_{\tilde{X}_{i}(G_{\alpha}^{i} \times \mathbf{R}^{2})}^{3}.
\endaligned
\end{equation}



\noindent Therefore by $(\ref{5.38.1})$ - $(\ref{5.38.3})$, $(\ref{5.39.1.2})$, $(\ref{5.39.1.3})$, and $(\ref{5.39.3})$,

\begin{equation}
2^{l_{2} - 2i} \int_{G_{\beta}^{l_{2}}} \int \int F_{2}(t,y) \frac{(x - y)}{|x - y|} \cdot Im[\bar{v} (\nabla - i \xi(t)) v](t,x) dx dy dt = O(2^{l_{2} - i} \| v_{0} \|_{L^{2}}^{2} (1 + \| u \|_{\tilde{X}_{i}(G_{\alpha}^{i} \times \mathbf{R}^{2})}^{4})
\end{equation}

\begin{equation}\label{5.39.4}
+ 2^{l_{2} - 2i} \int_{G_{\beta}^{l_{2}}} \int \int P_{\leq l_{2} - 10} Im[(\overline{u_{l}})^{2} u_{h}^{2} + (\overline{P_{\xi(t), \leq l_{2}} u_{h}})^{2} u_{l}^{2}](t, y) \frac{(x - y)}{|x - y|} \cdot Im[\bar{v} (\nabla - i \xi(t)) v](t,x) dx dy dt.
\end{equation}

\begin{equation}\label{5.40}
(\ref{5.39.4}) = 2^{l_{2} - 2i} \int_{G_{\beta}^{l_{2}}} \int \int P_{\leq l_{2} - 10} Im[(\overline{u_{l}})^{2} (u_{h}^{2} - (P_{\xi(t), \leq l_{2}} u_{h})^{2})](t, y) \frac{(x - y)}{|x - y|} \cdot Im[\bar{v} (\nabla - i \xi(t)) v](t,x) dx dy dt.
\end{equation}

\noindent Next make a Galilean transformation (see theorem $\ref{t1.6}$). First observe that $(\ref{5.40})$ is obviously translation invariant.

\begin{equation}
\aligned
(\ref{5.40}) = 2^{l_{2} - 2i}  \int_{G_{\beta}^{l_{2}}} \int \int P_{\leq l_{2} - 10} Im[(\overline{u_{l}})^{2} (u_{h}^{2} - (P_{\xi(t), \leq l_{2}} u_{h})^{2})](t, y + 2t \xi(G_{\beta}^{l_{2}})) \times \\ \frac{(x - y)}{|x - y|} \cdot Im[\bar{v} (\nabla - i \xi(t)) v](t, x + 2t \xi(G_{\beta}^{l_{2}})) dx dy dt.
\endaligned
\end{equation}

\noindent Next observe that

\begin{equation}\label{Gal1}
Im[\bar{v} (\nabla - i \xi(t)) v](t, x + 2t \xi(G_{\beta}^{l_{2}})) = Im[(\overline{e^{-ix \cdot \xi(G_{\beta}^{l_{2}})} v}) (\nabla - i \xi(t) + i \xi(G_{\beta}^{l_{2}}))e^{-ix \cdot \xi(G_{\beta}^{l_{2}})} v](t,x + 2t \xi(G_{\beta}^{l_{2}})).
\end{equation}

\noindent Also, observe that

\begin{equation}\label{Gal2}
\aligned
Im[(\overline{u_{l}})^{2} u_{h}^{2}](t, y + 2 t \xi(G_{\beta}^{l_{2}})) = Im[(\overline{e^{-ix \cdot \xi(G_{\beta}^{l_{2}})} u_{l}})^{2} (e^{-ix \cdot \xi(G_{\beta}^{l_{2}})} u_{h})^{2}](t, y + 2 t \xi(G_{\beta}^{l_{2}})) \\
= Im[(\overline{P_{\xi(t) - \xi(G_{\beta}^{l_{2}}), \leq l_{2} - 5} e^{-ix \cdot \xi(G_{\beta}^{l_{2}})} u})^{2} (P_{\xi(t) - \xi(G_{\beta}^{l_{2}}), \geq l_{2} - 5} e^{-ix \cdot \xi(G_{\beta}^{l_{2}})} u)^{2}](t, y + 2t \xi(G_{\beta}^{l_{2}})).
\endaligned
\end{equation}

\noindent Making a similar computation with the $\overline{u_{l}}^{2} (P_{\xi(t), \leq l_{2}} u_{h})^{2}$ term, the Galilean transformation, $(\ref{Gal1})$, and $(\ref{Gal2})$ imply that we can assume that $\xi(G_{\beta}^{l_{2}}) = 0$ in $(\ref{5.40})$. By $(\ref{2.20})$ this implies $|\xi(t)| << 2^{l_{2}}$ for all $t \in G_{\beta}^{l_{2}}$.\vspace{5mm}
 
\noindent To estimate $(\ref{5.40})$ we will make an inverse Fourier transform and integrate by parts in time.

\begin{equation}
\int_{G_{\beta}^{l_{2}}} \int \int P_{\leq l_{2} - 10} Im[(\overline{u_{l}})^{2} (u_{h}^{2} - (P_{\xi(t), \leq l_{2}} u_{h})^{2})] \frac{(x - y)}{|x - y|} \cdot Im[\bar{v} (\nabla - i \xi(t)) v](t,x) dx dy dt
\end{equation}

\begin{equation}\label{5.41}
\aligned
= (2 \pi)^{-1} \int_{G_{\beta}^{l_{2}}} \int \int \frac{(x - y)}{|x - y|} \cdot Im[\overline{v}(t,x) (\nabla - i \xi(t)) v(t,x)] \\
\times [\int \int \int \int \phi(\frac{(\eta_{1} + \eta_{2} + \eta_{3} + \eta_{4})}{2^{l_{2} - 10}}) e^{iy \cdot (\eta_{1} + \eta_{2} + \eta_{3} + \eta_{4})} \hat{u}_{l}(t,\eta_{3}) \hat{u}_{l}(t,\eta_{4}) \hat{\bar{u}}_{h}(t,\eta_{1}) \hat{\bar{u}}_{h}(t,\eta_{2}) \\
\times (1 - \phi(2^{-l_{2}} (\eta_{1} - \xi(t))) \phi(2^{-l_{2}} (\eta_{2} - \xi(t)))) d\eta_{1} d\eta_{2} d\eta_{3} d\eta_{4}] dx dy dt.
\endaligned
\end{equation}

\noindent Let 

\begin{equation}
q(\eta) = |\eta_{1}|^{2} + |\eta_{2}|^{2} - |\eta_{3}|^{2} - |\eta_{4}|^{2}.
\end{equation}

\noindent Since $|\eta_{3}|, |\eta_{4}| \leq 2^{l_{2} - 4}$, $||\eta_{3}|^{2} + |\eta_{4}|^{2}| \leq 2^{2l_{2} - 8}$. Also, on the support of $(1 - \phi(2^{-l_{2}} (\eta_{1} - \xi(t))) \phi(2^{-l_{2}} (\eta_{2} - \xi(t)))) \phi(\frac{( \eta_{1} + \eta_{2} + \eta_{3} + \eta_{4})}{2^{l_{2} - 10}})$, $|\eta_{1}|^{2} + |\eta_{2}|^{2} > 2^{2l_{2}}$, and $q(\eta) > 2^{2l_{2} - 2}$. Therefore, on the support of $ (1 - \phi(2^{-l_{2}} (\eta_{1} - \xi(t))) \phi(2^{-l_{2}} (\eta_{2} - \xi(t)))) \phi(2^{-l_{2} + 4} \eta_{3}) \phi(2^{-l_{2} + 4} \eta_{4})$,




\begin{equation}\label{5.41.2}
\frac{1}{q(\eta)} = \frac{1}{|\eta_{1}|^{2} + |\eta_{2}|^{2}} \sum_{j = 0}^{\infty} (-1)^{j} (\frac{|\eta_{3}|^{2} + |\eta_{4}|^{2}}{|\eta_{1}|^{2} + |\eta_{2}|^{2}})^{j}
\end{equation}

\noindent converges, and by theorem $C - M$ of \cite{GMS}, $\frac{1}{q(\eta)}$ is a convergent sum of terms whose operator norm is $\lesssim \frac{1}{\eta_{1}^{2} + \eta_{2}^{2}} \sim \frac{1}{|\eta_{1}| |\eta_{2}|}$ (since $|\eta_{3} + \eta_{4}| \leq 2^{l_{2} - 3}$ and $|\eta_{1}|, |\eta_{2}| \gtrsim 2^{l_{2}}$). Integrating by parts in time, if $G_{\beta}^{l_{2}} = [a, b]$,

\begin{equation}
 \aligned
\int_{G_{\beta}^{l_{2}}} \int \int \int \int \int \int \frac{1}{i q(\eta)} (\frac{d}{dt} e^{it q(\eta)}) \frac{(x - y)}{|x - y|} \cdot Im [\bar{v}(t,x) (\nabla - i \xi(t)) v(t,x)]  \\ \times \phi(\frac{(\eta_{1} + \eta_{2} + \eta_{3} + \eta_{4})}{2^{l_{2} - 10}})  e^{iy \cdot (\eta_{1} + \eta_{2} + \eta_{3} + \eta_{4})}  (1 - \phi(2^{-l_{2}} (\eta_{1} - \xi(t))) \phi(2^{-l_{2}} (\eta_{2} - \xi(t)))) \\ \times [e^{-it |\eta_{1}|^{2}} \hat{\bar{u}}_{h}(t,\eta_{1}) e^{-it |\eta_{2}|^{2}} \hat{\bar{u}}_{h}(t,\eta_{2}) e^{it |\eta_{3}|^{2}} \hat{u}_{l}(t,\eta_{3}) e^{it |\eta_{4}|^{2}} \hat{u}_{l}(t,\eta_{4})] d\eta_{1} d\eta_{2} d\eta_{3} d\eta_{4} dx dy dt,
\endaligned
\end{equation}

\begin{equation}\label{5.42}
\aligned
= \int \int \int \int \int \int e^{it q(\eta)} \frac{1}{iq(\eta)} \frac{(x - y)}{|x - y|} \cdot Im [\bar{v}(t,x) (\nabla - i \xi(t)) v(t,x)] \\ \times \phi(\frac{(\eta_{1} + \eta_{2} + \eta_{3} + \eta_{4})}{2^{l_{2} - 10}})  e^{iy \cdot (\eta_{1} + \eta_{2} + \eta_{3} + \eta_{4})}  (1 - \phi(2^{-l_{2}} (\eta_{1} - \xi(t))) \phi(2^{-l_{2}} (\eta_{2} - \xi(t))))  \\  \times [e^{-it |\eta_{1}|^{2}} \hat{\bar{u}}_{h}(t,\eta_{1}) e^{-it |\eta_{2}|^{2}} \hat{\bar{u}}_{h}(t,\eta_{2}) e^{it |\eta_{3}|^{2}} \hat{u}_{l}(t,\eta_{3}) e^{it |\eta_{4}|^{2}} \hat{u}_{l}(t,\eta_{4})] d\eta_{1} d\eta_{2} d\eta_{3} d\eta_{4} dx dy |_{a}^{b}
\endaligned
\end{equation}

\begin{equation}\label{5.43}
\aligned
- \int_{a}^{b} \int \int \int \int \int \int \frac{1}{iq(\eta)} e^{it q( \eta)} \frac{\partial}{\partial t} \frac{(x - y)}{|x - y|} \cdot Im [\bar{v}(t,x) (\nabla - i \xi(t)) v(t,x)]  \\ \times \phi(\frac{(\eta_{1} + \eta_{2} + \eta_{3} + \eta_{4})}{2^{l_{2} - 10}})  e^{iy \cdot (\eta_{1} + \eta_{2} + \eta_{3} + \eta_{4})}  (1 - \phi(2^{-l_{2}} (\eta_{1} - \xi(t))) \phi(2^{-l_{2}} (\eta_{2} - \xi(t)))) \\ \times [e^{-it |\eta_{1}|^{2}} \hat{\bar{u}}_{h}(t,\eta_{1}) e^{-it |\eta_{2}|^{2}} \hat{\bar{u}}_{h}(t,\eta_{2}) e^{it |\eta_{3}|^{2}} \hat{u}_{l}(t,\eta_{3}) e^{it |\eta_{4}|^{2}} \hat{u}_{l}(t,\eta_{4})] d\eta_{1} d\eta_{2} d\eta_{3} d\eta_{4} dx dy dt
\endaligned
\end{equation}

\begin{equation}\label{5.43.1}
\aligned
- \int_{a}^{b} \int \int \int \int \int \int \frac{1}{iq(\eta)} e^{it q( \eta)} \frac{(x - y)}{|x - y|} \cdot Im [\bar{v}(t,x) (\nabla - i \xi(t)) v(t,x)]  \\ \times \frac{\partial}{\partial t} [\phi(\frac{(\eta_{1} + \eta_{2} + \eta_{3} + \eta_{4})}{2^{l_{2} - 10}})  e^{iy \cdot (\eta_{1} + \eta_{2} + \eta_{3} + \eta_{4})}  (1 - \phi(2^{-l_{2}} (\eta_{1} - \xi(t))) \phi(2^{-l_{2}} (\eta_{2} - \xi(t))))] \\ \times [e^{-it |\eta_{1}|^{2}} \hat{\bar{u}}_{h}(t,\eta_{1}) e^{-it |\eta_{2}|^{2}} \hat{\bar{u}}_{h}(t,\eta_{2}) e^{it |\eta_{3}|^{2}} \hat{u}_{l}(t,\eta_{3}) e^{it |\eta_{4}|^{2}} \hat{u}_{l}(t,\eta_{4})] d\eta_{1} d\eta_{2} d\eta_{3} d\eta_{4} dx dy dt
\endaligned
\end{equation}

\begin{equation}\label{5.43.2}
\aligned
- \int_{a}^{b} \int \int \int \int \int \int \frac{1}{i q(\eta)} e^{itq(\eta)} \frac{(x - y)}{|x - y|} \cdot Im[\bar{v} (\nabla - i \xi(t)) v] (t,x) \\ \times \int \int \int \int \phi(\frac{(\eta_{1} + \eta_{2} + \eta_{3} + \eta_{4})}{2^{l_{2} - 10}})  e^{iy \cdot (\eta_{1} + \eta_{2} + \eta_{3} + \eta_{4})}   (1 - \phi(2^{-l_{2}} (\eta_{1} - \xi(t))) \phi(2^{-l_{2}} (\eta_{2} - \xi(t)))) \\ \frac{\partial}{\partial t}[ e^{-it |\eta_{1}|^{2}} \hat{\bar{u}}_{h}(t, \eta_{1}) e^{-it |\eta_{2}|^{2}} \hat{\bar{u}}_{h}(t, \eta_{2}) e^{it |\eta_{3}|^{2}} \hat{u}_{l}(t, \eta_{3}) e^{it |\eta_{4}|^{2}} \hat{u}_{l}(t, \eta_{4})] d\eta_{1} d\eta_{2} d\eta_{3} d\eta_{4} dx dy dt.
\endaligned
\end{equation}

\noindent \textbf{Remark:} To simplify notation, $L_{t}^{p} L_{x}^{q}$ will always indicate $L_{t}^{p} L_{x}^{q}(G_{\beta}^{l_{2}} \times \mathbf{R}^{2})$ unless otherwise indicated.\vspace{5mm}

\noindent By H{\"o}lder's inequality in Fourier space, conservation of mass, and Plancherel's theorem,

\begin{equation}
2^{l_{2} - 2i} (\ref{5.42}) \lesssim 2^{-2l_{2}} 2^{2l_{2}} 2^{i} 2^{l_{2} - 2i} \| v_{0} \|_{L^{2}}^{2} = 2^{l_{2} - i} \| v_{0} \|_{L^{2}}^{2} \lesssim \| v_{0} \|_{L^{2}}^{2}.
\end{equation}

\noindent Next take $(\ref{5.43})$.

\begin{equation}
\aligned
\frac{\partial}{\partial t} Im[\bar{v}(t,x) (\nabla - i \xi(t)) v(t,x)] \\ = \xi'(t) |v(t,x)|^{2} + Re[\bar{v}(t,x) (\nabla - i \xi(t)) \Delta v(t,x)] - Re[\Delta \bar{v}(t,x) (\nabla - i \xi(t)) v(t,x)] \\ = \xi_{j}'(t) |v(t,x)|^{2} + \partial_{k} Re[\bar{v}(t,x) (\nabla_{j} - i \xi_{j}(t)) \partial_{k} v(t,x)] - \partial_{k} Re[\partial_{k} \bar{v}(t,x) (\partial_{j} - i \xi_{j}(t)) v(t,x)].
\endaligned
\end{equation}

\noindent Therefore, taking the leading order term in $(\ref{5.41.2})$,

\begin{equation}\label{5.43.3}
(\ref{5.43}) = \int_{G_{\beta}^{l_{2}}} \int \int Im[(u_{l}^{2}) ((\frac{1}{|\nabla|} u_{h})^{2} - (P_{\xi(t), \leq l_{2}} \frac{1}{|\nabla|} u_{h})^{2})](t,y) \frac{(x - y)}{|x - y|} \cdot \xi'(t) |v(t,x)|^{2} dx dy dt
\end{equation}

\begin{equation}\label{5.43.4}
+ \int_{G_{\beta}^{l_{2}}} \int \int Im[(u_{l}^{2}) ((\frac{1}{|\nabla|} u_{h})^{2} - (P_{\xi(t), \leq l_{2}} \frac{1}{|\nabla|} u_{h})^{2})](t,y) \frac{(x - y)_{j}}{|x - y|} \partial_{k} Re[\overline{v(t,x)} (\nabla_{j} - i \xi_{j}(t)) \partial_{k} v(t,x)] dx dy dt
\end{equation}

\begin{equation}\label{5.43.5}
- \int_{G_{\beta}^{l_{2}}} \int \int Im[(u_{l}^{2}) ((\frac{1}{|\nabla|} u_{h})^{2} - (P_{\xi(t), \leq l_{2}} \frac{1}{|\nabla|} u_{h})^{2})](t,y) \frac{(x - y)_{j}}{|x - y|}  \partial_{k} Re[\overline{\partial_{k} v(t,x)} (\partial_{j} - i \xi_{j}(t)) v(t,x)].
\end{equation}

\noindent By $(\ref{2.20})$ and the Sobolev embedding theorem,

\begin{equation}
(\ref{5.43.3}) \lesssim \| \frac{1}{|\nabla|} u_{h} \|_{L_{t}^{\infty} L_{x}^{2}(G_{\beta}^{l_{2}} \times \mathbf{R}^{2})}^{2} \| u_{l} \|_{L_{t,x}^{\infty}(G_{\beta}^{l_{2}} \times \mathbf{R}^{2})}^{2} \| v(t,x) \|_{L_{t}^{\infty} L_{x}^{2}(G_{\beta}^{l_{2}} \times \mathbf{R}^{2})}^{2} (\int_{G_{\beta}^{l_{2}}} |\xi'(t)| dt) \lesssim \epsilon_{3} \epsilon_{1}^{-1} 2^{l_{2}} \| v_{0} \|_{L^{2}}^{2}.
\end{equation}

\noindent Integrating $(\ref{5.43.4})$ and $(\ref{5.43.5})$ by parts in space,

\begin{equation}
\aligned
(\ref{5.43.4}) = -  \int_{G_{\beta}^{l_{2}}} \int \int Im[(u_{l}^{2}) &((\frac{1}{|\nabla|} u_{h})^{2} - (P_{\xi(t), \leq l_{2}} \frac{1}{|\nabla|} u_{h})^{2})](t,y) \\ &[\frac{\delta_{jk}}{|x - y|} - \frac{(x - y)_{j} (x - y)_{k}}{|x - y|^{3}}] Re[\bar{v}(t,x) (\nabla_{j} - i \xi_{j}(t)) \partial_{k} v(t,x)] dx dy dt,
\endaligned
\end{equation}

\noindent so by the Hardy - Littlewood - Sobolev lemma, $(\ref{2.37})$, and $(\ref{2.38})$, the Sobolev embedding theorem, the fact that $G_{\beta}^{l_{2}} \subset G_{\alpha}^{i}$, and $l_{2} \leq i$,

\begin{equation}
\lesssim 2^{2i} \| \frac{1}{|\nabla|} u_{h} \|_{L_{t}^{3} L_{x}^{6}(G_{\beta}^{l_{2}} \times \mathbf{R}^{2})}^{2} \| u_{l} \|_{L_{t}^{\infty} L_{x}^{4}(G_{\beta}^{l_{2}} \times \mathbf{R}^{2})}^{2} \| v \|_{L_{t}^{6} L_{x}^{3}(G_{\beta}^{l_{2}} \times \mathbf{R}^{2})}^{2} \lesssim 2^{2i - l_{2}} \| v_{0} \|_{L^{2}}^{2} \| u \|_{\tilde{X}_{i}(G_{\alpha}^{i} \times \mathbf{R}^{2})}^{2}.
\end{equation}

\noindent Making a similar computation for $(\ref{5.43.5})$ implies

\begin{equation}
2^{l_{2} - 2i} (\ref{5.43}) \lesssim \| v_{0} \|_{L^{2}}^{2} \| u \|_{\tilde{X}_{i}(G_{\alpha}^{i} \times \mathbf{R}^{2})}^{2}.
\end{equation}

\noindent Next take $(\ref{5.43.1})$. By the product rule and

\begin{equation}
\frac{\partial}{\partial t} \phi(2^{-l_{2}}(\eta_{1} - \xi(t))) = 2^{-l_{2}} (\nabla \phi)(2^{-l_{2}}(\eta_{1} - \xi(t))) \cdot \xi'(t),
\end{equation}

\noindent so by the Sobolev embedding theorem,

\begin{equation}
2^{l_{2} - 2i} (\ref{5.43.1}) \lesssim 2^{-i} \| v_{0} \|_{L^{2}}^{2} \| \frac{1}{|\nabla|} u_{h} \|_{L_{t}^{\infty} L_{x}^{2}}^{2} \| u_{l} \|_{L_{t,x}^{\infty}}^{2} (\int_{G_{\beta}^{l_{2}}} |\xi'(t)| dt) \lesssim 2^{l_{2} - i} \| v_{0} \|_{L^{2}}^{2}.
\end{equation}

\noindent Finally take $(\ref{5.43.2})$. $e^{-it |\eta_{1}|^{2}} \bar{u}_{h}(t, \eta_{1})$ is the Fourier transform of $e^{it \Delta} \bar{u}_{h}$, and by direct computation

\begin{equation}
\frac{\partial}{\partial t} (e^{it \Delta} \bar{u}_{h}) = e^{it \Delta} P_{h} F(u) + (\frac{d}{dt} P_{\xi(t), \geq l_{2} - 5}) u.
\end{equation}

\noindent Likewise, $e^{it |\eta_{3}|^{2}} \hat{u}_{l}(t, \eta_{3})$ is the Fourier transform of $e^{-it \Delta} u_{l}$, and

\begin{equation}
\frac{\partial}{\partial t} (e^{-it \Delta} u_{l}) = e^{-it \Delta} P_{l} F(u) + (\frac{d}{dt} P_{\xi(t), \leq l_{2} - 5}) u.
\end{equation}

\noindent Therefore,

\begin{equation}\label{5.44}
\aligned
|\int_{a}^{b} \int \int \frac{(x - y)}{|x - y|} \cdot Im[\bar{v} (\nabla - i \xi(t)) v] (t,x) \\ \times \int \int \int \int \phi(\frac{(\eta_{1} + \eta_{2} + \eta_{3} + \eta_{4})}{2^{l_{2} - 10}})  e^{iy \cdot (\eta_{1} + \eta_{2} + \eta_{3} + \eta_{4})}  \frac{e^{it q( \eta)}}{q(\eta)} (1 - \phi(2^{-l_{2}} (\eta_{1} - \xi(t))) \phi(2^{-l_{2}} (\eta_{2} - \xi(t)))) \\ \frac{\partial}{\partial t}[ e^{-it |\eta_{1}|^{2}} \hat{\bar{u}}_{h}(t, \eta_{1}) e^{-it |\eta_{2}|^{2}} \hat{\bar{u}}_{h}(t, \eta_{2}) e^{it |\eta_{3}|^{2}} \hat{u}_{l}(t, \eta_{3}) e^{it |\eta_{4}|^{2}} \hat{u}_{l}(t, \eta_{4})] d\eta_{1} d\eta_{2} d\eta_{3} d\eta_{4} dx dy dt|
\endaligned
\end{equation}

\begin{equation}\label{5.45}
\aligned
\lesssim 2^{i - l_{2}} \| v_{0} \|_{L^{2}}^{2} \| \frac{1}{|\nabla|} u_{h} \|_{L_{t}^{\infty} L_{x}^{2}}^{2} \| P_{\xi(t), \leq l_{2}} u \|_{L_{t,x}^{\infty}}^{2} (\int_{G_{\beta}^{l_{2}}} |\xi'(t)| dt) \\
+ 2^{i} \| v_{0} \|_{L^{2}}^{2} \| \frac{1}{|\nabla|} u_{h} \|_{L_{t,x}^{4}}^{2} \| P_{l} (u_{h}^{3}) \|_{L_{t,x}^{2}} \| u_{l} \|_{L_{t,x}^{\infty}} + 2^{i} \| v_{0} \|_{L^{2}}^{2} \| \frac{1}{|\nabla|} u_{h} \|_{L_{t}^{6} L_{x}^{3}}^{2} \| P_{l} (u_{l}^{3}) \|_{L_{t}^{2} L_{x}^{4}} \| u_{l} \|_{L_{t}^{6} L_{x}^{12}} \\
+ 2^{i} \| v_{0} \|_{L^{2}}^{2} \| \frac{1}{|\nabla|} u_{h} \|_{L_{t}^{3} L_{x}^{6}} \| \frac{1}{|\nabla|} P_{\xi(t), \geq l_{2} - 5} u \|_{L_{t}^{3} L_{x}^{6}} \| u_{l} \|_{L_{t}^{6} L_{x}^{12}}^{4} + 2^{i} \| v_{0} \|_{L^{2}}^{2} \| \frac{1}{\Delta} u_{h} \|_{L_{t,x}^{4}} \| u_{h} \|_{L_{t,x}^{4}}^{3} \| u_{l} \|_{L_{t,x}^{\infty}}^{2}.
\endaligned
\end{equation}

\noindent By conservation of mass and the Sobolev embedding theorem,

\begin{equation}\label{5.45.1}
\lesssim 2^{i - l_{2}} \| v_{0} \|_{L^{2}}^{2} \| \frac{1}{|\nabla|} u_{h} \|_{L_{t}^{\infty} L_{x}^{2}}^{2} \| P_{\xi(t), \leq l_{2}} u \|_{L_{t,x}^{\infty}}^{2} (\int_{G_{\beta}^{l_{2}}} |\xi'(t)| dt) \lesssim 2^{i} \| v_{0} \|_{L^{2}}^{2}.
\end{equation}

\noindent Next, by $(\ref{2.38})$ and the Sobolev embedding theorem,

\begin{equation}\label{5.45.2}
2^{i} \| v_{0} \|_{L^{2}}^{2} \| \frac{1}{|\nabla|} u_{h} \|_{L_{t,x}^{4}}^{2} \| P_{l} (u_{h}^{3}) \|_{L_{t,x}^{2}} \| u_{l} \|_{L_{t,x}^{\infty}} + 2^{i} \| v_{0} \|_{L^{2}}^{2} \| \frac{1}{\Delta} u_{h} \|_{L_{t,x}^{4}} \| u_{h} \|_{L_{t,x}^{4}}^{3} \| u_{l} \|_{L_{t,x}^{\infty}}^{2} \lesssim 2^{i} \| v_{0} \|_{L^{2}}^{2} \| u \|_{\tilde{X}_{i}(G_{\alpha}^{i} \times \mathbf{R}^{2})}^{4}.
\end{equation}

\noindent Finally by $(\ref{2.37})$, $(\ref{2.38})$, and the Sobolev embedding theorem,

\begin{equation}\label{5.46}
2^{i} \| v_{0} \|_{L^{2}}^{2} \| \frac{1}{|\nabla|} u_{h} \|_{L_{t}^{6} L_{x}^{3}}^{2} \| P_{l} (u_{l}^{3}) \|_{L_{t}^{2} L_{x}^{4}} \| u_{l} \|_{L_{t}^{6} L_{x}^{12}}
+ 2^{i} \| v_{0} \|_{L^{2}}^{2} \| \frac{1}{|\nabla|} u_{h} \|_{L_{t}^{3} L_{x}^{6}} \| \frac{1}{|\nabla|} P_{\xi(t), \geq l_{2} - 5} u \|_{L_{t}^{3} L_{x}^{6}} \| u_{l} \|_{L_{t}^{6} L_{x}^{12}}^{4}
\end{equation}

\begin{equation}\label{5.47}
\lesssim 2^{i - 2 l_{2}} \| v_{0} \|_{L^{2}}^{2} \| u \|_{\tilde{X}_{i}(G_{\alpha}^{i} \times \mathbf{R}^{2})}^{2} \| u_{l} \|_{L_{t}^{6} L_{x}^{12}(G_{\beta}^{l_{2}} \times \mathbf{R}^{2})}^{4}.
\end{equation}

\noindent By $(\ref{2.38})$, conservation of mass, and the Sobolev embedding theorem

\begin{equation}\label{5.48}
\| u_{l} \|_{L_{t}^{6} L_{x}^{12}(G_{\beta}^{l_{2}} \times \mathbf{R}^{2})}^{2} \lesssim \sum_{0 \leq j_{1} \leq l_{2} - 5} 2^{j_{1}/2} \| P_{\xi(t), j_{1}} u \|_{L_{t}^{3} L_{x}^{6}(G_{\beta}^{l_{2}} \times \mathbf{R}^{2})} \lesssim 2^{l_{2}/2} \| u \|_{\tilde{X}_{l_{2}}(G_{\beta}^{l_{2}} \times \mathbf{R}^{2})}.
\end{equation}

\noindent Therefore, since $l_{2} \leq i$, $G_{\beta}^{l_{2}} \subset G_{\alpha}^{i}$,

\begin{equation}\label{5.49}
(\ref{5.47}) \lesssim 2^{i} \| v_{0} \|_{L^{2}}^{2} \| u \|_{\tilde{X}_{i}(G_{\alpha}^{i} \times \mathbf{R}^{2})}^{4}.
\end{equation}

\noindent Collecting $(\ref{5.45.1})$, $(\ref{5.45.2})$, and $(\ref{5.49})$,

\begin{equation}\label{5.50}
2^{l_{2} - 2i} (\ref{5.43.2}) \lesssim 2^{l_{2} - i} \| v_{0} \|_{L^{2}}^{2} (1 + \| u \|_{\tilde{X}_{i}(G_{\alpha}^{i} \times \mathbf{R}^{2})}^{4}).
\end{equation}

\noindent This finally completes the proof of lemma $\ref{appendix}$. $\Box$

\end{document}